\documentclass[a4]{amsart}

\input xypic
\xyoption{all}
\usepackage{tikz-cd}
\usepackage{amssymb} 
\usepackage{enumitem} 
\usepackage{scalerel} 

\usepackage{appendix}
\newtheorem{theorem}{Theorem}[section]
\newtheorem{proposition}[theorem]{Proposition}
\newtheorem{lemma}[theorem]{Lemma}
\newtheorem{corollary}[theorem]{Corollary}

\theoremstyle{definition}
\newtheorem{definition}[theorem]{Definition}
\newtheorem{remark}[theorem]{Remark}
\newtheorem{example}[theorem]{Example}

\newcommand{\tinyI}{\ensuremath{\scriptscriptstyle I}} 
\newcommand{\scaleI}{\ensuremath{\scaleto{I}{3pt}}} 
\newcommand{\sscaleI}{\ensuremath{\scaleto{I}{2.2pt}}} 

\newcommand{\cohlgy}[1]{\ensuremath{H^{*}(#1)}} 

\newcounter{bean}

\newcommand{\namedright}[3]{\ensuremath{#1\stackrel{#2}
 {\longrightarrow}#3}}
\newcommand{\nameddright}[5]{\ensuremath{#1\stackrel{#2}
 {\longrightarrow}#3\stackrel{#4}{\longrightarrow}#5}}

\newcommand{\larrow}{\relbar\!\!\relbar\!\!\rightarrow} 
\newcommand{\llarrow}{\relbar\!\!\relbar\!\!\larrow} 
\newcommand{\lllarrow}{\relbar\!\!\relbar\!\!\llarrow} 
\newcommand{\llllarrow}{\relbar\!\!\relbar\!\!\lllarrow} 

\newcommand{\lnamedright}[3]{\ensuremath{#1\stackrel{#2}
 {\larrow}#3}}
\newcommand{\lnameddright}[5]{\ensuremath{#1\stackrel{#2}
 {\larrow}#3\stackrel{#4}{\larrow}#5}} 
 
\newcommand{\llnamedright}[3]{\ensuremath{#1\stackrel{#2}
 {\llarrow}#3}}
\newcommand{\llnameddright}[5]{\ensuremath{#1\stackrel{#2}
 {\llarrow}#3\stackrel{#4}{\llarrow}#5}} 
 
\newcommand{\lllnamedright}[3]{\ensuremath{#1\stackrel{#2}
 {\lllarrow}#3}}
\newcommand{\lllnameddright}[5]{\ensuremath{#1\stackrel{#2}
 {\lllarrow}#3\stackrel{#4}{\lllarrow}#5}}

\newcommand{\qqed}{\hfill\Box}

\newcommand{\Z}{\mathbb{Z}} 
\newcommand{\N}{\mathbb{N}} 
\newcommand{\Q}{\mathbb{Q}} 
\newcommand{\C}{\mathbb{C}} 

\begin{document} 


\title{Weighted polyhedral products and Steenrod's problem} 
\author{Tseleung So}
\address{School of Mathematical Sciences, University of Southampton, Southampton 
               SO17 1BJ, United Kingdom}
\email{larry.so.tl@gmail.com} 
\author{Donald Stanley}
\address{Department of Mathematics and Statistics, University of Regina, Regina, SK S4S 0A2, Canada}
\email{Donald.Stanley@uregina.ca}
\author{Stephen Theriault}
\address{School of Mathematical Sciences, University of Southampton, Southampton 
               SO17 1BJ, United Kingdom}
\email{S.D.Theriault@soton.ac.uk} 



\begin{abstract} 
We construct a weighted version of polyhedral products and compute its cohomology in special cases. This is applied to resolve Steenrod's cohomology realization problem in a case related to products of spheres. 
\end{abstract}

\maketitle

\section{Introduction} 
\label{sec:intro} 

Steenrod~\cite{S} asked which graded algebras could be realized as the cohomology 
of spaces. Concentrating on the polynomial case, he showed that a graded polynomial algebra $\mathbb{Z}[x]$ can be realized in 
this way if and only if the degree of $x$ is $2$ or $4$, these being the cohomologies of 
the classifying spaces of $S^{1}$ and $S^{3}$ respectively. Other immediate examples 
of polynomial algebras with more than one generator that can be realized are the cohomologies 
of the classifying spaces of the classical matrix groups $U(n)$ and $Sp(n)$ for $n\geq 2$. 
The intimate connection between this realization problem and the existence of finite 
loop spaces was the germ of a grand program of research that played an important 
developmental role in homotopy theory for thirty years. Allowing the ground ring to be 
the integers modulo $p$ and using the machinery of $p$-compact groups, the 
program culminated with the complete resolution of the polynomial version of Steenrod's problem by Notbohm~\cite{N} 
for odd primes and Andersen-Grodal~\cite{AG} for $p=2$ (and other more general rings). 

The more general version of Steenrod's original problem allows for tensor products of 
polynomial algebras and exterior algebras, and quotients of them by appropriate relations. 
When $R=\mathbb{Z}[x_{1},\ldots,x_{m}]/I$, each generator~$x_{i}$ has degree $2$, 
and $I$ is an ideal generated by squarefree monomials, so that $R$ is a graded Stanley-Reisner ring, Davis-Januszkiewicz \cite{DJ} showed that these algebras can 
be realized as the cohomology of certain polyhedral products. 
Trevisan~\cite{T} showed that $R$ can also be realized if we allow $I$ to be generated by monomials that are not necessarily square free.
Bahri, Bendersky, Cohen 
and Gitler~\cite{BBCG3} shortly after gave a different proof. The constructions of  Davis-Januszkiewicz and Trevisan also work when the generators are of degrees both $2$ and $4$. If we are willing to mod out by torsion 
 the first two authors~\cite{SS} recently showed that spaces can be constructed realizing 
 $R=\mathbb{Z}[x_{1},\ldots,x_{m}]/I$ where $I$ is monomial and the $x_i$ have arbitrary even degree.
 When not moding out by torsion and allowing arbitrary degrees the problem is considerably more difficult, but  Takeda \cite{Tak} has classified which square 
free monomial algebras are realizable, under the condition that  
generators $a,b$ with $|a|=|b|=2^i$ for some $i>1$ satisfy $ab=0$. Takeda and the second author have also related the Stanley-Reisner case to a graph colouring problem \cite{ST}.~\smallskip

\noindent 
{\bf Steenrod's problem for sphere product algebras.}
Fix a positive integer $m$ 
and let $[m]=\{1,\ldots,m\}$. Associate to each $i\in [m]$ a degree~$d_{i}$. For a subring $R\subset \mathbb{Q}$, we call an $R$-algebra that is isomorphic to the cohomology of a product of spheres 
$H^*(\prod_{i \in [m]} S^{d_i},R)\cong \otimes_{i\in [m]} R[x_{n_i}]/(x_{d_i}^2)$ an \emph{$R$-sphere product algebra}. 
We consider a special case that gives sphere product algebras a ``weighted" multiplication.  Call 
$\mathfrak{c}=\{\mathfrak{c}_{\sigma}\mid\sigma\subseteq [m]\}$ a coeffcient sequence (Definition \ref{coeff_sequ}). Let $A(\mathfrak{c})$ be the graded commutative algebra defined 
as a module by 
\[A(\mathfrak{c})=\bigoplus_{\sigma\subseteq [m]} R\langle a_{\sigma}\rangle\] 
where $a_{\sigma}$ has degree $d_{\sigma}=\sum_{i\in\sigma} d_{i}$, with the (weighted) multiplication 
determined by the formula 
\[\prod_{a\in\sigma} a_{\{i\}}=\mathfrak{c}_{\sigma}a_{\sigma}.\] 
For the multiplication to make sense the integers $\mathfrak{c}_{\sigma}$ must satisfy certain conditions. 
At the very least this requires that if $\sigma$ is the disjoint union of $\sigma_{1}$ and $\sigma_{2}$ 
then $\mathfrak{c}_{\sigma_{1}} \mathfrak{c}_{\sigma_{2}}\mid \mathfrak{c}_{\sigma}$. 
As normalizing conditions, we also assume 
that $\mathfrak{c}_{\emptyset}=1$ and $\mathfrak{c}_{\{i\}}=1$ for $1\leq i\leq m$, 
and $\mathfrak{c}_{\sigma}$ does not have $u$ 
as a factor for any non-identity unit $u$ in $R$. 
Algebraically all coefficient sequences $\mathfrak{c}=\{\mathfrak{c}_{\sigma}\mid\sigma\subseteq [m]\}$ 
give rise to valid weighted multiplications and, up to algebra isomorphism, all weighted multiplications come from  normalized coefficient sequences (Lemma \ref{all_coeff}). 
A topological question is to determine which of the weighted sphere product 
algebras~$A(\mathfrak{c})$ can be realized as the cohomology of a space. 

\smallskip  

\noindent 
{\bf Weighted polyhedral products}. 
The most novel and important part of the paper is the construction of a new family of spaces, 
called weighted polyhedral products, in order to geometrically realize the algebras $A(\mathfrak{c})$ for a 
large family $\mathcal{F}$ of coefficient sequences $\mathfrak{c}$. To describe $\mathcal{F}$, we introduce an 
intermediate object called a power sequence (Definition \ref{dfn_power_sequence}) which to each $\sigma\subseteq [m]$ associates a sequence 
$c^{\sigma}=(c^{\sigma}_{1},\ldots,c^{\sigma}_{m})$. These sequences have the property 
that if $\tau\subseteq\sigma$ then $c_{i}^{\tau}$ divides $c_{i}^{\sigma}$ for each $1\leq i\leq m$. 
Let $\mathcal{F}$ be the family of coefficient sequences
$\mathfrak{c}=\{\mathfrak{c}_{\sigma}\}$ with the property that for some power sequence $c^{\sigma}$,
$\mathfrak{c}_{\sigma}=\prod_{i\in\sigma} c_{i}^{\sigma}$ for every $\sigma\subseteq [m]$. 

Our construction generalizes the notion of polyhedral products, which
have been the subject of intense recent attention. They are 
a family of spaces that have wide application across a range of mathematical disciplines, 
notably toric topology and geometric group theory  (see the survey~\cite{BBC} for more details). 
Explicitly, let $K$ be a simplicial complex on the vertex set $[m]$. To each vertex in $K$ 
associate a pair of pointed spaces $(X_{i},A_{i})$, where $A_{i}$ is a pointed 
subspace of $X_{i}$. For each face $\sigma$ in~$K$, let 
\[
(\underline{X},\underline{A})^{\sigma}=\{(x_1,\ldots,x_m)\in\prod^m_{i=1}X_i|x_i\in A_i\text{ whenever }i\notin\sigma\}.
\]
The \emph{polyhedral product} $(\underline{X},\underline{A})^{K}$ is defined to be
\[
(\underline{X},\underline{A})^{K}=\bigcup_{\sigma\subseteq K}(\underline{X},\underline{A})^{\sigma}\subseteq\prod^m_{i=1}X_i.
\]
Note that $(\underline{X},\underline{A})^{K}$ is a subspace of the product $\prod_{i=1}^{m} X_{i}$. 

To define a weighted polyhedral product we start with pairs that have compatible 
``power maps". This occurs, for example, if $X_{i}$ and $A_{i}$ have compatible $H$-structures, 
or if they have compatible co-$H$-structures. Then for a power sequence
$c^{\sigma}=(c_{1}^{\sigma},\ldots,c_{m}^{\sigma})$ we associate to each
$c_{i}^{\sigma}$ the power map of degree $c_{i}^{\sigma}$ on the pair $(X_{i},A_{i})$. 
Regarding~$K$ as filtered by its full $k$-skeletons, $K_{0}\subseteq K_{1}\cdots\subseteq K_{m}$, 
the weighted polyhedral product $(\underline{X},\underline{A})^{K,c}$ can be thought of as 
starting with the usual polyhedral product and using these power maps to introduce weightings (degree maps) for the 
``attaching maps" in passing from $(\underline{X},\underline{A})^{K_{i}}$ to 
$(\underline{X},\underline{A})^{K_{i+1}}$. The precise definition is lengthy and given in 
Section~\ref{sec:construction}. In particular, for these pairs the usual polyhedral 
product is the special case of the weighted polyhedral product corresponding to the weighting 
$c^{\sigma}=(1,\ldots,1)$ for all $\sigma\in K$.    

As weighted polyhedral products generalize polyhedral products, it is interesting to consider 
whether they satisfy analogous properties. We first show that polyhedral products can be recovered 
as a special case of weighted polyhedral products  for a family of power sequences 
satisfying a minimality condition. 
We then identify a relation to higher Whitehead products (Section~\ref{sec:properties}), 
prove functoriality with 
respect to maps of power couples (Section~\ref{sec:functoriality1}) and functoriality 
with respect to projections to full subcomplexes (Section~\ref{sec:functoriality2}), establish 
a homotopy decomposition of $\Sigma(\underline{X},\underline{\ast})^{K,c}$ 
as a wedge of summands (Section \ref{sec:suspsplitting}) that is analogous to the homotopy decomposition 
of the suspension of the usual polyhedral product in~\cite{BBCG1}, and calculate 
the cohomology of $(\underline{X},\underline{\ast})^{K,c}$ (Section \ref{coho_compute}). 
\smallskip 

\noindent 
\textbf{Steenrod's problem and weighted polyhedral products.} 
The cohomological description of $(\underline{X},\underline{\ast})^{K,c}$ lets us return 
to Steenrod's problem and the geometric realization of $A(\mathfrak{c})$. 

\begin{theorem} 
   \label{geomrealization} (Theorem \ref{Steenrodpowerseq})
   Let $\mathfrak{c}=\Phi(c)\in\mathcal{F}$ (see~\eqref{eqn_Phi map}) and let $A(\mathfrak{c})$ be defined as above. Then 
   \[A(\mathfrak{c})\cong\cohlgy{(\underline{X},\underline{\ast})^{\Delta^{m-1},c}}\] 
   where $X_{i}=S^{d_{i}}$. 
\end{theorem} 

Theorem~\ref{geomrealization} solves Steenrod's problem in the case of weighted sphere product algebras 
$A(\mathfrak{c})$ when~$\mathfrak{c}$ is a coefficient sequence in the family $\mathcal{F}$. 
In fact, we prove a more general result by calculating the cohomology ring  
$H^*((\underline{X},\underline{\ast})^{K,c})$ for any simplicial complex $K$ and power sequence $c$,
where each~$X_i$ is a suspension with  torsion-free homology (Theorem~\ref{weighted_exterior_realization}). As graded modules,
$H^*((\underline{X},\underline{\ast})^{K,c})$ is the tensor product of $H^*(X_i)$'s quotiented by the generalized Stanley-Reisner ideal, resembling the polyhedral product case. However, the multiplication is more complicated and is ``twisted'' by coefficients that are determined by $c$.
Theorem~\ref{geomrealization} is the special case when $X_i=S^{d_i}$ and $K=\Delta^{m-1}$.

The family $\mathcal{F}$ can be compared to the family of all coefficient sequences. 
We show in Lemma~\ref{no_surj} that there are more coefficient sequences than those in $\mathcal{F}$, 
but the two families are closely linked in the sense that the coefficient sequences form a commutative 
monoid and $\mathcal{F}$ generates the same group completion as is generated by all coefficient sequences 
(Corollary~\ref{completion}). 

One future direction is to determine whether Steenrod's problem can be resolved for $A(\mathfrak{c})$ 
for all coefficient sequences. In the two generator case all possible $A(\mathfrak{c})$ (Remark \ref{use2}) can be realized by our construction, but in the three generator case our construction does not give all possible $A(\mathfrak{c})$ (Lemma~\ref{no_surj}).
However, using other methods the authors along with Williams have constructed spaces that realize all weighted sphere product algebras \cite{SSTW} when $m=3$. 

We also look at another interpretation of our results. If $B$ is a finitely generated projective 
$R$-algebra and $A\subset B$ is an $R$-subalgebra, then $A$ is an \emph{order} in $B$ if 
$A\otimes_R \mathbb{Q}=B$. In our case, the weighted sphere product algebra $A(\mathfrak{c})$ 
is an order in a sphere product algebra and we can ask of all such orders are weighted sphere product algebras. It turns out that when $m=3$ all orders in a sphere product 
algebra $B$ are weighted sphere product algebras if and only if $B$ is graded 
exterior~\cite{SSTW}.~\smallskip  

Another future direction is to further develop properties of weighted polyhedral products. Important 
properties of polyhedral products have been established and work continues at great pace to reveal more. 
It would be interesting to determine analogues of these properties in the more general case of 
weighted polyhedral products. Possibilities include determining the homotopy type of 
$\Sigma(\underline{X},\underline{A})^{K,c}$ and the cohomology of $(\underline{X},\underline{A})^{K,c}$ 
for general pairs $(X,A)$, extending known results for the usual polyhedral products and the work in this paper for pairs $(X,\ast)$; determining the 
homotopy type of $(CA,A)^{K,c}$ for special families of simplicial complexes $K$, such as the weighted versions of moment-angle complexes $(D^2,S^1)^{K,c}$ and real moment-angle complexes $(D^1,S^0)^{K,c}$; and determining 
how the weighted polyhedral products for pairs $(S^{1},\ast)$ and $(\mathbb{R}P^{\infty},\ast)$ 
inform on right-angled Artin and right-angled Coxeter groups, as do the polyhedral products 
for those pairs.  
As seen in Section \ref{relWH} higher order Whitehead products can be constructed using weighted polyhedral products which also give universal spaces for more general higher order operations. Similar constructions could be done for Toda brackets and Massey products. 

\medskip 
\noindent 
\textbf{Glossary  of notation}. 
It may be helpful to list the main notational conventions used in the paper. In particular, 
we use the Buchstaber-Panov convention of denoting faces of simplicial complexes by 
lower case Greek letters and subsets of $\{1,\ldots,m\}$ by capital Roman letters. 
\medskip 

\begin{itemize} 
\item $[m]=\{1,\ldots,m\}$, an underlying set;  
\item $\Delta^{m-1}$, the $m$-simplex; 
\item $K,L$, simplicial complexes with vertex set $[m]$; 
\item $\sigma,\tau\in K$, faces of the simplicial complex; 
\item $\partial\sigma$, the boundary of the face $\sigma$; 
\item $I,J\subseteq [m]$, subsets of $[m]$; 
\item $c$, a power sequence: Definition~\ref{dfn_power_sequence}; 
\item $\mathfrak{c}$, a coefficient sequence: Definition~\ref{coeffseq}; 
\item $(\underline{X},\underline{A})^{\sigma}$, a polyhedral product for a simplex, introduced 
         in Definition~\ref{dfn_poly prod}; 
\item $(\underline{X},\underline{A})^{K}$, a polyhedral product for a simplicial complex:  
         Definition~\ref{dfn_poly prod}; 
\item $(\underline{X},\underline{A})^{\sigma,c}$, a weighted polyhedral product for a simplex: 
         Definition~\ref{dfn_weighted poly prod system}; 
\item $(\underline{X},\underline{A})^{K,c}$, a weighted polyhedral product for a simplicial complex: 
         Definition~\ref{dfn_weighted poly prod}; 
\item $\eta(\sigma),\eta(\partial\sigma),\eta(K)$, a map from the polyhedral product for $\sigma,\partial\sigma,K$
         respectively to the associated weighted polyhedral product for $\sigma,\partial\sigma,K$: 
         Definitions~\ref{dfn_weighted poly prod system} and~\ref{dfn_weighted poly prod}; 
\item $\underline{c}^{\sigma/\tau}$, a product of power maps on $(\underline{X},\underline{A})^{K}$ for 
         $\tau\subseteq\sigma$: Definition~\ref{dfn_c^sigma_tau}; 
\item $\imath_{\sigma}^{\tau}, \imath_{\partial\sigma}^{\sigma}$, the inclusion of the polyhedral product for 
         $\sigma, \partial\sigma$ respectively into the polyhedral product for $\tau,\sigma$, where 
         $\sigma\subseteq\tau$: Definitions~\ref{dfn_weighted poly prod system} and \ref{dfn_poly prod};
\item $\imath_{\sigma,c}^{\tau,c}, \imath_{\partial\sigma,c}^{\sigma,c}$, the inclusion of the weighted polyhedral product for 
         $\sigma, \partial\sigma$ respectively into the weighted polyhedral product for $\tau,\sigma$, where 
         $\sigma\subseteq\tau$: Definition~\ref{dfn_weighted poly prod system};      
\item $K_{I}$, the full subcomplex of $K$ on the vertex set $I$: start of Section~\ref{sec:functoriality2};   
\item $\overline{K}_{I}$, the full subcomplex of $K$ on the vertex set $I$ but now 
         regarded as being on the vertex set $[m]$: start of Section~\ref{sec:functoriality2};  
\item $\iota_{L}^{K}$, the inclusion of the polyhedral product for $L$ into the polyhedral product 
         for $K$: start of Section~\ref{sec:functoriality2}; 
\item $\iota_{L,c}^{K,c}$, the inclusion of the weighted polyhedral product for $L$ into the weighted 
         polyhedral product for $K$; Proposition~\ref{prop-projectLK};
\item $p_{\sigma}^{\tau}, p_{K}^{K_{I}}$, the projection of the polyhedral product for 
         $\sigma,K$ respectively to the polyhedral product for $\tau,K_{I}$, where 
         $\tau\subseteq\sigma$: start of Section~\ref{sec:functoriality2}; 
\item $p_{\sigma,c}^{\tau,c}, p_{K,c}^{K_{I},c}$, the projection of the weighted polyhedral product for 
         $\sigma,K$ respectively to the weighted polyhedral product for $\tau,K_{I}$, where 
         $\tau\subseteq\sigma$: Lemma~\ref{3diagramlemma} and 
         Proposition~\ref{prop-projectLK} respectively;        
\item $\underline{X}^{\wedge\sigma}$, the smash product $X_{i_{1}}\wedge\cdots\wedge X_{i_{k}}$, 
         where $\sigma=(i_{1},\ldots,i_{k})$: start of Section~\ref{sec:suspsplitting}; 
\item $q_{\sigma}$, the quotient map from $(\underline{X},\underline{\ast})^{K}$ to 
         $\underline{X}^{\wedge\sigma}$: start of Section~\ref{sec:suspsplitting}; 
\item $q_{\sigma,c}$, the quotient map from $(\underline{X},\underline{\ast})^{K,c}$   
         to $X^{\wedge\sigma}$: Equation~(\ref{diagram_A=* pushout}); 
\item $\Lambda(\underline{Y},c)$, a weighted algebra: Definition~\ref{def_lambda_weighted_graded}; 
\item $A(\mathfrak{c})$, a more specialized weighted sphere product algebra: Definition~\ref{Acdef}. 
\end{itemize}

\section*{Acknowledgment}
The first author was supported by Pacific Institute for the Mathematical Sciences (PIMS) Postdoctoral Fellowship, NSERC Discovery Grant and NSERC RGPIN-2020-06428. The second author is supported by NSERC RGPIN-05466-2020.

\section{Construction of weighted polyhedral products} 
\label{sec:construction} 

We work with compactly generated spaces, products will also be taken in that category. For us a (relative) CW complex is built by adding cells along not necessarily cellular maps. A sequence of spaces and maps $A\stackrel{f}{\rightarrow} B\rightarrow C$ is a cofibration sequence if $f$ is a cofibration and $C$ is homeomorphic to the quotient $B/A$. We refer to $C$ as the cofiber of $f$.
From now on let $(X,A)$ denote a pair of pointed spaces such that $A$ is a subspace of~$X$ 
and the inclusion $A\rightarrow X$ is a pointed map. A map of pairs $f:(X,A)\longrightarrow(Y,B)$ is a pointed 
map $f:X\longrightarrow Y$ such that $f(A)\subset B$. Equivalently, there is a commutative diagram
\[
\xymatrix{
A\ar[r]^-{f|_A}\ar[d]	&B\ar[d]\\
X\ar[r]^-{f}					&Y.
}
\] 
The formal definition of a polyhedral product is as follows. 

\begin{definition}\label{dfn_poly prod} 
Let $(\underline{X},\underline{A})=\{(X_i,A_i)\}^m_{i=1}$ be a sequence of $m$ pairs of spaces and 
let $K$ be a simplicial complex on $[m]=\{1,\ldots,m\}$. For each $\sigma\in K$, let $(\underline{X},\underline{A})^{\sigma}$ be the subspace of $\prod^m_{i=1}X_i$ defined by
\[
(\underline{X},\underline{A})^{\sigma}=\{(x_1,\ldots,x_m)\in\prod^m_{i=1}X_i\mid x_i\in A_i\text{ whenever }i\notin\sigma\}.
\]
The \emph{polyhedral product} $(\underline{X},\underline{A})^{K}$ is defined as
\[
(\underline{X},\underline{A})^{K}=\bigcup_{\sigma\subseteq K}(\underline{X},\underline{A})^{\sigma}\subseteq\prod^m_{i=1}X_i.
\]

For $\tau\subset \sigma$ and $L\subset K$, we use the notation $\imath_\tau^\sigma\colon (\underline{X},\underline{A})^{\tau}\rightarrow (\underline{X},\underline{A})^{\sigma}$
and $\imath_L^K:(\underline{X},\underline{A})^L\rightarrow (\underline{X},\underline{A})^K$ for the inclusions. More generally if $L$ is a simplicial complex on $I\subset [m]$ all of whose simplices are simplices in $K$ we also use $\imath_L^K:(\underline{X},\underline{A})^L\rightarrow (\underline{X},\underline{A})^K$ to denote the inclusion. 
\end{definition}

\begin{definition}\label{dfn_poly prod functoriality}
Let $(\underline{X},\underline{A})=\{(X_i,A_i)\}^m_{i=1}$ and $(\underline{Y},\underline{B})=\{(Y_i,B_i)\}^m_{i=1}$ be sequences of $m$ pairs of spaces and 
let $\underline{f}=\{f_i:(X_i,A_i)\longrightarrow(Y_i,B_i)\}^m_{i=1}$ be a sequence of maps of pairs. The \emph{induced map} 
$\underline{f}^{K}:(\underline{X},\underline{A})^{K}\longrightarrow(\underline{Y},\underline{B})^{K}$ is given by the restriction 
of $\prod^m_{i=1}f_i:\prod^m_{i=1}X_i\longrightarrow\prod^m_{i=1}Y_i$ to $(\underline{X},\underline{A})^K$.
\end{definition}

Alternatively there is a categorical construction of a polyhedral product $(\underline{X},\underline{A})^K$.  We also let $K$ denote the face category of $K$ 
whose objects are the simplices of $K$ and morphisms are their inclusions, and let $TOP_*$ be the category of pointed topological spaces. 
Given a sequence of pairs of spaces $(\underline{X},\underline{A})$, define the functor $(\underline{X},\underline{A})^{\bullet}:K\longrightarrow TOP_*$ by 
sending $\sigma\in K$ to $(\underline{X},\underline{A})^{\sigma}$ and 
sending $\tau\rightarrow\sigma$ in $K$ to the inclusion $(\underline{X},\underline{A})^{\tau}\longrightarrow (\underline{X},\underline{A})^{\sigma}$.

Note that $(\underline{X},\underline{A})^{K}=\mbox{colim}_{\sigma\subseteq K}(\underline{X},\underline{A})^{\sigma}$ and that the $\imath_L^K$ are the induced map between the colimits.

\begin{definition}\label{dfn_N action}  
A \emph{power couple} is a pair of spaces $(X,A)$ equipped with a collection of maps of pairs $\{\rho_a:(X,A)\to(X,A)\}_{a\in\N}$, called \emph{power maps}, such that $\rho_1$ is the identity map and $\rho_a\circ\rho_b=\rho_{ab}$ for any $a,b\in\N$. We say that $f:(X,A)\to(Y,B)$ is a \emph{map between power couples} if it is a map of pairs $f:(X,A)\to(Y,B)$ making the diagrams
\[ 
\xymatrix{
X\ar[r]^-{f}\ar[d]_-{\rho_a}	&Y\ar[d]^-{\widetilde{\rho}_a}\\
X\ar[r]^-{f}					&Y
} 
\qquad 
\xymatrix{
A\ar[r]^-{f|_A}\ar[d]_-{\rho_a|_A}	&B\ar[d]^-{\widetilde{\rho}_a|_B}\\
A\ar[r]^-{f|_A}						&B
}
\]
commute for each $a\in\N$, where $\rho_a$ and $\widetilde{\rho}_a$ are power maps for $(X,A)$ and $(Y,B)$.
\end{definition}

\begin{example}\label{ex_monoid} 
Let $(X,A)$ be a pair of spaces where $X$ and $A$ are topological monoids with their identities being the basepoints and the inclusion $A\to X$ being multiplicative. Then $(X,A)$ is a power couple where the power map $\rho_a$ is the $a$-fold self-multiplication
\[
\rho_a(x)=\underbrace{x\cdot\ldots\cdot x}_{a\text{ times}}.
\]
\end{example}

\begin{example}\label{ex_suspension}
Identify $S^1$ with the unit circle $\{e^{ti}|t\in\mathbb{R}\}$ in $\mathbb{C}$. For $a\in\N$ define the power map on $(S^1,\ast)$ to be $e^{ti}\mapsto e^{ati}$. Let $(X,A)$ be a pair of spaces. Since $\Sigma X\cong S^1\wedge X$, the pair $(\Sigma X,\Sigma A)$ is a power couple where the power map $\rho_a$ is given by
\[
\rho_a(e^{ti}\wedge x)=e^{ati}\wedge x.
\]
\end{example}

In the special case when $(X,A)$ is a pointed relative CW-complex $(X,A)$ is called a \emph{CW-power couple}. Let $(Y,B)$ be another CW-power couple. A \emph{map between CW-power couples} $f:(X,A)\to(Y,B)$ is a map between power couples. 

\begin{example}
The pair $(\mathbb{C},\mathbb{C}^*)$, equipped with power maps $\rho_a(z)=z^a$ for $a\in\N$, is a power couple but not a CW-power couple. On the other hand, the homotopy equivalent pair $(D^2,S^1)$ is a CW-power couple.
\end{example}


\begin{definition}\label{dfn_poly prod power map}
Let $(\underline{X},\underline{A})=\{(X_i,A_i)\}^m_{i=1}$ be a sequence of power couples. Let $\underline{a}=\{a_1,\ldots,a_m\}$ be a sequence of positive integers. If $K$ is a simplicial complex on $[m]$ then for any $\sigma\in K$ we call the map
\[\underline{a}(\sigma)\colon\namedright{(\underline{X},\underline{A})^{\sigma}}  
     {}{(\underline{X},\underline{A})^{\sigma}}\] 
defined by the product map $\rho_{a_{1}}\times\cdots\times\rho_{a_{m}}$ (which we will also denote as $(a_1, \cdots, a_m)$ for convenience) a \emph{power map of $(\underline{X},\underline{A})^{\sigma}$}, where $\rho_{a_{i}}$ is the power map of $X_{i}$ 
if $i\in\sigma$ and the power map of $A_{i}$ if $i\notin\sigma$. Taking the colimit over the faces of $K$ 
gives a \emph{power map}
\[
\underline{a}(K):(\underline{X},\underline{A})^K\longrightarrow(\underline{X},\underline{A})^K.
\] 
\end{definition} 

Let $\Delta^{m-1}$ be the abstract simplicial complex on $m$ vertices whose faces consist of 
all ordered subsequences of $\{1,\ldots,m\}$, including the emptyset. 

\begin{definition}\label{dfn_power_sequence}
Fix a positive integer $m$. A \emph{power sequence} $c$ is a map
\[
c:\Delta^{m-1}\longrightarrow\N^m,\qquad
\sigma\mapsto(c^{\sigma}_1,\ldots,c^{\sigma}_m)
\]
such that $c^{\sigma}_i=1$ for $i\notin\sigma$ and $c^{\tau}_i$ divides $c^{\sigma}_i$ for $\tau\subseteq\sigma$ and $1\leq i\leq m$.
\end{definition} 

Note that in some cases we only use the power sequence restricted to a simplicial complex $K\subset \Delta^{m-1}$, and that any power sequence defined on $K$ can be extended to $\Delta^{m-1}$.
The division requirement in the definition of a power sequence implies that 
$c_{i}^{\{i\}}$ divides $c_{i}^{\sigma}$ for every $1\leq i\leq m$ and every $\sigma\in\Delta^{m-1}$. 
This leads to the special case when the power sequence is ``generated" by the sets $\{i\}$ 
for $1\leq i\leq m$. 

\begin{definition} 
A power sequence $c$ is called \emph{minimal} if for every $\sigma\in\Delta^{m-1}$ we have 
$(c_{i}^{\sigma},\ldots,c_{m}^{\sigma})$ with $c_{i}^{\sigma}=1$ if $i\notin\sigma$ and 
$c_{i}^{\sigma}=c_{i}^{\{i\}}$ for $i\in\sigma$. 
\end{definition}

\begin{example} 
\label{m=3ex} 
Let $m=3$ and let 
\(c\colon\Delta^{2}\longrightarrow\mathbb{N}^{3}\) 
be the power sequence defined minimally by $c_{1}^{\{1\}}=p$, $c_{2}^{\{2\}}=q$ and 
$c_{3}^{\{3\}}=r$ for natural numbers $p,q,r$. Then $c$ sends 
\[\{1\}\mapsto (p,1,1)\qquad \{2\}\mapsto (1,q,1)\qquad \{3\}\mapsto (1,1,r)\]  
\[(1,2)\mapsto (p,q,1)\qquad (1,3)\mapsto (p,1,r)\qquad (2,3)\mapsto (1,q,r)\qquad (1,2,3)\mapsto (p,q,r).\] 
This example could be adjusted to a non-minimal power sequence, by for example, changing 
$(2,3)\mapsto (1,q,sr)$ and $(1,2,3)\mapsto (p,q,tsr)$ for any $s>1$ and $t\geq 1$. 
\end{example} 

\begin{definition}\label{dfn_c^sigma_tau}
Let $(\underline{X},\underline{A})=\{(X_i,A_i)\}^m_{i=1}$ be a sequence of power couples, let $c$ be a power sequence and let $K$ be a simplicial complex on $[m]$. If $\tau$ and $\sigma$ are faces of $\Delta^{m-1}$ such that $\tau\subseteq\sigma$, then $\displaystyle\frac{c^{\sigma}_i}{c^{\tau}_i}$ is a positive integer for each $i$. Let
\[
\underline{c}^{\sigma/\tau}=\underline{c}^{\sigma/\tau}(K):(\underline{X},\underline{A})^K\longrightarrow (\underline{X},\underline{A})^K
\]
be the power map $\displaystyle\bigg(\frac{c^{\sigma}_1}{c^{\tau}_1},\ldots,\frac{c^{\sigma}_m}{c^{\tau}_m}\bigg)$ in Definition~\ref{dfn_poly prod power map}.
\end{definition} 
 
\begin{example} 
In the minimal case of Example~\ref{m=3ex}, taking $\sigma=(1,2,3)$ and $\tau=(1,3)$ 
gives $\underline{c}^{\sigma/\tau}=(1,q,1)$ while the non-minimal case gives $\underline{c}^{\sigma/\tau}=(1,q,t)$. 
\end{example} 

\begin{lemma}\label{lemma_degree map composite}
For faces $\mu\subseteq\tau\subseteq\sigma$ in $\Delta^{m-1}$, the composites
\[
(\underline{X},\underline{A})^K\overset{\underline{c}^{\sigma/\tau}}{\longrightarrow}(\underline{X},\underline{A})^K\overset{\underline{c}^{\tau/\mu}}{\longrightarrow}(\underline{X},\underline{A})^K\qquad\mbox{and}\qquad 
(\underline{X},\underline{A})^K\overset{\underline{c}^{\tau/\mu}}{\longrightarrow}(\underline{X},\underline{A})^K\overset{\underline{c}^{\sigma/\tau}}{\longrightarrow}(\underline{X},\underline{A})^K
\]
both equal $\underline{c}^{\sigma/\mu}:(\underline{X},\underline{A})^K\longrightarrow (\underline{X},\underline{A})^K$.
\end{lemma}

\begin{proof}
The power map on each $(X_i,A_i)$ is associative and commutative. Since
\[
\frac{c^{\sigma}_i}{c^{\mu}_i}=\frac{c^{\sigma}_i}{c^{\tau}_i}\cdot\frac{c^{\tau}_i}{c^{\mu}_i}=\frac{c^{\tau}_i}{c^{\mu}_i}\cdot\frac{c^{\sigma}_i}{c^{\tau}_i}
\]
for $1\leq i\leq m$, the lemma follows.
\end{proof}

\subsection{Construction of $(\underline{X},\underline{A})^{K,c}$}

Given a power sequence $c$, we will construct a space $(\underline{X},\underline{A})^{K,c}$, called the weighted polyhedral product, as a generalization of $(\underline{X},\underline{A})^K$. 

\begin{definition}\label{dfn_weighted poly prod system}
Given a positive integer $m$, a sequence $(\underline{X},\underline{A})=\{(X_i,A_i)\}^m_{i=1}$ of power couples and a power sequence $c$, a \emph{weighted polyhedral product system} $\{(\underline{X},\underline{A})^{\bullet,c},\eta(\bullet)\}$ consists of a functor $(\underline{X},\underline{A})^{\bullet,c}:\Delta^{m-1}\longrightarrow TOP_*$ and a pointed map $\eta(\sigma):(\underline{X},\underline{A})^{\sigma}\longrightarrow (\underline{X},\underline{A})^{\sigma,c}$ for each simplex $\sigma\in\Delta^{m-1}$ 
satisfying: 
\begin{enumerate}[label=(\roman*)]
\item	$(\underline{X},\underline{A})^{\emptyset,c}=(\underline{X},\underline{A})^{\emptyset}=\prod^m_{i=1}A_i$ and $\eta(\emptyset):(\underline{X},\underline{A})^{\emptyset}\longrightarrow(\underline{X},\underline{A})^{\emptyset,c}$ is the identity map;
\item\label{dgrm_weighted PP pushout for sigma=1}
for $1\leq i\leq m$, the space $(\underline{X},\underline{A})^{\{i\},c}$ and the maps $\eta(\{i\})$ and $\imath_{\emptyset,c}^{\{i\},c}$ are defined by the pushout 
\[
\xymatrix{
(\underline{X},\underline{A})^{\emptyset}\ar[r]^{\imath_{\emptyset}^{\{ i \}}}\ar[d]^-{
\eta(\partial \{i\})}
&(\underline{X},\underline{A})^{\{i\}}\ar[d]^-{\eta(\{i\})}\\
(\underline{X},\underline{A})^{\emptyset,c}\ar[r]^-{\imath^{\{i\},c}_{\emptyset,c}}								&(\underline{X},\underline{A})^{\{i\},c}
}
\]
\item\label{dgrm_weighted PP pushout for partial sigma}
for $\sigma\in\Delta^{m-1}$ with $|\sigma|\geq1$, the space $(\underline{X},\underline{A})^{\sigma,c}$ and the maps $\eta(\sigma)$ and $\imath_{\partial\sigma,c}^{\sigma,c}$ are defined by the pushout 
\[
\xymatrix{
(\underline{X},\underline{A})^{\partial\sigma}\ar[r]^-{\imath^{\sigma}_{\partial\sigma}}\ar[d]^-{\eta(\partial\sigma)}	&(\underline{X},\underline{A})^{\sigma}\ar[d]^-{\eta(\sigma)}\\
(\underline{X},\underline{A})^{\partial\sigma,c}\ar[r]^-{\imath^{\sigma, c}_{\partial\sigma, c}}		&(\underline{X},\underline{A})^{\sigma,c}
}
\]
where, for $\tau\rightarrow\sigma$ and $\imath^{\sigma, c}_{\tau, c}:(\underline{X},\underline{A})^{\tau,c}\longrightarrow(\underline{X},\underline{A})^{\sigma,c}$ the map determined by $(\underline{X},\underline{A})^{\bullet,c}$: 
\begin{enumerate}[label=(\alph*)] 
\item	$(\underline{X},\underline{A})^{\partial\sigma,c}= \mbox{colim}_{\tau\subsetneq\sigma}(\underline{X},\underline{A})^{\tau,c}$ for $\sigma\in\Delta^{m-1}$;
\item	$\eta(\partial\sigma):(\underline{X},\underline{A})^{\partial\sigma}\longrightarrow(\underline{X},\underline{A})^{\partial\sigma,c}$ is the colimit of composites
\[
(\underline{X},\underline{A})^{\tau}\overset{\underline{c}^{\sigma/\tau}}
{\longrightarrow}(\underline{X},\underline{A})^{\tau}\overset{\eta(\tau)}
{\longrightarrow}(\underline{X},\underline{A})^{\tau,c}\overset{\imath_{\tau, c}^{\partial\sigma, c}}
{\longrightarrow}(\underline{X},\underline{A})^{\partial\sigma,c}
\]
over $\tau\subsetneq\sigma$, for $\underline{c}^{\sigma/\tau}$ as in Definition~\ref{dfn_c^sigma_tau} and 
$\imath^{\partial\sigma,c}_{\tau,c}$ the inclusion into the colimit. 
\end{enumerate}
\end{enumerate}
\end{definition}

\begin{remark}
Note using (iii)(b) and the fact that  $\eta(\emptyset):(\underline{X},\underline{A})^{\emptyset}\longrightarrow(\underline{X},\underline{A})^{\emptyset,c}$ is the identity map it follows that $\eta(\partial \{ i\})=\eta(\emptyset)\circ\underline{c}^{\{i\}/\emptyset}=\underline{c}^{\{i\}/\emptyset}$. 
We can also compute that $\underline{c}^{\{i\}/\emptyset}$ is the product map $1\times\cdots\times c^{\{i\}}_i\times\cdots\times 1:\prod^m_{j=1}A_j\to\prod^m_{j=1}A_j$. For clarity we have included the special case (ii).
\end{remark}

A construction is needed to show that weighted polyhedral product systems exist and are well-defined. This is the point of the next lemma. 

\begin{lemma}\label{lemma_construct weighted poly prod}
Let $m$ be a positive integer, let $(\underline{X},\underline{A})=\{(X_i,A_i)\}^m_{i=1}$ be a sequence of power couples and let $c$ be a power sequence. Then a weighted polyhedral product system $\{(\underline{X},\underline{A})^{\bullet,c},\eta(\bullet)\}$ exists and satisfies: 
\begin{enumerate}[label=(\roman*)]
\setcounter{enumi}{3}
\item\label{inclusion as functor}	for $\mu\subseteq\tau\subseteq\sigma\in\Delta^{m-1}$, $\imath^{\sigma, c}_{\sigma, c}$ 
is the identity map and the composite $\imath^{\sigma, c}_{\tau, c}\circ\imath^{\tau ,c}_{\mu, c}$ 
equals $\imath^{\sigma, c}_{\mu, c}$;
\item	for $\tau\subsetneq\sigma\in\Delta^{m-1}$, $\imath^{\sigma, c}_{\tau, c}$ factors as the 
composite $\imath^{\sigma, c}_{\partial\sigma, c}\circ\imath^{\partial\sigma, c}_{\tau, c}$;
\item	for $\tau\subseteq\sigma\in\Delta^{m-1}$ there is a commutative diagram
\[
\xymatrix{
(\underline{X},\underline{A})^{\tau}\ar[r]^-{\imath^{\partial\sigma}_{\tau}}\ar[d]^-{\underline{c}^{\sigma/\tau}}	&(\underline{X},\underline{A})^{\partial\sigma}\ar[r]^-{\imath^{\sigma}_{\partial\sigma}}\ar[dd]^-{\eta(\partial\sigma)}	&(\underline{X},\underline{A})^{\sigma}\ar[dd]^-{\eta(\sigma)}\\
(\underline{X},\underline{A})^{\tau}\ar[d]^-{\eta(\tau)}	&	&\\
(\underline{X},\underline{A})^{\tau,c}\ar[r]^-{\imath^{\partial\sigma,c}_{\tau,c}}	&(\underline{X},\underline{A})^{\partial\sigma,c}\ar[r]^-{\imath^{\sigma,c}_{\partial\sigma,c}}	&(\underline{X},\underline{A})^{\sigma,c}
}
\]
\end{enumerate}
\end{lemma}

Note that condition~\ref{inclusion as functor} is the same as saying that the weighted polyhedral product system is a functor, we have added it because it is convenient in the proof. 

\begin{proof}
We construct $(\underline{X},\underline{A})^{\sigma,c}, \imath^{\sigma,c}_{\tau,c}$ and $\eta(\sigma)$ by induction. Let $P(i)$ be the statement``$(\underline{X},\underline{A})^{\sigma,c}, \imath^{\sigma,c}_{\tau,c}$ and $\eta(\sigma)$ exist and satisfy Properties (i) - (vi) for $|\sigma|\leq i$''. To begin, for a simplex $\sigma$ with $|\sigma|\leq 1$ we define $(\underline{X},\underline{A})^{\sigma,c},\imath^{\sigma,c}_{\emptyset,c}$ and $\eta(\sigma)$ by Properties~(i) and~(ii). Thus $P(1)$ is true.

Assume that $P(i)$ is true for $i<n\geq 2$. Take $\sigma\in\Delta^{m-1}$ with cardinality $n$. By the inductive assumption we can define $(\underline{X},\underline{A})^{\partial\sigma,c},\imath^{\partial\sigma,c}_{\tau,c}$ and $\eta(\partial\sigma)$ for $\tau\subsetneq\sigma$ by $(b),(c)$ in Definition~\ref{dfn_weighted poly prod system}~(iii). In particular, to see that $\eta(\partial\sigma)$ is well-defined let $\tau,\mu\subsetneq\sigma$ and consider the following diagram
\begin{equation}\label{diagram_eta restriction}
\xymatrix{
(\underline{X},\underline{A})^{\tau\cap\mu}\ar@{=}[r]\ar[dd]_-{\underline{c}^{\sigma/\tau\cap\mu}}
	&(\underline{X},\underline{A})^{\tau\cap\mu}\ar[r]^-{\imath^{\tau}_{\tau\cap\mu}}\ar[d]_-{\underline{c}^{\sigma/\tau}}	&(\underline{X},\underline{A})^{\tau}\ar[d]^-{\underline{c}^{\sigma/\tau}}\\
	&(\underline{X},\underline{A})^{\tau\cap\mu}\ar[r]^-{\imath^{\tau}_{\tau\cap\mu}}\ar[d]_-{\underline{c}^{\tau/\tau\cap\mu}}	&(\underline{X},\underline{A})^{\tau}\ar[dd]^-{\eta(\tau)}\\
(\underline{X},\underline{A})^{\tau\cap\mu}\ar@{=}[r]\ar[d]_-{\eta(\tau\cap\mu)}	&(\underline{X},\underline{A})^{\tau\cap\mu}\ar[d]_-{\eta(\tau\cap\mu)}	&\\
(\underline{X},\underline{A})^{\tau\cap\mu,c}\ar@{=}[r]	&(\underline{X},\underline{A})^{\tau\cap\mu,c}\ar[r]^-{\imath^{\tau,c}_{\tau\cap\mu,c}}	&(\underline{X},\underline{A})^{\tau,c}
}
\end{equation}
The top left rectangle commutes due to Lemma~\ref{lemma_degree map composite}. The top right square commutes since the power maps 
commute with the inclusions $A_i\rightarrow X_i$. The bottom left square commutes tautologically. The bottom right rectangle commutes due to Property~(vi) and the 
inductive hypothesis. Thus the outer rectangle commutes and we can define $\eta(\partial\sigma)$ to be the map induced by the colimit construction.

Define $(\underline{X},\underline{A})^{\sigma,c},\imath^{\sigma,c}_{\partial\sigma,c}$ and $\eta(\sigma)$ by the pushout diagram in Property~(iii) and for $\tau\subsetneq\sigma$ 
define $\imath^{\sigma,c}_{\tau,c}=\imath_{\partial\sigma,c}^{\tau,c}\circ\imath_{\tau,c}^{\partial\sigma,c}$ so that Property~(v) is satisfied. For $\mu\subseteq\tau\subseteq\sigma$, 
by the inductive assumption and definitions of $\imath^{\sigma,c}_{\tau,c}$ and $\imath^{\partial\sigma,c}_{\tau,c}$ we have
\[
\imath^{\sigma,c}_{\tau,c}\circ\imath^{\tau,c}_{\mu,c}=(\imath^{\sigma,c}_{\partial\sigma,c}\circ\imath^{\partial\sigma,c}_{\tau,c})\circ\imath^{\tau,c}_{\mu,c}=\imath^{\sigma,c}_{\partial\sigma,c}\circ(\imath^{\partial\sigma,c}_{\tau,c}\circ\imath^{\tau,c}_{\mu,c})=\imath^{\sigma,c}_{\partial\sigma,c}\circ\imath^{\partial\sigma,c}_{\mu,c}=\imath^{\sigma,c}_{\mu,c}.
\]
Therefore Property~(iv) holds. Next, consider the diagram from Lemma~\ref{lemma_construct weighted poly prod} $(vi)$. The left rectangle commutes by definition of $\eta(\partial\sigma)$ in Definition~\ref{dfn_weighted poly prod system}~(iii) and the right rectangle commutes due to Property~(iii)(c). The commutativity of the whole diagram implies that Property~(vi) holds. Therefore $P(n)$ is true and the induction is complete.
\end{proof} 

\begin{definition}\label{dfn_weighted poly prod}
Let $(\underline{X},\underline{A})=\{(X_i,A_i)\}^m_{i=1}$ be a sequence of power couples, let $c$ be a power sequence and let $K$ be a simplicial complex on $[m]$.  The \emph{weighted polyhedral product} $(\underline{X},\underline{A})^{K}$ is the colimit
\[
(\underline{X},\underline{A})^{K,c}=\mbox{colim}_{\sigma\subseteq K}(\underline{X},\underline{A})^{\sigma,c}.
\]
Further, suppose $\omega\subseteq\Delta^{m-1}$ is the smallest simplex containing $K$. The \emph{associated map} 
\[\eta(K):(\underline{X},\underline{A})^{K}\to(\underline{X},\underline{A})^{K,c}\] 
is the colimit of the composites
\[
(\underline{X},\underline{A})^{\sigma}\overset{\underline{c}^{\omega/\sigma}}{\longrightarrow}(\underline{X},\underline{A})^{\sigma}\overset{\eta(\sigma)}{\longrightarrow}(\underline{X},\underline{A})^{\sigma,c}\overset{\imath^{K,c}_{\sigma,c}}{\longrightarrow}(\underline{X},\underline{A})^{K,c}
\]
over $\sigma\subseteq K$. Here $(\underline{X},\underline{A})^{\sigma,c}$ and $\eta(\sigma)$ are from Definition~\ref{dfn_weighted poly prod system} and $\imath^{K,c}_{\sigma,c}$ is the inclusion into the colimit. 
\end{definition} 

Note this definition generalizes the one for $K=\partial\sigma$ above. Even more generally if $L\subset K$ we get the maps between colimits coming from the restriction of the functor
$$
\imath_L^K\colon (\underline{X},\underline{A})^{L,c}=\mbox{colim}_{\sigma\subseteq L}(\underline{X},\underline{A})^{\sigma,c}\rightarrow (\underline{X},\underline{A})^{K,c}=\mbox{colim}_{\sigma\subseteq K}(\underline{X},\underline{A})^{\sigma,c},
$$
and if $I\subset [m]$ and $V(L)\subset I$
$$ 
\imath_{L}^{\overline {L}}\colon (\underline{X},\underline{A})^{L,c}\rightarrow (\underline{X},\underline{A})^{\overline{L},c}
$$
 coming from the inclusions $(\underline{X},\underline{A})^{\sigma, c}\rightarrow (\underline{X},\underline{A})^{\overline{\sigma}, c}=(\underline{X},\underline{A})^{\sigma, c}\times
 \prod_{i\not\in I} A_i$
  at the base points. In this case $\imath_L^K$ is the composition $\imath_{\overline{L}}^K\circ \imath_L^{\overline{L}}$. 
This is one place where we need the $A_i$ to have basepoints. 

\subsection{An alternative description of $(\underline{X},\underline{A})^{K,c}$}

We give another description of the weighted polyhedral product directly in terms of a colimit. 
We are grateful to an anonymous referee for suggesting this approach. 
This is a direct construction without induction and allows an easy description of maps between various 
$(\underline{X},\underline{A})^{K,c}$ (Proposition \ref{life_made_easier}), it also would give alternative ways of describing the projections in Section \ref{sec:functoriality2}.
The inductive description is used for showing homotopy invariance, computing the cohomology and getting the suspension decomposition since homotopic properties are easier to understand for pushouts than for general colimits. 

\begin{definition}
For a simplicial complex $K$, 
consider the coequalizer 
$$
F_c(K)=colim \coprod_{\tau\rightarrow \sigma\in K} (\underline{X},\underline{A})_{\tau\rightarrow \sigma}^{\tau}
\rightrightarrows \coprod_{\sigma\in K} (X,A)_{\sigma}^{\sigma}$$
where $ (\underline{X},\underline{A})_{\tau\rightarrow \sigma}^{\tau}=(\underline{X},\underline{A})_{\tau}^{\tau}=(\underline{X},\underline{A})^{\tau}$ with the subscripts letting us know which factor we are in, the top map on each factor is $i_{\tau}^{\sigma}\colon (\underline{X},\underline{A})_{\tau\rightarrow \sigma}^{\tau} \rightarrow
(X,A)_{\sigma}^{\sigma}$ and the bottom map is 
$\underline{c}^{\sigma/\tau}\colon (\underline{X},\underline{A})_{\tau\rightarrow \sigma}^{\tau}\rightarrow (\underline{X},\underline{A})_{\tau}^{\tau}$. 

Note that if $L\subset K$ then we get an inclusion of diagrams and so  maps, $F_c(L)\rightarrow F_c(K)$ which are  natural (in all inclusions of subcomplexes). 
We define a map $A(K)\colon F_c(K)\rightarrow (\underline{X},\underline{A})^{K,c}$ by its restriction to each factor of the colimit and then will prove it is an isomorphism in Lemma \ref{it_is_an_iso}

$$A(K)(\tau)\colon (\underline{X},\underline{A})_{\tau}^{\tau}\stackrel{\eta(\tau)}{\rightarrow} 
(X,A)^{\tau,c}\stackrel{\text{include}}{\rightarrow} (X,A)^{K,c}$$
and, although this is implied by the above formula,

$$A(K)(\tau\rightarrow \sigma)\colon (\underline{X},\underline{A})_{\tau\rightarrow \sigma}^{\tau} 
\stackrel{=}{\rightarrow} (\underline{X},\underline{A})^{\tau}\stackrel{i_{\tau}^{\sigma}}{\rightarrow} (X,A)^{\sigma}
\stackrel{\eta(\sigma)}{\rightarrow} (\underline{X},\underline{A})^{\sigma,c}
\stackrel{\text{include}}{\rightarrow} (\underline{X},\underline{A})^{K,c}.$$

Note that the first map in the sequence is into a piece of the colimit that gives the upper left corner of the diagram in (iii) of Definition \ref{dfn_weighted poly prod system}.

\end{definition}

\begin{lemma}\label{getting_the_maps}
The maps above determine a map $A(K)\colon F_c(K)\rightarrow (\underline{X},\underline{A})^{K,c}$ that is natural in the simplicial complex variable.  In other words for any $L\subset K$ the following diagram commutes

\[
\xymatrix{
F_c(L)\ar[r]^-{A(L)}\ar[d]_-{\text{include}}	&(\underline{X},\underline{A})^{L,c}\ar[d]^-{\text{include}}\\
F_c(K)\ar[r]_-{A(K)}		&(\underline{X},\underline{A})^{K,c}
}
\]
where $F_c(L)\rightarrow F_c(K)$ is induced by the inclusion of colimit diagrams mentioned above. 

In addition the following diagram commutes. 

\[
\xymatrix{
(\underline{X},\underline{A})_{\tau}^{\tau}\ar[r]^-{\eta(\tau)}\ar[d]_-{\text{include}}	
&(\underline{X},\underline{A})^{\tau,c}
\\
F_c(\tau)\ar[ur]_-{A(\tau)}	
}
\]

\end{lemma}

\begin{proof}
The compatibilities of the maps leading to a map from the colimit is straightforward using $(iii)(b)$ of Definition \ref{dfn_weighted poly prod system}. The commuting of the two diagrams is straightforward. 
\end{proof}

The map $A(K)$ also has a natural inverse $B(K)$ which is constructed recursively. 

\begin{lemma}\label{it_is_an_iso}
Each $A(K)$ is an isomorphism. 
\end{lemma}

\begin{proof}
We somewhat explicitly construct the inverse $B(K)$ to $A(K)$ while noticing that the construction is natural. 
From the definition (\ref{dfn_weighted poly prod}) maps from $(\underline{X},\underline{A})^{K,c}$ are determined by the maps from 
$(\underline{X},\underline{A})^{\sigma,c}$ with $\sigma\in K$.
It is an exercise with colimits to check that $F_c(K)=colim_{\sigma\in K} F_c(\sigma)$. Since also $(\underline{X},\underline{A})^{K,c}=colim_{\sigma\in K} (\underline{X},\underline{A})^{\sigma,c}$ and as noted in Lemma  \ref {getting_the_maps} the maps $F_c$ are compatible with inclusions of subcomplexes, $F_c(\sigma)$ being an isomorphism for each $\sigma\in K$ will imply $F_c(K)$ is an isomorphism. 

So in defining $B(K)$ and checking it is the inverse of $A(K)$ we can just look at the case when $K=\sigma$ is a simplex. Note that the $K=\emptyset$ case of the lemma is true. So to do an inductive argument we can also assume that 
$A(\tau)$ has inverse $B(\tau)$ for all $\tau\subsetneq \sigma$. Using the pushout 
of Definition \ref{dfn_weighted poly prod system} that describes $(\underline{X},\underline{A})^{\sigma,c}$, $B(\sigma)$ is determined by 
the maps
$$B(\sigma)(\tau,c)\colon (\underline{X},\underline{A})^{\tau,c}\stackrel{(A(\tau))^{-1}}{\longrightarrow} F(\tau)\stackrel{\text{include}}{\longrightarrow} F(\sigma)
$$ 
where  $\tau\subsetneq \sigma$, 
together with 
for each $\tau\subset \sigma$
$$
B(\sigma)(\tau)\colon (\underline{X},\underline{A})^{\tau}\stackrel{=}{\longrightarrow} (\underline{X},\underline{A})^{\tau}_{\tau\rightarrow \sigma}\stackrel{\text{include}}{\longrightarrow} F(\sigma)
$$ 
We just need to check these define a map from the pushout and then that $A(\sigma)$ and $B(\sigma)$ are inverses of each other. That these define a map from the pushout follows from the commutativity of the following two diagrams.

\[
\xymatrix{
(\underline{X},\underline{A})^{\tau}\ar[r]^-{i_{\tau}^{\sigma}}\ar[dd]_-{=}	
&(\underline{X},\underline{A})^{\sigma}\ar[d]^-{=}	
\\
& 
(\underline{X},\underline{A})_{\sigma\rightarrow \sigma}^{\sigma}\ar[d]^-{=}	
\\
(\underline{X},\underline{A})_{\tau}^{\tau\rightarrow \sigma} \ar[r]^-{i_{\tau}^{\sigma}}
\ar[dr]_-{\text{include}} &(\underline{X},\underline{A})_{\sigma}^{\sigma} \ar[d]^-{\text{include}}
\\ & F_c(\sigma)
}
\]
here the upper rectangle trivially commute and the lower triangle commutes by the compatibility property of colimits. 

\[
\xymatrix{
(\underline{X},\underline{A})^{\tau}\ar[r]^-{\underline{c}^{\sigma/\tau}}\ar[d]_-{=}	
& (\underline{X},\underline{A})^{\tau}\ar[r]^-{\eta(\tau)} & (\underline{X},\underline{A})^{\tau,c}\ar[dr]^-{A(\tau)^{-1}}	
\\
(\underline{X},\underline{A})_{\tau\rightarrow\sigma}^{\tau} \ar[r]^-{\underline{c}^{\sigma/\tau}}
& (\underline{X},\underline{A})_{\tau}^{\tau} \ar[ur]_-{\eta(\tau)} \ar[rr]_-{\text{include}}
&& F_c(\tau)\ar[r]_-{\text{include}} & F_c(\sigma)
}
\]
here the commuting triangle comes from the induction hypothesis and the second commuting diagram of Lemma 
\ref{getting_the_maps}. So indeed the formulas give us an induced map 
$B(\sigma)\colon  (\underline{X},\underline{A})^{\sigma,c}\rightarrow F_c(\sigma)$. 

Note that $A(\sigma)B(\sigma)|_{ (\underline{X},\underline{A})^{\tau,c}}=\mbox{inclusion}\colon { (\underline{X},\underline{A})^{\tau,c}}\rightarrow (\underline{X},\underline{A})^{\sigma,c}$ 
since we have defined $B(\sigma)(\tau)$ using $A(\tau)^{-1}$ and  
$A(\sigma)B(\sigma)|_{ (\underline{X},\underline{A})^{\sigma}}=\eta(\sigma)\colon { (\underline{X},\underline{A})^{\sigma}}\rightarrow (\underline{X},\underline{A})^{\sigma,c}$. So $A(\sigma)B(\sigma)=\mbox{identity}$. 

When  $\tau\subsetneq \sigma$ that $B(\sigma)A(\sigma)|_{(\underline{X},\underline{A})_{\tau}^{\tau}}\colon {(\underline{X},\underline{A})_{\tau}^{\tau}}
\rightarrow F_c(\sigma)$ is the inclusion follows since $B(\tau)$ and $A(\tau)$ are inverses and by the compatibility of $A$ and $B$ with inclusions of simplices. That $B(\sigma)A(\sigma)|_{(\underline{X},\underline{A})_{\sigma}^{\sigma}}$ is the inclusion follows from the commutativity of the following diagram

\[
\xymatrix{
(\underline{X},\underline{A})_{\sigma}^{\sigma}\ar[r]^-{=}\ar[dr]_-{\eta(\sigma)}	
& (\underline{X},\underline{A})^{\sigma}
\ar[r]^-{=}
\ar[d]^-{\eta(\sigma)}
&  (\underline{X},\underline{A})_{\sigma\rightarrow \sigma}^{\sigma}
\ar[r]^-{=}
&(\underline{X},\underline{A})_{\sigma}^{\sigma} \ar[d]^-{\text{include}}
\\
& (\underline{X},\underline{A})^{\sigma,c} \ar[rr]^-{B(\sigma)}
&& F_c(\sigma)
}
\]
with the rectangle commuting by the definition of $B(\sigma)$. 
This completes the proof that $A(\sigma)$ is an isomorphism with inverse $B(\sigma)$ and the induction step of the lemma thus completing the lemma's proof. 
\end{proof}

The above description of the weighted polyhedral product makes it easy to construct maps between them. 

\begin{proposition}\label{life_made_easier}
Suppose we have two power sequences $c,c'\colon K \rightarrow \N^m $ and a map $a\colon K\rightarrow\N^m $ such that for every $\tau\subset\sigma\in K$, $c'^{\sigma/\tau}a(\sigma)=a(\tau)c^{\sigma/\tau}$ then there is a map $M(a)\colon F_c(K)\rightarrow F_{c'}(K)$ induced by the maps $ (\underline{X},\underline{A})_{\sigma}^{\sigma}\stackrel{a(\sigma)}{\rightarrow}(\underline{X},\underline{A})_{\sigma}^{\sigma}$ and  $ (\underline{X},\underline{A})_{\tau\rightarrow \sigma}^{\tau}\stackrel{a(\sigma)}{\rightarrow}(\underline{X},\underline{A})_{\tau\rightarrow\sigma}^{\tau}$.
\end{proposition}
\begin{proof}
Straightforward. 
\end{proof}

\begin{example}
Let $L$ be the $LCM$ of the $a(\sigma)$ and $a'$ given by $a'(\sigma)=\frac{L}{a(\sigma)}$. Then 
$M(a')M(a)$ induces a mulitplication self map of $F_c(K)$ that factors through $F_{c'}(K)$.  We ponder if $(X,A)^K$ will have a homotopy exponent if and only if $(X,A)^{K,c}$ has one. 
\end{example}

\subsection{Some examples}

We now give several concrete examples of weighted polyhedral products. The first shows that 
weighted polyhedral products include the usual polyhedral product as a special case.

\begin{example}\label{polyprod_reduction}
If the power sequence $c$ is minimal with $c_{i}^{\{i\}}=1$ for all $1\leq i\leq m$ then 
$c(\sigma)=(1,\ldots,1)$ for all $\sigma\in\Delta^{m-1}$ and Definitions~\ref{dfn_weighted poly prod system} 
and~\ref{dfn_weighted poly prod} imply that 
$(\underline{X},\underline{A})^{K,c}=(\underline{X},\underline{A})^{K}$ and $\eta(K)$ is the identity map.  
\end{example} 

For integers $n,k\geq 2$, the \emph{Moore space} $P^{n}(k)$ is the homotopy cofibre of the degree $k$ map on $S^{n-1}$. 

\begin{example}\label{Moore} 
Let $m=1$ and $(\underline{X},\underline{A})=\{(D^2,S^1)\}$. Since $S^{1}$ has power maps given by rotating the circle, and these can be extended to maps of the disc, the pair $(D^{2},S^{1})$ has power maps. Let $c:\Delta^{0}\to\N$ be a power sequence with $c_{1}^{\{1\}}=k$ for some positive integer $k$. Observe that $(\underline{D}^{2},\underline{S}^{1})^{\emptyset}=(\underline{D}^{2}, \underline{S}^{1})^{\emptyset,c}=S^1$, the map $\underline{c}^{\{1\}/\emptyset}$ is $S^{1}\stackrel{k}{\longrightarrow} S^{1}$, so by Definition~\ref{dfn_weighted poly prod system}~(ii),  
\[
(\underline{D}^{2},\underline{S}^{1})^{\{1\},c}=\begin{cases}
D^2					&\text{if } k=1;\\
P^2(k)	&\text{if } k>1.
\end{cases} 
\]
\end{example} 

\begin{example} 
\label{prodex} 
Consider pairs of relative $CW$-complexes $(\Sigma X_{1},\ast)$ and $(\Sigma X_{2},\ast)$. By 
Example~\ref{ex_suspension} each of these are power couples. Let 
$c\colon\Delta^{1}\longrightarrow\mathbb{N}^{2}$ be a minimal power sequence with 
$c_{1}^{\{1\}}=p$ and $c_{2}^{\{2\}}=q$. So $c^{\{1\}}=(p,1)$, $c^{\{2\}}=(1,q)$ and 
$c^{(1,2)}=(p,q)$. In terms of the usual polyhedral 
product, we have 
\[(\underline{\Sigma X},\underline{\ast})^{\emptyset}=\ast\times\ast\qquad   
      (\underline{\Sigma X},\underline{\ast})^{\{1\}}=\Sigma X_{1}\times\ast\qquad 
      (\underline{\Sigma X},\underline{\ast})^{\{2\}}=\ast\times\Sigma X_{2}\] 
where we have deliberately written all spaces as Cartesian products. 
By Definition~\ref{dfn_weighted poly prod system}, 
$(\underline{\Sigma X},\underline{\ast})^{\emptyset,c}=\ast\times\ast$, so the pushout in 
Definition~\ref{dfn_weighted poly prod system}~(ii) implies that   
\[(\underline{\Sigma X},\underline{\ast})^{\{1\},c}=\Sigma X_{1}\times\ast\qquad 
      (\underline{\Sigma X},\underline{\ast})^{\{2\},c}=\ast\times\Sigma X_{1}\] 
and $\eta(\{1\})$, $\eta(\{2\})$ are the identity maps. Observe that  
$(\underline{\Sigma X},\underline{\ast})^{\partial(1,2)}=\Sigma X_{1}\vee\Sigma X_{2}$ 
and 
$(\underline{\Sigma X},\underline{\ast})^{\partial(1,2),c}=\Sigma X_{1}\vee\Sigma X_{2}$. 
By Definition~\ref{dfn_c^sigma_tau}, 
\[c^{(1,2)}/c^{\{1\}}=(1,q)\qquad c^{(1,2)}/c^{\{2\}}=(p,1).\] 
In particular, $c^{(1,2)}/c^{\{1\}}$ induces multiplication by $1\times q$ on $\Sigma X_{1}\times\ast$ 
while $c^{(1,2)}/c^{\{2\}}$ induces multiplication by $p\times 1$ on $\ast\times\Sigma X_{2}$. 
Therefore, by (c) in Definition~\ref{dfn_weighted poly prod system}~(iii), the map $\eta(\partial(1,2))$ is the 
wedge sum of degree maps 
$\Sigma X_{1}\vee\Sigma X_{2}\stackrel{1\vee 1}{\longrightarrow}\Sigma X_{1}\vee\Sigma X_{2}$. 
Hence, by Definition~\ref{dfn_weighted poly prod system}~(iii) there is a pushout 
\begin{equation} 
\label{12po} 
\xymatrix{
    \Sigma X_{1}\vee\Sigma X_{2}\ar[r]^-{\imath}\ar[d]^-{1\vee 1}	
         & \Sigma X_{1}\times\Sigma X_{2}\ar[d]^-{\eta(1,2)}\\
     \Sigma X_{1}\vee\Sigma X_{2}\ar[r] 
         &(\underline{\Sigma X},\underline{\ast})^{(1,2),c} 
}
\end{equation}  
where $\imath$ is the inclusion. Thus 
$(\underline{\Sigma X},\underline{\ast})^{(1,2),c}\simeq\Sigma X_{1}\times\Sigma X_{2}$ 
and $\eta(1,2)$ is the identity map.  
\end{example} 

\begin{example} 
\label{pqWhitehead} 
It is instructive to vary the power sequence in Example~\ref{prodex} so that it is not minimal. This time let  
$c^{\{1\}}=(1,1)$, $c^{\{2\}}=(1,1)$ and $c^{(1,2)}=(p,q)$. As before, 
$(\underline{\Sigma X},\underline{\ast})^{\partial(1,2)}=\Sigma X_{1}\vee\Sigma X_{2}$ 
and 
$(\underline{\Sigma X},\underline{\ast})^{\partial(1,2),c}=\Sigma X_{1}\vee\Sigma X_{2}$. 
But now Definition~\ref{dfn_c^sigma_tau} implies that 
\[c^{(1,2)}/c^{\{1\}}=(p,q)\qquad c^{(1,2)}/c^{\{2\}}=(p,q).\] 
So $c^{(1,2)}/c^{\{1\}}$ induces multiplication by $p\times q$ on $\Sigma X_{1}\times\ast$ 
while $c^{(1,2)}/c^{\{2\}}$ induces multiplication by $p\times q$ on $\ast\times\Sigma X_{2}$. 
Therefore, by (c) in Definition~\ref{dfn_weighted poly prod system}~(iii), the map $\eta(\partial(1,2))$ is the 
wedge sum of degree maps 
$\Sigma X_{1}\vee\Sigma X_{2}\stackrel{p\vee q}{\longrightarrow}\Sigma X_{1}\vee\Sigma X_{2}$. 
Hence, by Definition~\ref{dfn_weighted poly prod system}~(iii) there is a pushout 
\begin{equation} 
\label{12po} 
\xymatrix{
    \Sigma X_{1}\vee\Sigma X_{2}\ar[r]^-{\imath}\ar[d]^-{p\vee q}	
         & \Sigma X_{1}\times\Sigma X_{2}\ar[d]^-{\eta(1,2)}\\
     \Sigma X_{1}\vee\Sigma X_{2}\ar[r] 
         &(\underline{\Sigma X},\underline{\ast})^{(1,2),c}
}
\end{equation}  
where $\imath$ is the inclusion. Along the top row there is a homotopy cofibration 
\[\nameddright{\Sigma(X_{1}\wedge X_{2})}{[\imath_{1},\imath_{2}]}{\Sigma X_{1}\vee\Sigma X_{2}} 
     {\imath}{\Sigma X_{1}\times\Sigma X_{2}}\] 
where $[\imath_{1},\imath_{2}]$ is the Whitehead product of the maps $\imath_{1}$ and $\imath_{2}$ including 
$\Sigma X_{1}$ and $\Sigma X_{2}$ respectively into the wedge $\Sigma X_{1}\vee\Sigma X_{2}$. 
Assuming $X_{1}$ and $X_{2}$ are path-connected, the naturality of the Whitehead product and~\cite[Proposition 3.4]{Arkowitz} imply that 
$(p\vee q)\circ [\imath_{1},\imath_{2}]\simeq [p\cdot \imath_{1},q\cdot \imath_{2}]\simeq pq[\imath_{1},\imath_{2}]$. 
Thus the pushout~(\ref{12po}) implies that there is a homotopy cofibration 
\[\lllnameddright{\Sigma(X_{1}\wedge X_{2})}{pq[\imath_{1},\imath_{2}]}{\Sigma X_{1}\vee\Sigma X_{2}}{} 
     {(\underline{\Sigma X},\underline{\ast})^{(1,2),c}},\] 
identifying $(\underline{\Sigma X},\underline{\ast})^{(1,2),c}$ as the homotopy cofibre of $pq$ times 
the Whitehead product $[\imath_{1},\imath_{2}]$. Further, the pushout~(\ref{12po}) implies that the 
map $\eta(1,2)$ is degree $p$ when restricted to $\Sigma X_{1}$, degree $q$ when restricted 
to $\Sigma X_{2}$, and degree~$1$ on the cohomology of 
$\widetilde{H}^{\ast}(\Sigma X_{1}\wedge\Sigma X_{2})$.   
\end{example} 

We next give an example that shows that the weighted polyhedral product is not an invariant 
of the homotopy type of the underlying pairs of spaces. Recall that a map of pairs 
$f\colon (X,A)\to(Y,B)$ is a \emph{homotopy equivalence of pairs} if there exists a map of 
pairs $g\colon (Y,B)\to(X,A)$ and homotopies $H:X\times[0,1]\to X$ and $K:Y\times[0,1]\to Y$ 
such that $H_{0}$ is the identity map on $X$, $H_{1}=g\circ f$, and the restriction of $H$ to 
$A\times [0,1]$ has image in $A$, while $K_{0}$ is the identity map on $Y$, $K_{1}=f\circ g$, 
and the restriction of~$K$ to $B\times [0,1]$ has image in $B$. For example, the pairs 
$(D^{2},S^{1})$ and $(\mathbb{C},\mathbb{C}^{\ast})$ are homotopy equivalent. 

For ordinary polyhedral products, there is a homotopy equivalence $(D^2,S^1)^K\simeq(\C,\C^*)^K$ 
for any simplicial complex $K$~\cite[Theorem 4.7.5]{BP}. However this need not be true for weighted 
polyhedral products. 

\begin{example}\label{counterexample}  
We will show that $(D^2,S^1)^{\{1\},c}$ and $(\C,\C^*)^{\{1\},c}$ are not homotopy equivalent when $c^{\{1\}}_1>1$. 
Write $c^{\{1\}}_1=k$ for short. By Example~\ref{Moore}, $(D^2,S^1)^{\{1\},c}\cong P^2(k)$, 
which is not contractible. We will show that $(\C,\C^*)^{\{1\},c}$ is contractible, and hence  
not homotopy equivalent to $(D^2,S^1)^{\{1\},c}$. As there is a homotopy equivalence of pairs 
$(\C,\C^*)\to(D^2,D^2-\{0\})$, it suffices to show that $(D^2,D^2-\{0\})^{\{1\},c}$ is contractible.

For $a\in\N$ let $\rho_a$ be the power map on $(D^2,D^2-\{0\})$. By 
Definition~\ref{dfn_weighted poly prod system}~(ii), $(D^2,D^2-\{0\})^{\{1\},c}$ is defined by the pushout
\[
\xymatrix{
D^2-\{0\}\ar[r]^-{\text{include}}\ar[d]_-{\rho_{k}}	&D^2\ar[d]^-{\eta(\{1\})}\\
D^2-\{0\}\ar[r]^-{\jmath}								&(D^2,D^2-\{0\})^{\{1\},c}.
}
\]
Therefore 
\begin{eqnarray*}
(D^2,D^2-\{0\})^{\{1\},c}
&\cong&D^2/z\sim w\qquad\text{if }z=e^{\frac{2\gamma i\pi}{k}w}\text{ for some }\gamma\in\Z\\
&\cong&D^2/\Z_{k}.
\end{eqnarray*}
Consider the commutative diagram
\[
\xymatrix{
D^2-\{0\}\ar[r]^-{\text{include}}\ar[d]_-{\rho_{k}}	&D^2\ar[d]^-{\imath}\ar@/^1pc/[ddr]^-{\rho_{k}}	&\\
D^2-\{0\}\ar[r]^-{\jmath}\ar@/_1pc/[drr]_-{\text{include}}	&(D^2,D^2-\{0\})^{\{1\},c}\ar@{-->}[dr]^-{f}	&\\
	&	&D^2
}
\]
where $f$ is the induced pushout map. Observe that $f$ is surjective since $\rho_{k}:D^2\to D^2$ is surjective, 
and the commutativity of the lower triangle implies that $f$ is injective. Thus $f$ is a bijection. 
A continuous bijection $g:X\to Y$ is a homeomorphism if $X$ is compact and $Y$ is Hausdorff. Since $\imath:D^2\to(D^2,D^2-\{0\})^{\{1\},c}$ is surjective, $(D^2,D^2-\{0\})^{\{1\},c}$ is compact. Clearly, $D^{2}$ 
is Hausdorff. Therefore $f:(D^2,D^2-\{0\})^{\{1\},c}\to D^2$ is a homeomorphism, implying that $(D^2,D^2-\{0\})^{\{1\},c}$ is contractible.
\end{example}

\section{Initial properties of weighted polyhedral products} 
\label{sec:properties} 

In this section we consider the behaviour of minimal power sequences as a way of 
relating polyhedral products and weighted polyhedral products in the case of pairs $(X,\ast)$, 
and then relate higher Whitehead products to weighted polyhedral products in the case 
of pairs $(\Sigma X,\ast)$. 

\subsection{Minimal power sequences} 
In Example~\ref{prodex}, the identity between $(\underline{X},\underline{\ast})^{\Delta^{1},c}$ 
and $(\underline{X},\underline{\ast})^{\Delta^{1}}$ for a minimal sequence $c$ was no accident. 
It is an instance of a general property.

\begin{lemma} 
   \label{minimal} 
   Let $c$ be a minimal power sequence. Then   
   \(\namedright{(\underline{X},\underline{\ast})^{K}}{\eta(K)}{(\underline{X},\underline{\ast})^{K,c}}\) 
   is the identity map. 
\end{lemma}

\begin{proof} 
We first show that if $\sigma\in K$ then 
\(\namedright{(\underline{X},\underline{\ast})^{\sigma}}{\eta(\sigma)}{(\underline{X},\underline{\ast})^{\sigma,c}}\) 
is the identity map. If $\sigma=\emptyset$ then by Definition~\ref{dfn_weighted poly prod system}~(i), 
$(\underline{X},\underline{\ast})^{\emptyset}=(\underline{X},\underline{\ast})^{\emptyset,c}=\ast$. 
If $\vert\sigma\vert\geq 1$ we proceed by induction. 
If $\vert\sigma\vert=1$ then $\sigma=\{i\}$ for some $1\leq i\leq m$, so by 
Definition~\ref{dfn_weighted poly prod system}~(ii) there is a pushout 
\[
\xymatrix{
(\underline{X},\underline{\ast})^{\emptyset}\ar[r]\ar[d]^-{\underline{c}^{\{i\}/\emptyset}}	&(\underline{X},\underline{\ast})^{\{i\}}\ar[d]^-{\eta(\{i\})}\\
(\underline{X},\underline{\ast})^{\emptyset,c}\ar[r]^-{\imath^{\{i\},c}_{\emptyset,c}}								&(\underline{X},\underline{\ast})^{\{i\},c}
}
\]
where $\underline{c}^{\{i\}/\emptyset}$ is the product map $1\times\cdots\times c^{\{i\}}_i\times\cdots\times 1:\prod^m_{j=1}\ast\to\prod^m_{j=1}\ast$. Since $\underline{c}^{\{i\}/\emptyset}$ is the identity map, 
the pushout implies that $\eta(\{i\})$ is the identity map. 

Now suppose that $\vert\sigma\vert=n$ for some $1<n\leq m$. By inductive hypothesis, 
if $\tau\subsetneq\sigma$ then 
\(\namedright{(\underline{X},\underline{\ast})^{\tau}}{\eta(\tau)}{(\underline{X},\underline{\ast})^{\tau,c}}\) 
is the identity map. By definition of the polyhedral product, 
$(\underline{X},\underline{\ast})^{\partial\sigma}= 
     \mbox{colim}_{\tau\subsetneq\sigma}(\underline{X},\underline{\ast})^{\tau}$, 
and by definition of the weighted polyhedral product,  
$(\underline{X},\underline{\ast})^{\partial\sigma,c}= 
     \mbox{colim}_{\tau\subsetneq\sigma}(\underline{X},\underline{\ast})^{\tau,c}$. 
By part~(c) of Definition~\ref{dfn_weighted poly prod system}~(iii), the map 
\(\namedright{(\underline{X},\underline{\ast})^{\partial\sigma}}{\eta(\partial\sigma)} 
      {(\underline{X},\underline{\ast})^{\partial\sigma,c}}\) 
is the colimit of composites 
\begin{equation} 
  \label{minind} 
(\underline{X},\underline{\ast})^{\tau}\overset{\underline{c}^{\sigma/\tau}}{\longrightarrow}(\underline{X},\underline{\ast})^{\tau}\overset{\eta(\tau)}{\longrightarrow}(\underline{X},\underline{\ast})^{\tau,c}\overset{\imath_{\tau,c}^{\partial\sigma,c}}{\longrightarrow}(\underline{X},\underline{\ast})^{\partial\sigma,c}
\end{equation}    
for $\tau\subsetneq\sigma$ where $\underline{c}^{\sigma/\tau}$ is as in Definition~\ref{dfn_c^sigma_tau} 
and $\iota_{\tau,c}^{\partial\sigma,c}$ is the inclusion into the colimit. By definition, 
$\underline{c}^{\sigma/\tau}=\bigg(\frac{c^{\sigma}_1}{c^{\tau}_1},\ldots,\frac{c^{\sigma}_m}{c^{\tau}_m}\bigg)$. 
Since $c$ is minimal there are three cases: if $i\in\tau$ then $i\in\sigma$ so 
$\frac{c_{i}^{\sigma}}{c_{i}^{\tau}}=\frac{c_{i}^{\{i\}}}{c_{i}^{\{i\}}}=1$; if $i\notin\tau$ but $i\in\sigma$ 
then $\frac{c_{i}^{\sigma}}{c_{i}^{\tau}}=\frac{c_{i}^{\{i\}}}{1}=c_{i}^{\{i\}}$; 
and if $i\notin\tau$ and $i\notin\sigma$ then 
$\frac{c_{i}^{\sigma}}{c_{i}^{\tau}}=\frac{1}{1}=1$. Thus the map 
\(\namedright{(\underline{X},\underline{\ast})^{\tau}}{\underline{c}^{\sigma/\tau}}{(\underline{X},\underline{\ast})^{\tau}}\) 
is the product $\prod_{i=1}^{m}\frac{c_{i}^{\sigma}}{c_{i}^{\tau}}$ where each factor is $1$ 
except when $i\notin\tau$ but $i\in\sigma$. But in those cases the space in coordinate $i$ 
in $(\underline{X},\underline{\ast})^{\tau}$ is $\ast$, implying that $c_{i}^{\{i\}}$ is also 
the identity map. Thus $\underline{c}^{\sigma/\tau}$ is the identity map. Therefore, by~(\ref{minind}) 
the map $\eta(\partial\sigma)$ between colimits is determined by the maps $\eta(\tau)$ 
for $\tau\subsetneq\sigma$, all of which are assumed to be the identity map. Hence 
$\eta(\partial\sigma)$ is the identity map, and this implies from the pushout in 
Definition~\ref{dfn_weighted poly prod system}~(iii) that $\eta(\sigma)$ is the identity map. 

Turning to $K$, suppose that $\omega\in\Delta^{m-1}$ is the smallest simplex containing $K$. By definition, $(\underline{X},\underline{\ast})^{K,c}=\mbox{colim}_{\sigma\subseteq K}(\underline{X},\underline{\ast})^{\sigma,c}$ and the map $\eta(K)\colon(\underline{X},\underline{\ast})^{K}\to(\underline{X},\underline{\ast})^{K,c}$ is the colimit of the composites
\[
(\underline{X},\underline{\ast})^{\sigma}\overset{\underline{c}^{\omega/\sigma}}{\longrightarrow}(\underline{X},\underline{\ast})^{\sigma}\overset{\eta(\sigma)}{\longrightarrow}(\underline{X},\underline{\ast})^{\sigma,c}\overset{\imath^{K,c}_{\sigma,c}}{\longrightarrow}(\underline{X},\underline{\ast})^{K,c}
\]
over $\sigma\subseteq K$. Arguing as above, the minimality of $c$ implies that $\underline{c}^{\omega/\sigma}$ is the identity map. Thus the fact that $\eta(\sigma)$ is the identity map for all $\sigma\subseteq K$ implies that $\eta(K)$ is also the identity map. 
\end{proof}

%
%
%

\subsection{Relation to higher Whitehead products}\label{relWH}
Example~\ref{pqWhitehead} indicates that multiples of Whitehead products play a role 
in the attaching maps for weighted polyhedral products. It is reasonable to ask to what extent 
higher Whitehead products play a similar role. To define terms, fix $m\geq 2$ and let 
$FW(\Sigma X_{1},\ldots,\Sigma X_{m})$ be the subspace of 
$\Sigma X_{1}\times\cdots\times\Sigma X_{m}$ defined by  
\[FW(\Sigma X_{1},\ldots,\Sigma X_{m})=\{(x_{1},\ldots,x_{m})\in\Sigma X_{1}\times\cdots\times\Sigma X_{m}\mid 
           \mbox{at least one $x_{i}$ is $\ast$}\}.\] 
The space $FW(\Sigma X_{1},\ldots,\Sigma X_{m})$ is called the \emph{fat wedge}. 
By~\cite{P} there is a homotopy cofibration 
\[\nameddright{\Sigma^{m-1}(X_{1}\wedge\cdots\wedge X_{m})}{\phi_{m}}{FW(\Sigma X_{1},\ldots,\Sigma X_{m})} 
       {}{\Sigma X_{1}\times\cdots\times\Sigma X_{m}}\] 
where $\phi_{m}$ is called the \emph{universal higher Whitehead product}. Note that when $m=2$ we have 
$FW(\Sigma X_{1},\Sigma X_{2})=\Sigma X_{1}\vee\Sigma X_{2}$ and $\phi_{2}$ is the 
usual Whitehead product.

\begin{definition}\label{dfn_higher WH prod}
Let $Y$ be a pointed space and for $1\leq i\leq m$ let $f_i:\Sigma X_i\longrightarrow Y$ be a pointed map. Suppose that there is an extension of the wedge sum $\bigvee^m_{i=1}f_i:\bigvee^m_{i=1}\Sigma X_i\longrightarrow Y$ to a map $f:FW(\Sigma X_1,\ldots,\Sigma X_m)\longrightarrow Y$. Then the composite
\[
\Sigma^{m-1}(X_1\wedge\cdots\wedge X_m)\overset{\phi_m}{\longrightarrow}FW(\Sigma X_1,\ldots,\Sigma X_m)\overset{f}{\longrightarrow}Y
\]
is called a \emph{higher Whitehead product of $f_1,\ldots,f_m$}. 
\end{definition} 

Note that there may be another extension
\(\namedright{f'\colon FW(\Sigma X_{1},\ldots,\Sigma X_{m})}{}{Y}\) 
of $\bigvee^m_{i=1}f_i$ and the homotopy classes of the higher Whitehead products $f'\circ\phi_m$ and $f\circ\phi_m$ may be different. The set of higher Whitehead products of $f_1,\ldots,f_m$ is denoted by $[f_1,\ldots,f_m]$. 

In terms of polyhedral products, 
$FW(\Sigma X_{1},\ldots,\Sigma X_{m})=(\underline{\Sigma X},\underline{\ast})^{\partial\Delta^{m-1}}$, 
$\Sigma X_{1}\times\cdots\times\Sigma X_{m}=(\underline{\Sigma X},\underline{\ast})^{\Delta^{m-1}}$, 
and the inclusion 
\(\namedright{FW(\Sigma X_{1},\ldots,\Sigma X_{m})}{}{\Sigma X_{1}\times\cdots\times\Sigma X_{m}}\) 
is the map of polyhedral products induced by the inclusion 
\(\namedright{\partial\Delta^{m-1}}{}{\Delta^{m-1}}\). 
Thus for a weighted polyhedral product system based on a sequence of power couples 
$(\underline{\Sigma X},\underline{\ast})$, if $\sigma=(i_{1},\ldots,i_{k})$ is a face of 
$\Delta^{m-1}$ then the pushout in Definition~\ref{dfn_weighted poly prod system}~(iii) 
expands to a homotopy cofibration diagram 
\begin{equation} 
\label{chigherWh} 
\xymatrix{
\Sigma^{k-1}(X_{i_{1}}\wedge\cdots\wedge X_{i_{k}})\ar[rr]^-{\phi_{k}}\ar@{=}[d]
& &  (\underline{\Sigma X},\underline{\ast})^{\partial\sigma}\ar[r]^-{\imath^{\sigma}_{\partial\sigma}}\ar[d]^-{\eta(\partial\sigma)}	&(\underline{\Sigma X},\underline{\ast})^{\sigma}\ar[d]^-{\eta(\sigma)}\\
\Sigma^{k-1}(X_{i_{1}}\wedge\cdots\wedge X_{i_{k}})\ar[rr]^-{\eta(\partial\sigma)\circ\phi_{k}} 
& & (\underline{\Sigma X},\underline{\ast})^{\partial\sigma,c}\ar[r]^-{\imath^{\sigma,c}_{\partial\sigma,c}}		&(\underline{\Sigma X},\underline{\ast})^{\sigma,c}.
}
\end{equation} 
The bottom row of~(\ref{chigherWh}) shows that the ``attaching map" for the ``top complex" of 
$(\underline{\Sigma X},\underline{\ast})^{\sigma,c}$ is $\eta(\partial\sigma)\circ\phi_{k}$.  
We will show that this attaching map is a higher Whitehead product related to the power sequence $c$.  

To compress notation, for $1\leq i\leq m$, let 
$\imath_i:\Sigma X_i=(\underline{\Sigma X},\underline{\ast})^{\{i\},c}\to 
     (\underline{\Sigma X},\underline{\ast})^{\partial\sigma,c}$ 
be the inclusion $\imath^{\partial\sigma,c}_{\{i\},c}$  in Definition~\ref{dfn_weighted poly prod system}.

\begin{lemma} 
   If $\sigma=(i_{1},\ldots,i_{k})\in\Delta^{m-1}$ and the power sequence $c$ has the property 
   that $c_{i}^{\{i\}}=1$ for all $i\in\sigma$ then $\eta(\partial\sigma)\circ\phi_k$ is a higher 
   Whitehead product in\vspace{2mm} 
   $\bigg[c^{\sigma}_{i_1}\imath_{i_1},\ldots,c^{\sigma}_{i_k}\imath_{i_k}\bigg]$.
\end{lemma}

\begin{proof}
By Definition~\ref{dfn_higher WH prod} it suffices to show that $\eta(\partial\sigma)$ is an extension of the wedge sum
\[
\bigvee^k_{s=1}c^{\sigma}_{i_s}\imath_{i_s}\colon\bigvee^k_{s=1}\Sigma X_{i_s}\longrightarrow(\underline{\Sigma X},\underline{\ast})^{\partial\sigma,c}.
\] 

Applying Lemma~\ref{lemma_construct weighted poly prod}~(vi) to each vertex inclusion 
\(\{i\}\rightarrow\partial\sigma\) 
implies that there exists a commutative diagram 
\[
\xymatrix{
\Sigma X_{i_s}\ar[r]^-{\text{incl}}\ar[d]_-{c^{\sigma}_{i_s}}	
  & (\underline{\Sigma X},\underline{\ast})^{\partial\sigma}=FW(\Sigma X_{i_{1}},\ldots,\Sigma X_{i_{k}})\ar[dd]^-{\eta(\partial\sigma)}\\
\Sigma X_{i_s}\ar[d]_-{\eta(\{i\})}	&\\
\Sigma X_{i_s}\ar[r]^-{\imath_{i_s}} &(\underline{\Sigma X},\underline{\ast})^{\partial\sigma,c}.
}
\]
Since each $c_{i}^{\{i\}}=1$, Definition~\ref{dfn_weighted poly prod system}~(ii) implies that 
$\eta(\{i\}):\Sigma X_t\to\Sigma X_t$ is the identity map. 
Therefore, if $V$ is the vertex set of $\sigma$ then 
$(\underline{\Sigma X},\underline{\ast})^{V}=(\underline{\Sigma X},\underline{\ast})^{V,c}=\bigvee^k_{s=1}\Sigma X_{i_s}$ 
and there is a commutative diagram  
\[
\xymatrix{
\bigvee^k_{s=1}\Sigma X_{i_s}\ar[r]^-{\text{incl}}\ar[d]_-{\bigvee^k_{s=1}c^{\sigma}_{i_s}}	&FW(\Sigma X_1,\ldots,\Sigma X_m)\ar[d]^-{\eta(\partial\sigma)}\\
\bigvee^k_{s=1}\Sigma X_{i_s}\ar[r]^-{\bigvee^k_{s=1}\imath_{i_s}}			&(\underline{\Sigma X},\underline{\ast})^{\partial\sigma,c}.
}
\]
Thus $\eta(\partial\sigma)$ is an extension of $\bigvee^k_{s=1}c^{\sigma}_{i_s}$.
\end{proof}

\section{Functoriality of weighted polyhedral products with respect to maps of power couples} 
\label{sec:functoriality1} 

In this section we show that, for a given power sequence $c$, the weighted polyhedral product is functorial with 
respect to maps of power couples. Recall that a map of pairs $f\colon (X,A)\to(Y,B)$ is a \emph{homeomorphism of pairs} if there exists a map of pairs $g\colon (Y,B)\to(X,A)$ such that $g\circ f$ is the identity map on $X$ and restricts to the identity map on $A$ while $f\circ g$ is the identity map on~$Y$ and restricts to the identity map on $B$. A homotopy equivalence of pairs was defined just before Example~\ref{counterexample}.

In what follows, the restriction of the map 
\(\llnamedright{\prod_{i=1}^{m} X_{i}}{\prod_{i=1}^{m}\varphi_{i}}{\prod_{i=1}^{m} Y_{i}}\)  
to $(\underline{X},\underline{A})^{\sigma}$ will be written as 
\(\llnamedright{(\underline{X},\underline{A})^{\sigma}}{\prod^m_{i=1}\varphi_i}{(\underline{X},\underline{A})^{\sigma}}\).

\begin{proposition}\label{lemma_map between weight poly prod}
Let $c$ be a power sequence, let $(\underline{X},\underline{A})=\{(X_i,A_i)\}^m_{i=1}$ and $(\underline{Y},\underline{B})=\{(Y_i,B_i)\}^m_{i=1}$ be sequences of power couples, and let $\{(\underline{X},\underline{A})^{\bullet,c},\eta(\bullet)\}$ and $\{(\underline{Y},\underline{B})^{\bullet,c},\widetilde{\eta}(\bullet)\}$ be the associated weighted polyhedral products systems. Suppose that $\{\varphi_i:(X_i,A_i)\longrightarrow(Y_i,B_i)\}^m_{i=1}$ is a sequence of maps between power couples. Then there exists a unique collection of maps 
\[\{\varphi(\sigma):(\underline{X},\underline{A})^{\sigma,c}\longrightarrow(\underline{Y},\underline{B})^{\sigma,c}|\sigma\in\Delta^{m-1}\}\] 
such that, for each $\tau\subseteq\sigma\in\Delta^{m-1}$, the following diagrams commute
\begin{enumerate}[label=(\alph*)]
\begin{minipage}{0.45\linewidth}
\item\label{diagram_unique weight poly prod1}
\hspace{0.5cm}
$\xymatrix{
(\underline{X},\underline{A})^{\tau,c}\ar[r]^-{\varphi(\tau)}\ar[d]^-{\imath^{\sigma,c}_{\tau,c}}	&(\underline{Y},\underline{B})^{\tau,c}\ar[d]^-{\widetilde{\imath}^{\sigma,c}_{\tau,c}}\\
(\underline{X},\underline{A})^{\sigma,c}\ar[r]^-{\varphi(\sigma)}		&(\underline{Y},\underline{B})^{\sigma,c}
}$
\end{minipage}
\hfill
\begin{minipage}{0.45\linewidth}
\item\label{diagram_unique weight poly prod2}
\hspace{0.5cm}
$\xymatrix{
(\underline{X},\underline{A})^{\sigma}\ar[r]^-{\prod^m_{i=1}\varphi_i}\ar[d]^-{\eta(\sigma)}	&(\underline{Y},\underline{B})^{\sigma}\ar[d]^-{\widetilde{\eta}(\sigma)}\\
(\underline{X},\underline{A})^{\sigma,c}\ar[r]^-{\varphi(\sigma)}		&(\underline{Y},\underline{B})^{\sigma,c}.
}$
\end{minipage}
\end{enumerate}
Further, if the maps $\varphi_i$ have the additional property of being homeomorpihsms of pairs for $1\leq i\leq m$ then the maps $\varphi(\sigma)$ are homeomorphisms of pairs for every $\sigma\in\Delta^{m-1}$, and if the pairs $(X_i,A_i)$ are CW-power couples and the maps $\varphi_i$ are homotopy equivalences of pairs for $1\leq i\leq m$ then the maps~$\varphi(\sigma)$ are homotopy equivalences of pairs for every $\sigma\in\Delta^{m-1}$. 
\end{proposition}

\begin{proof}
We construct the sequence of maps $\{\varphi(\sigma)|\sigma\in\Delta^{m-1}\}$ by induction on the cardinality of $\sigma$. Let $P(i)$ be the statement that ``there exist unique pointed maps $\varphi(\sigma):(\underline{X},\underline{A})^{\sigma,c}\longrightarrow(\underline{Y},\underline{A})^{\sigma,c}$ that make Diagrams~\ref{diagram_unique weight poly prod1} and~\ref{diagram_unique weight poly prod2} commute for $|\sigma|\leq i$''. When $\vert\sigma\vert=0$ then $\sigma=\emptyset$, and by Definition~\ref{dfn_weighted poly prod system}, we have
\[
(\underline{X},\underline{A})^{\emptyset}=(\underline{X},\underline{A})^{\emptyset,c}=\prod^m_{j=1}A_j 
\qquad\mbox{and}\qquad 
(\underline{Y},\underline{B})^{\emptyset}=(\underline{Y},\underline{B})^{\emptyset,c}=\prod^m_{j=1}B_j,
\]
and $\eta(\emptyset)$ and $\tilde{\eta}(\emptyset)$ are the identity maps. Define $\varphi(\emptyset)$ to be $\prod^m_{j=1}\varphi_j$. Then Diagrams~\ref{diagram_unique weight poly prod1} and~\ref{diagram_unique weight poly prod2} coincide and become
\[
\xymatrix{
\prod^m_{j=1}A_j\ar[r]^-{\prod^m_{j=1}\varphi_j}\ar[d]^-{=}	&\prod^m_{j=1}B_j\ar[d]^-{=}\\
\prod^m_{j=1}A_j\ar[r]^-{\prod^m_{j=1}\varphi_j}			&\prod^m_{j=1}B_j.
}
\]
This clearly commutes and the map $\varphi(\emptyset)=\prod^m_{j=1}\varphi_j$ is the unique map making the diagram commute, so $P(0)$ is true.

Assume that $P(i)$ is true for $i<n$. Take $\sigma\in\Delta^{m-1}$ with $|\sigma|=n$. For any $\mu\subseteq\tau\subsetneq\sigma$ consider the diagram 
\begin{equation}\label{diagram_map between weight poly prod well def}
\xymatrix{
(\underline{X},\underline{A})^{\mu,c}\ar[r]^-{\varphi(\mu)}\ar[d]^-{\imath^{\tau,c}_{\mu,c}}	&(\underline{Y},\underline{B})^{\mu,c}\ar[r]^-{\tilde{\imath}^{\partial\sigma,c}_{\mu,c}}\ar[d]^-{\tilde{\imath}^{\tau,c}_{\mu,c}}	&(\underline{Y},\underline{B})^{\partial\sigma,c}\ar[d]^-{=}\\
(\underline{X},\underline{A})^{\tau,c}\ar[r]^-{\varphi(\tau)}	&(\underline{Y},\underline{B})^{\tau,c}\ar[r]^-{\tilde{\imath}^{\partial\sigma,c}_{\tau,c}}	&(\underline{Y},\underline{B})^{\partial\sigma,c}
}
\end{equation}
The left square commutes due to the inductive hypothesis and the right square commutes due to the definition of $\tilde{\imath}^{\partial\sigma,c}_{\tau,c}$ as the inclusion into the colimit. Let 
\[\varphi(\partial\sigma):(\underline{X},\underline{A})^{\partial\sigma,c}\longrightarrow(\underline{Y},\underline{B})^{\partial\sigma,c}\] 
be the colimit of the composites
\[
(\underline{X},\underline{A})^{\tau,c}\overset{\varphi(\tau)}{\longrightarrow}(\underline{Y},\underline{B})^{\tau,c}\overset{\tilde{\imath}^{\partial\sigma,c}_{\tau,c}}{\longrightarrow}(\underline{Y},\underline{B})^{\partial\sigma,c}
\]
over $\tau\subsetneq\sigma$. Then~(\ref{diagram_map between weight poly prod well def}) implies that $\varphi(\partial\sigma)$ is well-defined and it makes the diagram
\begin{equation} 
\label{tildeimathvarphi} 
\xymatrix{
(\underline{X},\underline{A})^{\tau,c}\ar[r]^-{\varphi(\tau)}\ar[d]^-{\imath^{\partial\sigma,c}_{\tau,c}}	&(\underline{Y},\underline{B})^{\tau,c}\ar[d]^-{\widetilde{\imath}^{\partial\sigma,c}_{\tau,c}}\\
(\underline{X},\underline{A})^{\partial\sigma,c}\ar[r]^-{\varphi(\partial\sigma)}		&(\underline{Y},\underline{B})^{\partial\sigma,c}
}
\end{equation} 
commute. For any $\tau\subsetneq\sigma$, consider the diagram
\[
\xymatrix{
(\underline{X},\underline{A})^{\tau}\ar[r]^-{\prod^m_{j=1}\varphi_j}\ar[d]^-{\underline{c}^{\sigma/\tau}}	&(\underline{Y},\underline{B})^{\tau}\ar[d]^-{\underline{c}^{\sigma/\tau}}\ar[r]^-{\tilde{\imath}^{\partial\sigma}_{\tau}}	&(\underline{Y},\underline{B})^{\partial\sigma}\ar[dd]^-{\widetilde{\eta}(\partial\sigma)}\\
(\underline{X},\underline{A})^{\tau}\ar[r]^-{\prod^m_{j=1}\varphi_j}\ar[d]^-{\eta(\tau)}	&(\underline{Y},\underline{B})^{\tau}\ar[d]^-{\widetilde{\eta}(\tau)}		&\\
(\underline{X},\underline{A})^{\tau,c}\ar[r]^-{\varphi(\tau)}		&(\underline{Y},\underline{B})^{\tau,c}\ar[r]^-{\tilde{\imath}^{\partial\sigma,c}_{\tau,c}}	&(\underline{Y},\underline{B})^{\partial\sigma,c}.
}
\]
The top left square commutes due to the hypothesis that each $\varphi_{j}$ is a map of power couples. The bottom left square commutes since $\vert\tau\vert<n$ and the square is an instance of Diagram~\ref{diagram_unique weight poly prod2} in the inductive hypothesis. The right rectangle commutes due to the definition of $\widetilde{\eta}(\partial\sigma)$ in Definition~\ref{dfn_weighted poly prod system}~(iii). Thus the entire diagram commutes, and therefore from the outer rectangle we obtain $\tilde{\imath}_{\tau,c}^{\partial\sigma,c}\circ\varphi(\tau)\circ\eta(\tau)\circ \underline{c}^{\sigma/\tau}=\widetilde{\eta}(\partial\sigma)\circ\tilde{\imath}_{\tau}^{\partial\sigma}\circ\prod_{j=1}^{m}\varphi_{j}$. 
By~(\ref{tildeimathvarphi}), $\tilde{\imath}_{\tau,c}^{\partial\sigma,c}\circ\varphi(\tau)=\varphi(\partial\sigma)\circ\imath_{\tau,c}^{\partial\sigma,c}$, so we obtain a commutative square 
\[
\xymatrix{
(\underline{X},\underline{A})^{\tau}\ar[rr]^-{\tilde{\imath}_{\tau}^{\partial\sigma}\circ\prod^m_{j=1}\varphi_j}\ar[d]^-{\imath_{\tau,c}^{\partial\sigma,c}\circ\eta(\tau)\circ \underline{c}^{\sigma/\tau}}	& &(\underline{Y},\underline{B})^{\partial\sigma}\ar[d]^-{\widetilde{\eta}(\partial\sigma)}\\
(\underline{X},\underline{A})^{\partial\sigma,c}\ar[rr]^-{\varphi(\partial\sigma)}	&	&(\underline{Y},\underline{B})^{\partial\sigma,c}.
}
\]
By Definition~\ref{dfn_weighted poly prod system}~(iii), $\eta(\partial\sigma)$ is the colimit of the maps 
$\imath_{\tau,c}^{\partial\sigma,c}\circ\eta(\tau)\circ \underline{c}^{\sigma/\tau}$ for all $\tau\subsetneq\sigma$, 
so we obtain a commutative diagram
\begin{equation} 
\label{leftface} 
\xymatrix{
(\underline{X},\underline{A})^{\partial\sigma}\ar[r]^-{\prod^m_{j=1}\varphi_j}\ar[d]^-{\eta(\partial\sigma)}	&(\underline{Y},\underline{B})^{\partial\sigma}\ar[d]^-{\widetilde{\eta}(\partial\sigma)}\\
(\underline{X},\underline{A})^{\partial\sigma,c}\ar[r]^-{\varphi(\partial\sigma)}		&(\underline{Y},\underline{B})^{\partial\sigma,c}.
}
\end{equation} 

Next define $\varphi(\sigma):(\underline{X},\underline{A})^{\sigma,c}\to(\underline{Y},\underline{B})^{\sigma,c}$ as follows. Consider the diagram
\begin{equation}\label{diagram_construct weighted poly prod map}
\xymatrix{
	&(\underline{X},\underline{A})^{\partial\sigma}\ar[dl]_-{\eta(\partial\sigma)}\ar[rr]\ar@{-}[d]	&	&(\underline{X},\underline{A})^{\sigma}\ar[dl]^-{\eta(\sigma)}\ar[dd]^(0.725){\prod^m_{j=1}\varphi_j}\\
(\underline{X},\underline{A})^{\partial\sigma,c}\ar[rr]\ar[dd]_-{\varphi(\partial\sigma)}	& \ar[d]^-{\prod^m_{j=1}\varphi_j} & (\underline{X},\underline{A})^{\sigma,c}\ar@{..>}[dd]	&\\
	&(\underline{Y},\underline{B})^{\partial\sigma}\ar[dl]_-{\widetilde{\eta}(\partial\sigma)}\ar[rr]	&	&(\underline{Y},\underline{B})^{\sigma}\ar[dl]^-{\widetilde{\eta}(\sigma)}\\
(\underline{Y},\underline{B})^{\partial\sigma,c}\ar[rr]	&	&(\underline{Y},\underline{B})^{\sigma,c}	&
}
\end{equation}
where all maps pointing to the right are inclusions and the dotted arrow will be defined momentarily. The left face commutes by~(\ref{leftface}). The rear face clearly commutes. The top and the bottom faces are two instances of the pushout diagram in Definition~\ref{dfn_weighted poly prod system}~(iii). The universal property of the pushout in the bottom face therefore implies that there is a unique map 
\[\varphi(\sigma)\colon\namedright{(\underline{X},\underline{A})^{\sigma,c}}{}{(\underline{Y},\underline{B})^{\sigma,c}}\] 
(the dotted map in the diagram) that makes the whole diagram commute. Notice that the front and the right faces are Diagrams~\ref{diagram_unique weight poly prod1} and~\ref{diagram_unique weight poly prod2} in the assertion of the Lemma. 
Hence $P(n)$ is true and the induction is complete.

Next we show that if the maps $\varphi_j$ are homeomorphisms for $1\leq j\leq m$ then the maps $\varphi(\sigma)$ are also homeomorphisms for all $\sigma\in\Delta^{m-1}$. This is done by induction on the cardinality of $\sigma$. When $\vert\sigma\vert=0$ then $\sigma=\emptyset$ and $\varphi(\emptyset)=\prod^m_{j=1}\varphi_j$ is a homeomorphism. Assume that $\varphi(\sigma)$ is a homeomorphism for any $\sigma$ with $|\sigma|<n$. Take $\sigma$ with $|\sigma|=n$. Since $\varphi(\partial\sigma)$ is the colimit of the composites $(\underline{X},\underline{A})^{\tau,c}\overset{\varphi(\tau)}{\to}(\underline{Y},\underline{B})^{\tau,c}\rightarrow(\underline{Y},\underline{B})^{\partial\sigma,\tilde{c}}$ over $\tau\subsetneq\sigma$ and each map~$\varphi(\tau)$ is a homeomorphism, so is $\varphi(\partial\sigma)$. In~(\ref{diagram_construct weighted poly prod map}) the three solid vertical maps are homeomorphisms, so the universal property of the pushout implies that the dotted map $\varphi(\sigma)$ is also a homeomorphism. Hence the induction is complete.

Finally we show that if the pairs $(X_j,A_j)$ are CW-power couples and the maps $\varphi_j$ are homotopy equivalences for $1\leq j\leq m$ then the maps~$\varphi(\sigma)$ are also homotopy equivalences for all $\sigma\in\Delta^{m-1}$. This is proved by induction on the cardinality of~$\sigma$. When $\vert\sigma\vert=0$ then $\sigma=\emptyset$ and $\varphi(\emptyset)=\prod^m_{j=1}\varphi_j$ is a homotopy equivalence. Assume that $\varphi(\sigma)$ is a homotopy equivalence for any $\sigma$ with $|\sigma|<n$. Take~$\sigma$ with $|\sigma|=n$. By Definition~\ref{dfn_weighted poly prod system}~(iii), for any $\tau\subseteq\mu\subseteq\sigma$ the map 
\[
\mbox{colim}_{\omega\subsetneq\tau}(\underline{X},\underline{A})^{\omega,c}=(\underline{X},\underline{A})^{\partial\tau,c}\overset{\imath^{\mu, c}_{\partial\tau, c}}{\longrightarrow}(\underline{X},\underline{A})^{\mu,c}
\]
is always a cofibration and its image is a closed subspace. Then the inductive hypothesis and the Homotopy Lemma (see~\cite[Theorem 4.4]{BBCG1} or~\cite[Proposition 3.7]{WZZ}) implies that $\varphi(\partial\sigma)$ is a homotopy equivalence. In~(\ref{diagram_construct weighted poly prod map}) the three solid vertical maps are homotopy equivalences, so the Gluing Lemma~\cite[Lemma 2.1.3]{MP} implies that the dotted map $\varphi(\sigma)$ is also a homotopy equivalence. Hence the induction is complete.
\end{proof}

\begin{remark} 
There is also a slight strengthening of Proposition~\ref{lemma_map between weight poly prod}.
Suppose that $(X_i,A_i)$ and $(Y_i,B_i)$ are CW-power couples for $1\leq i\leq m$. 
Since all spaces $A_i,X_i,B_i,Y_i$ are CW-complexes, so are $(\underline{X},\underline{A})^{\sigma,c}$ and $(\underline{Y},\underline{B})^{\sigma,c}$ for any $\sigma\in\Delta^{m-1}$ and $c$. Therefore if all the maps $\varphi_i$ are weak equivalences of pairs for $1\leq i\leq m$ then they are homotopy equivalences of pairs. By Lemma~\ref{lemma_map between weight poly prod}, the resulting maps $\varphi(\sigma)$ are also homotopy equivalences of pairs for all $\sigma\in\Delta^{m-1}$.
\end{remark} 

\begin{remark} 
Proposition~\ref{lemma_map between weight poly prod} implies that for a given power sequence $c$ the weighted polyhedral product system is uniquely defined.
\end{remark}

\section{Functoriality of weighted polyhedral products with respect to projections to full subcomplexes} 
\label{sec:functoriality2} 

Let $K$ be a simplicial complex on the vertex set $[m]$. A subcomplex $L\subseteq K$ is 
a \emph{full subcomplex} if every face of $K$ on the vertex set of $L$ is also a face of $L$. 
Turning this around, if $I\subseteq [m]$ then the full subcomplex $K_{I}$ of $K$ consists 
of all the faces of $K$ whose vertices lie in $I$. 

The inclusion of $I$ into $[m]$ induces a map of polyhedral products 
\(\iota_{K_{I}}^{K}\colon\namedright{(\underline{X},\underline{A})^{K_{I}}}{}{(\underline{X},\underline{A})^{K}}\). 
Denham and Suciu~\cite{DS} showed that this map has a left inverse 
\(p_{K}^{K_{I}}\colon\namedright{(\underline{X},\underline{A})^{K}}{}{(\underline{X},\underline{A})^{K_{I}}}\). 
To see this, consider the inclusion 
\(\namedright{(\underline{X},\underline{A})^{K}}{}{\prod_{i=1}^{m} X_{i}}\) 
induced by the inclusion 
\(\namedright{K}{}{\Delta^{m-1}}\). 
Since $K_{I}$ is a full subcomplex of $K$, the composite 
\(p'\colon\nameddright{(\underline{X},\underline{A})^{K}}{}{\prod_{i=1}^{m} X_{i}} 
     {\mbox{\tiny proj}}{\prod_{i\in I} X_{i}}\) 
has image $(\underline{X},\underline{A})^{K_{I}}$. The map $p_{K}^{K_{I}}$ is the restriction of $p'$ to its image. 

More generally, let $\overline{K}_{I}$ be the simplicial complex $K_{I}$ but regarded as 
having vertex set $[m]$ instead of $I$. Here, the elements in $[m]-I$ are ghost vertices. 
Then, by the definition of the polyhedral product, there is a homeomorphism 
$(\underline{X},\underline{A})^{\overline{K}_{I}}\cong(\underline{X},\underline{A})^{K_{I}}\times\prod_{i\notin I} A_{i}$, 
and the argument above implies that there is a commutative diagram 
\begin{equation} 
  \label{KIdgrm} 
  \diagram 
      (\underline{X},\underline{A})^{K_{I}}\times\prod_{i\notin I}  A_{i}\rto^-{\cong}\drrto_{\text{proj}}  
            & (\underline{X},\underline{A})^{\overline{K}_{I}}\rto^{\imath_{{\overline K_I}}^K}
            & (\underline{X},\underline{A})^{K}\dto^{p_{K}^{K_{I}}} \\ 
      & & (\underline{X},\underline{A})^{K_{I}}. 
  \enddiagram 
\end{equation}  
Note that all the maps in this diagram commute with any power maps $c^{\tau/\sigma}$. 
The purpose of this section is to show that there are analogues to $p_{K}^{K_{I}}$ and 
diagram~(\ref{KIdgrm}) in the case of weighted polyhedral products. 

As weighted polyhedral products are defined by colimits, it will be helpful to have a face-to-face 
analogue of the retractions above for the usual polyhedral product. Suppose that $\Delta^{m-1}$ 
has vertex set $[m]$. If $\tau\in\Delta^{m-1}$, let $V(\tau)$ be its 
vertex set. Regarding $\tau$ as a simplicial complex on $[m]$, it has ghost vertices $[m]-V(\tau)$. 
Fix a subset $I\subseteq [m]$. Suppose that $\tau\in\Delta^{m-1}$ 
and $\sigma\in\Delta^{\vert I\vert -1}$. Note that $V(\tau)\subseteq [m]$ and $V(\sigma)\subseteq I$; 
suppose also that $V(\tau)\cap I\subseteq V(\sigma)$. 
Define the map $p^{\sigma}_{\tau}$ by the composite 
$$
p_\tau^\sigma\colon (\underline X,\underline A)^{\tau} 
\stackrel{\rm reorder}{\cong}
\prod_{i \in V(\tau)} X_i \times \prod_{i\in [m]- V(\tau)} A_i
\stackrel{\rm project}{\longrightarrow}
\prod_{i \in V(\tau)\cap I} X_i \times \prod_{i\in I- V(\tau)\cap I} A_i
$$
$$
\stackrel{\rm include}{\longrightarrow}
\prod_{i \in V(\sigma)} X_i \times \prod_{i\in I- V(\sigma)} A_i
\stackrel{\rm reorder}{\cong} 
(\underline X,\underline A)^{\sigma}. 
$$
We refer to the map $p^{\sigma}_{\tau}$ as a projection. Arguing as for diagram~(\ref{KIdgrm}) 
shows the following. 

\begin{lemma}\label{classical} 
Fix $I\subseteq [m]$. Suppose that $\tau\in\Delta^{m-1}$, $\sigma\in\Delta^{\vert I\vert -1}$,  
and $V(\tau)\cap I\subseteq V(\sigma)$. If $\delta$ is a face of $\tau$ and $\tau_{\tinyI}$ is the 
full subcomplex of $\tau$ on $I$, then the following two diagrams commute
\begin{equation}\label{diagram2c}
\xymatrix{
(\underline{X}, \underline{A})^{\delta}\ar[r]^{p_{\delta}^{\sigma}}
\ar[d]_{\imath_{\delta}^{\tau}}
& 
(\underline{X}, \underline{A})^{\sigma}  \\
(\underline{X}, \underline{A})^{\tau}\ar[ur]_{p_{\tau}^{\sigma}}					
}
\end{equation}

\begin{equation}\label{diagram3c}
\xymatrix{
(\underline{X}, \underline{A})^{\tau}\ar[r]^{p_{\tau}^{\tau_{\sscaleI}}}
\ar[dr]_{p_{\tau}^{\sigma}}
& (\underline{X}, \underline{A})^{\tau_{\scaleI}}   
\ar[d]^{\imath_{\tau_{\sscaleI}}^{\sigma}}\\
&(\underline{X}, \underline{A})^{\sigma}.
}
\end{equation} 
~$\qqed$ 
\end{lemma} 

The map 
\(p_{K}^{K_{I}}\colon\namedright{(\underline{X},\underline{A})^{K}}{}{(\underline{X},\underline{A})^{K_{I}}}\) 
in~(\ref{KIdgrm}) can be obtained by putting the maps $p_{\tau}^{\sigma}$ together. More precisely it is the unique map such that 
restricted to $(\underline{X},\underline{A})^{\sigma}$ for $\sigma\in K$ it is the composition $\imath_{\sigma_I}^K p_{\sigma}^{\sigma_I}$. 

In the case of weighted polyhedral products we will define analogous maps $p_{\tau,c}^{\sigma, c}$, 
prove corresponding versions of~(\ref{diagram2c}) and~(\ref{diagram3c}), and put them together to 
construct a retraction 
$$ 
p_{K,c}^{K_{I},c}\colon (\underline X,\underline A)^{K,c}\rightarrow 
(\underline X,\underline A)^{K_I, c}. 
$$  
To cohere with Definition \ref{dfn_weighted poly prod system}, an additional diagram is needed 
to relate the projection for the usual polyhedral product with the one for the weighted polyhedral product. 
This is summed up in Lemma~\ref{3diagramlemma}. 

Some remarks on notation for power sequences are necessary. Let $K$ be a simplicial 
complex on~$[m]$ and let $c$ be a power sequence. Then $c$ has $m$ components. 
Suppose that $I\subseteq [m]$. The full sucomplex $K_{I}$ is on the vertex set $I$. Rather 
than writing $c_{I}$ to restrict to those components in $I$, we simply write $c$ with the 
understanding that when paired with $K_{I}$ only those components of $c$ in~$I$ are used. 
Thus we write $(\underline{X},\underline{A})^{K_{I},c}$. Also, if $\tau\in K$ then 
$\tau_{\tinyI}\in K_{I}$ and $\tau_{\tinyI}$ is a face of $\tau$. So it makes sense to 
consider, for example, $\underline{c}^{\tau/\tau_{\sscaleI}}$, with $m$ components, as in 
Definition~\ref{dfn_c^sigma_tau}. 

\begin{lemma}\label{3diagramlemma} 
Fix $I\subseteq [m]$. Suppose that $\tau\in\Delta^{m-1}$, $\sigma\in\Delta^{\vert I\vert -1}$,  
and $V(\tau)\cap I\subseteq V(\sigma)$. Then there exists a unique map 
$$
p_{\tau,c}^{\sigma, c}\colon 
(\underline X,\underline A)^{\tau,c}\longrightarrow 
(\underline X,\underline A)^{\sigma,c}
$$ 
such that, for any face $\delta$ of $\tau$, the following diagrams commute 
\begin{equation}\label{diagram1}
\xymatrix{
(\underline{X}, \underline{A})^{\tau}\ar[r]^{p_{\tau}^{\sigma}}
\ar[d]_{\eta(\tau)\circ \underline{c}^{\sigma/\tau_{\sscaleI}}}
& (\underline{X}, \underline{A})^\sigma   
\ar[d]^{\eta(\sigma)\circ \underline{c}^{\tau/\tau_{\sscaleI}}}\\
(\underline{X}, \underline{A})^{\tau,c}\ar[r]^{p_{\tau,c}^{\sigma,c}}					&(\underline{X}, \underline{A})^{\sigma,c}
}
\end{equation}

\begin{equation}\label{diagram2}
\xymatrix{
(\underline{X}, \underline{A})^{\delta,c}\ar[r]^{p_{\delta,c}^{\sigma,c}}
\ar[d]_{\imath_{\delta,c}^{\tau,c}}
& 
(\underline{X}, \underline{A})^{\sigma,c }  \\
(\underline{X}, \underline{A})^{\tau,c}\ar[ur]_{p_{\tau,c}^{\sigma,c}}					
}
\end{equation}

\begin{equation}\label{diagram3}
\xymatrix{
(\underline{X}, \underline{A})^{\tau,c}\ar[r]^{p_{\tau,c}^{\tau_{\sscaleI},c}}
\ar[dr]_{p_{\tau,c}^{\sigma,c}}
& (\underline{X}, \underline{A})^{\tau_{\scaleI}, c}   
\ar[d]^{\imath_{\tau_{\sscaleI},c}^{\sigma,c}}\\
&(\underline{X}, \underline{A})^{\sigma,c}.
}
\end{equation}
\end{lemma}

\begin{proof}
Let $P(n)$ be the statement of the proposition when $|\tau|\leq n$. The proof proceeds in 
two cases: in the special case when $\sigma=\tau_{\tinyI}$ and then the general case for any 
$\sigma\in\Delta^{\vert I\vert -1}$. 
\medskip 

\noindent 
\textit{Case 1}: 
Assume that $\sigma=\tau_{\tinyI}$. If $n=0$ then $\tau=\sigma=\emptyset$ and Diagram~(\ref{diagram1}) 
has the form
\begin{equation}
\xymatrix{
(\underline{X}, \underline{A})^{\emptyset}=\prod_{i=1}^{m} A_{i}
\ar[rr]^{p_{\emptyset}^{\emptyset_{I}}={\rm project}}
\ar[d]_{\eta(\emptyset)={\rm id}}
& & (\underline{X}, \underline{A})^{\emptyset} =\prod_{i\in I} A_{i}  
\ar[d]^{\eta(\emptyset)={\rm id}}\\
(\underline{X}, \underline{A})^{\emptyset,c}=\prod_{i=1}^{m} A_{i}
\ar[rr]_{p_{\emptyset,c}^{\emptyset_{I},c}={\rm project}}					
&&(\underline{X}, \underline{A})^{\emptyset,c}=\prod_{i\in I} A_{i}. 
}
\end{equation} 
So defining $p_{\tau,c}^{\sigma,c}$ as the projection $p_{\emptyset, c}^{\emptyset_{I}, c}$ results in a 
unique map that makes the diagram commute. 

For $n>0$ consider the following diagram 
\begin{equation}\label{rightpo} 
\xymatrix{
(\underline{X}, \underline{A})^{\delta}\ar[r]^{\imath}
\ar[d]_{\eta(\delta)\circ \underline{c}^{\tau/\delta}}
& (\underline{X}, \underline{A})^{\partial\tau}\ar[r]^{\imath}
\ar[d]_{\eta(\partial\tau)}
&
(\underline{X}, \underline{A})^\tau   
\ar[d]^{\eta(\tau)}\\
(\underline{X}, \underline{A})^{\delta,c}\ar[r]^{\imath}
& (\underline{X}, \underline{A})^{\partial\tau,c}\ar[r]^{\imath}
&(\underline{X}, \underline{A})^{\tau,c}  
}
\end{equation}
where the maps labelled $\imath$ are shorthand for the inclusions $\imath_{\delta}^{\partial\tau}$ and so on.  
Both squares commute by Definition~(\ref{dfn_weighted poly prod system})~(iii), with the right 
square being a pushout. We will use this pushout to construct 
$p_{\tau,c}^{\tau_{\scaleI},c}$ by constructing maps from its corners. 

First, for $\delta\subsetneq\tau$, by inductive hypothesis there are maps 
$p_{\delta,c}^{\tau_{\scaleI}, c}\colon 
(\underline{X}, \underline{A})^{\delta,c} \longrightarrow 
(\underline{X}, \underline{A})^{\tau_{\scaleI},c}$  
satisfying  Diagram~(\ref{diagram2}). Since 
$(\underline{X}, \underline{A})^{\partial\tau,c}=\mbox{colim}_{\delta\subsetneq \tau}(\underline{X}, \underline{A})^{\delta,c}$
we obtain a unique map 
$p_{\partial\tau,c}^{\tau_{\scaleI}, c}\colon 
(\underline{X}, \underline{A})^{\partial\tau,c}
 \longrightarrow 
(\underline{X}, \underline{A})^{\tau_{\scaleI},c}$  
such that 
\begin{equation}\label{diagram5}
\xymatrix{
(\underline{X}, \underline{A})^{\delta,c}\ar[dr]^{p_{\delta,c}^{\tau_{\sscaleI},c}}
\ar[d]_{\imath}
\\
(\underline{X}, \underline{A})^{\partial\tau, c}   
\ar[r]_{p_{\partial\tau,c}^{\tau_{\sscaleI},c}}&(\underline{X}, \underline{A})^{\tau_{\scaleI},c}
}
\end{equation}
commutes for all $\delta\subsetneq \tau$. 

Next, we show that the right square of the following diagram commutes: 
\begin{equation}\label{squacom}
\xymatrix{
(\underline{X}, \underline{A})^{\delta}\ar[r]^{\imath}
\ar[d]_{\eta(\delta)\circ \underline{c}^{\tau/\delta}}
& (\underline{X}, \underline{A})^{\partial\tau}\ar[r]^{\imath}
\ar[d]_{\eta(\partial\tau)}
&
(\underline{X}, \underline{A})^\tau   
\ar[d]^{\eta(\tau_I)\circ \underline{c}^{\tau/\tau_{\sscaleI}}\circ p_{\tau}^{\tau_{\sscaleI}}}\\
(\underline{X}, \underline{A})^{\delta,c}\ar[r]^{\imath}
& (\underline{X}, \underline{A})^{\partial\tau,c}
\ar[r]^{p_{\partial\tau,c}^{{\tau_{\sscaleI}}, c}}
&(\underline{X}, \underline{A})^{\tau_{\scaleI},c}.  
}
\end{equation} 
Since  
$(\underline{X}, \underline{A})^{\partial\tau}=\mbox{colim}_{\delta\subsetneq \tau}(\underline{X}, \underline{A})^{\delta}$ 
we just need to show that the outside rectangle in~(\ref{squacom}) commutes 
for each $\delta\subsetneq \tau$. Since Diagram~(\ref{diagram5}) states that 
$p_{\delta,c}^{\tau_{\scaleI},c}=p_{\partial\tau,c}^{\tau_{\scaleI},c}\circ \imath_{\delta,c}^{\partial\tau, c}$,
this means we need to show that 
\begin{equation}\label{etacpi} 
\eta(\tau_{\tinyI})\circ \underline{c}^{\tau/\tau_{\sscaleI}}\circ p_{\tau}^{\tau_{\scaleI}}\circ \imath_{\delta}^{\tau}
=p_{\delta,c}^{{\tau_{\scaleI}}, c}\circ\eta(\delta)\circ \underline{c}^{\tau/\delta}.
\end{equation} 

To show~(\ref{etacpi}), consider the diagram 
\begin{equation}
\xymatrix{
(\underline{X}, \underline{A})^{\delta}
\ar[r]^{\imath}
\ar[dr]^{p_{\delta}^{\delta_{\sscaleI}}}
\ar[d]_{\eta(\delta)\circ \underline{c}^{\tau/\delta}}
& (\underline{X}, \underline{A})^{\tau} 
\ar[dr]^{p_{\tau}^{\tau_{\sscaleI}}}   
\\
(\underline{X}, \underline{A})^{\delta,c} 
\ar[dr]^{p_{\delta,c}^{\delta_{\sscaleI},c}} 
&
(\underline{X}, \underline{A})^{\delta_{\scaleI}} 
\ar[d]^{\eta(\delta_I)\circ \underline{c}^{\tau/\delta_{\sscaleI}}} \ar[r]^{\imath} 
&
(\underline{X}, \underline{A})^{\tau_{\scaleI}} 
\ar[d]^{\eta(\tau_{\scaleI})\circ \underline{c}^{\tau/\tau_{\sscaleI}}} 
\\
& (\underline{X}, \underline{A})^{\delta_{\scaleI},c} \ar[r]^{\imath}
& (\underline{X}, \underline{A})^{\tau_{\scaleI},c}
}
\end{equation}
The top right diamond commutes since inclusion commutes with projection by Lemma \ref{classical}. 
Now please keep Lemma \ref{lemma_degree map composite} in mind.
The bottom left diamond commutes since 
$\underline{c}^{\tau/\delta_{\sscaleI}}=\underline{c}^{\tau/\delta}\circ\underline{c}^{\delta/\delta_{\sscaleI}}$ and 
by an instance of Diagram~(\ref{diagram1}) precomposed with $\underline{c}^{\tau/\delta}$. The bottom right square commutes by the factorization
$\underline{c}^{\tau/\delta_{\sscaleI}}=\underline{c}^{\tau/\tau_{\sscaleI}}\circ \underline{c}^{\tau_{\sscaleI}/\delta_{\sscaleI}}$ and by precomposing an instance of Diagram~(vi) in Lemma \ref{lemma_construct weighted poly prod} with $\underline{c}^{\tau/\tau_{\sscaleI}}$. 
The commutativity of the diagram therefore implies that its outer perimeter commutes, giving 
$\eta(\tau_{\tinyI})\circ \underline{c}^{\tau/\tau_{\sscaleI}}\circ p_{\tau}^{\tau_{\scaleI}}\circ \imath_{\delta}^{\tau}
=\imath_{\delta_{\scaleI},c}^{\tau_{\scaleI},c}\circ p_{\delta,c}^{\delta_{\scaleI}, c}\circ\eta(\delta)\circ \underline{c}^{\tau/\delta}$.  
Since $\imath_{\delta_{\scaleI},c}^{\tau_{\scaleI},c}\circ p_{\delta,c}^{\delta_{\scaleI}, c}=\imath_{\delta,c}^{\tau_{\scaleI},c}$, 
we obtain~(\ref{etacpi}), as required. 

Hence the right square in Diagram~(\ref{squacom}) communtes. From the pushout in the 
right square of~(\ref{rightpo}) we therefore obtain a unique map
$p_{\tau,c}^{\tau_{\scaleI},c} \colon (X,A)^{\tau,c}\longrightarrow (X,A)^{\tau_{\scaleI},c}$ 
making Diagrams~(\ref{diagram1}) and~(\ref{diagram2}) commute for all $|\tau|\leq n$, $\sigma=\tau_{\tinyI}$ and 
$\delta\subsetneq\tau$. Diagram~(\ref{diagram2}) also trivially commutes when $\delta=\tau$. 
Diagram~(\ref{diagram3}) trivially commutes since $\sigma=\tau_{\tinyI}$. 
This completes $P(n)$ when $\sigma=\tau_{\tinyI}$.
\medskip 

\noindent 
\textit{Case 2}: 
Suppose that $\sigma\in\Delta^{\vert I\vert -1}$ and $V(\tau)\cap I\subseteq V(\sigma)$. Define
$p_{\tau,c}^{\sigma,c}$ as the composite 
\begin{equation}\label{pgeneralsigma} 
p_{\tau,c}^{\sigma,c}=\imath_{\tau_{\scaleI},c}^{\sigma,c}\circ 
p_{\tau,c}^{\tau_{\scaleI},c}.
\end{equation} 
This is the unique possibility that makes Diagram~(\ref{diagram3}) commute. 
For $\delta\subsetneq\tau$ consider the following diagram
\begin{equation}
\xymatrix{ 
(\underline{X}, \underline{A})^{\delta,c} \ar[d]^{=} \ar[r]^{p_{\delta,c}^{\delta_{\sscaleI},c}} 
& (\underline{X}, \underline{A})^{\delta_{\scaleI,c}} \ar[d]^{\imath} \\
(\underline{X}, \underline{A})^{\delta,c} \ar[d]^{\imath} \ar[r]^{p_{\delta,c}^{\tau_{\sscaleI},c}} 
& (\underline{X}, \underline{A})^{\tau_{\scaleI},c} \ar[d]^{=} \\
(\underline{X}, \underline{A})^{\tau,c} \ar[r]^{p_{\tau,c}^{\tau_{\sscaleI},c}} \ar[dr]_{p_{\tau,c}^{\sigma,c}} 
& (\underline{X}, \underline{A})^{\tau_{\scaleI},c} \ar[d]^{\imath} \\ 
& (\underline{X}, \underline{A})^{\sigma,c}
}
\end{equation}
The top square commutes by the induction hypothesis, the middle square commutes by Case~1, 
and the bottom triangle commutes by the definition of $p_{\tau,c}^{\sigma,c}$ in~(\ref{pgeneralsigma}).  
This together with Diagram~(\ref{diagram3}) and the induction hypothesis imply that Diagram~(\ref{diagram2}) commutes.

Next consider the diagram 
\begin{equation}\label{pieta}
\xymatrix{
(\underline{X}, \underline{A})^{\tau}\ar[r]^{p_{\tau}^{\tau_{\sscaleI}}}
\ar[d]_{\eta(\tau)\circ \underline{c}^{\sigma/\tau_{\sscaleI}}}
& (\underline{X}, \underline{A})^{\tau_{\scaleI}}\ar[rr]^{\imath_{\tau_{\sscaleI}}^{\sigma}}
\ar[d]^{\eta(\tau_{\tinyI})\circ\underline{c}^{\sigma/\tau_{\sscaleI}}\circ\underline{c}^{\tau/\tau_{\sscaleI}}}
&&
(\underline{X}, \underline{A})^\sigma  
\ar[d]^{\eta(\sigma)\circ \underline{c}^{\tau/\tau_{\sscaleI}}}\\
(\underline{X}, \underline{A})^{\tau,c}\ar[r]^{p_{\tau,c}^{\tau_{\sscaleI},c}}
& (\underline{X}, \underline{A})^{\tau_{\scaleI},c}
\ar[rr]^{\imath_{\tau_{\sscaleI},c}^{\sigma,c}}
&&(\underline{X}, \underline{A})^{\sigma,c}  
}
\end{equation} 
The left square commutes by precomposing Diagram~(\ref{diagram1}) 
for $\sigma=\tau_{\sscaleI}$ by $c^{\sigma/\tau_{\sscaleI}}$ and the right square commutes
by precomposing the diagram in Lemma~\ref{lemma_construct weighted poly prod}~(vi) by $c^{\tau/\tau_{\sscaleI}}$.  
Since $\imath_{\tau_{\scaleI}}^\sigma\circ p_\tau^{\tau_{\scaleI}}=p_\tau^\sigma$ by Lemma~\ref{classical} and 
$\imath_{\tau_{\scaleI},c}^{\sigma,c}\circ p_{\tau,c}^{\tau_{\scaleI},c}=p_{\tau,c}^{\sigma,c}$ 
from Diagram~(\ref{diagram3}), the outer perimeter of~(\ref{pieta}) is Diagram~(\ref{diagram1}). 
This completes the induction step, thus completing the proof of the proposition. 
\end{proof} 

\begin{remark}\label{3diagramremark} 
A special case of~(\ref{diagram1}) is worth noting. If $\sigma\subseteq\tau$ then 
$\sigma=\tau_{\tinyI}$ for some $I$, so $\underline{c}^{\tau/\tau_{\sscaleI}}=\underline{c}^{\tau/\sigma}$ and 
$\underline{c}^{\sigma/\tau_{\sscaleI}}=\underline{c}^{\sigma/\sigma}=(1,\ldots,1)$, so~(\ref{diagram1}) simplifies to a 
commutative diagram 
\[
\xymatrix{
(\underline{X}, \underline{A})^{\tau}\ar[r]^{p_{\tau}^{\sigma}}
\ar[d]_{\eta(\tau)}
& (\underline{X}, \underline{A})^\sigma   
\ar[d]^{\eta(\sigma)\circ \underline{c}^{\tau/\sigma}}\\
(\underline{X}, \underline{A})^{\tau,c}\ar[r]^{p_{\tau,c}^{\sigma,c}}		
&(\underline{X}, \underline{A})^{\sigma,c}. 
} 
\]
\end{remark} 

The projections also satisfy a naturality property. 

\begin{lemma}\label{lemma-projectA} 
Fix $I\subseteq [m]$. Suppose that $\tau\in\Delta^{m-1}$ and $\delta$ is a face of $\tau$. 
Then the following diagram commutes: 
\begin{equation}
\xymatrix{
(\underline{X}, \underline{A})^{\delta,c}\ar[r]^{p_{\delta,c}^{\delta_{\sscaleI},c}}
\ar[d]_{\imath_{\delta,c}^{\tau,c}}
& (\underline{X}, \underline{A})^{\delta_{\scaleI},c}\
\ar[d]^{\imath_{\delta_{\sscaleI},c}^{\tau_{\sscaleI},c}}\\
(\underline{X}, \underline{A})^{\tau,c}\ar[r]^{p_{\tau,c}^{\tau_{\sscaleI},c}}
& (\underline{X}, \underline{A})^{\tau_{\scaleI},c}
}
\end{equation}
\end{lemma} 

\begin{proof} 
By~(\ref{diagram2}) there is a factorization 
$p_{\delta,c}^{\tau_{\scaleI},c}=p_{\tau,c}^{\tau_{\scaleI},c}\circ \imath_{\delta,c}^{\tau,c}$, 
and by~(\ref{diagram3}) there is a factorization 
$p_{\delta,c}^{\tau_{\scaleI},c}=\imath_{\delta_{\scaleI},c}^{\tau_{\scaleI},c}\circ p_{\delta,c}^{\delta_{\scaleI},c}$. Putting 
this together gives 
$\imath_{\delta_{\scaleI},c}^{\tau_{\scaleI},c}\circ p_{\delta,c}^{\delta_{\scaleI},c}=p_{\tau,c}^{\tau_{\scaleI},c}\circ \imath_{\delta,c}^{\tau,c}$, 
as asserted. 
\end{proof}  

Let $K$ be a simplicial complex on the vertex set $[m]$, let $c$ be a power sequence, 
and let $I\subseteq [m]$. Define 
\[\iota_{K_{I},c}^{K,c}\colon\namedright{(\underline{X},\underline{A})^{K_{I},c}}{}{(\underline{X},\underline{A})^{K,c}}\qquad 
p_{K,c}^{K_{I},c}\colon\namedright{(\underline{X},\underline{A})^{K,c}}{}{(\underline{X},\underline{A})^{K_{I},c}}\] 
as the colimit of the composites 
\[\namedright{(\underline{X},\underline{A})^{\sigma,c}}{\iota_{\sigma,c}^{\tau,c}}{(\underline{X},\underline{A})^{\tau,c}} 
      \rightarrow (\underline{X},\underline{A})^{K,c}\qquad 
\namedright{(\underline{X},\underline{A})^{\tau,c}}{p_{\tau,c}^{\sigma,c}}{(\underline{X},\underline{A})^{\sigma,c}} 
      \rightarrow (\underline{X},\underline{A})^{K_{I},c}\] 
respectively, where the colimit runs over the faces $\tau$ of $K$ with $V(\tau)\cap I\subseteq V(\sigma)$ for 
 faces $\sigma$ of $K_{I}$. 

\begin{proposition}\label{prop-projectLK}
Suppose that $L\subset K$ are simplicial complexes on the vertex set $[m]$. Then the following 
diagram commutes: 
\begin{equation}
\xymatrix{
(\underline{X}, \underline{A})^{L,c}\ar[r]^{p_{L,c}^{L_{\scaleI},c}}
\ar[d]_{\imath_{L,c}^{K,c}}
& (\underline{X}, \underline{A})^{L_I,c}\
\ar[d]^{\imath_{L_{\scaleI},c}^{K_{\scaleI},c}}\\
(\underline{X}, \underline{A})^{K,c}\ar[r]^{p_{K,c}^{K_{\scaleI},c}} 
& (\underline{X}, \underline{A})^{K_I,c}
}
\end{equation}
\end{proposition} 

\begin{proof} 
This follows from Lemma~\ref{lemma-projectA} by taking colimits over the face posets of $L$ and $K$. 
\end{proof} 

As a consequence, we obtain a retraction that is the analogue in the usual polyhedral product case of 
$(\underline{X},\underline{A})^{L}$ retracting off $(\underline{X},\underline{A})^{K}$ when $L$ 
is a full subcomplex of $K$. 

\begin{corollary} 
   \label{weightedprojection1} 
   Let $K$ be a simplicial complex on the vertex set $[m]$ and let $c$ be a power sequence. 
   If $L$ is a full subcomplex of $K$ then there is a map 
   \(p_{K,c}^{L,c}\colon\namedright{(\underline{X},\underline{A})^{K,c}}{}{(\underline{X},\underline{A})^{L,c}}\) 
   that is a left inverse to the inclusion 
   \(\namedright{(\underline{X},\underline{A})^{L,c}}{\imath_{L,c}^{K,c}}{(\underline{X},\underline{A})^{K,c}}\). 
\end{corollary} 

\begin{proof} 
As $L$ is a full subcomplex of $K$, there is a subset $I\subseteq [m]$ such that $L=K_{I}$. 
Since $L_{I}=K_{I}$, the inclusion $\imath_{L_{\scaleI},c}^{K_{\scaleI},c}$ in Proposition~\ref{prop-projectLK} 
is the identity map. Thus if we take $p_{L,c}^{K,c}$ to be $p_{K_{\scaleI},c}^{K,c}$ then the commutative 
diagram in Proposition~\ref{prop-projectLK} implies that $p_{K,c}^{L,c}$ is a left inverse 
of $\imath_{L,c}^{K,c}$.  
\end{proof} 

The retraction in Corollary~\ref{weightedprojection1} can be enhanced to an analogue 
in the usual polyhedral product case of the retraction in~(\ref{KIdgrm}). This first requires 
a lemma. 

Recall that if $L$ is a simplicial complex on a vertex set $I\subseteq [m]$ then $\overline{L}$ 
is the simplicial complex with the same faces as $L$ but regarded as having vertex set $[m]$. 
In other words, the elements in $[m]-I$ are ghost vertices for $\overline{L}$. 

\begin{lemma}\label{chomeo} 
Suppose that $I\subseteq [m]$ and $L$ a simplicial complex on $[m]$ such that $V(L)\subset I$. Then 
there is a homeomorphism 
$$
p_{\overline{L},c}^{L, c}\times \prod_{i\not\in I} p_{\overline{L},c}^{\{i\},c}\colon
(\underline{X},\underline{A})^{\overline{L},c}\rightarrow (\underline{X},\underline{A})^{L,c}
\times \prod_{i\not\in I} A_i. 
$$ 
\end{lemma} 

\begin{proof}
If $L=\emptyset$ then, by Definition~\ref{dfn_weighted poly prod system}~(i), both sides are 
$\prod_{i\in [m]} A_i$ and each of the maps $p_{\overline{L},c}^{L, c}$ and $p_{\overline{L},c}^{\{i\},c}$ 
for $i\notin I$ is the identity. 

Assume the lemma is true when $L$ has fewer than $n$ simplices. 
Now suppose that $L$ has $n$ simplices. There are two cases. 
\medskip 

\noindent 
\textit{Case 1}: If $L$ is not a simplex, consider the composite
$$
(\underline{X},\underline{A})^{\overline{L},c}=\mbox{colim}_{\overline{\tau}\in\overline{L}} (\underline{X},\underline{A})^{\overline{\tau},c}
\stackrel{\Phi}{\rightarrow}
\mbox{colim}_{\overline{\tau}\in\overline{L}} ((\underline{X},\underline{A})^{\overline{\tau}_I,c}\times \prod_{i\not\in I} A_i)
$$
$$
\stackrel{\Theta}{\rightarrow}
(\mbox{colim}_{\tau\in L} (\underline{X},\underline{A})^{\tau,c})\times \prod_{i\not\in I} A_i
=(\underline{X},\underline{A})^{L,c}\times \prod_{i\not\in I} A_i. 
$$
Here $\Phi$ is given by swapping coordinates and $\Theta$ is given by factoring the product $\prod_{i\notin I}A_i$. For $\Phi$, by inductive hypothesis there are homeomorphisms 
$(X,A)^{\overline{\tau},c}\stackrel{\cong}{\rightarrow}
(X,A)^{\overline{\tau}_{\scaleI},c}\times \prod_{i\not\in I} A_i$ 
for each $\overline{\tau}\in\overline{L}$. As $\Phi$ is a colimit of homeomorphisms, it too is a 
homeomorphism. For $\Theta$, in general, if spaces are compactly generated then there is a homeomorphism 
$\mbox{colim}_{j\in\mathcal{J}}(A_{j}\times B)\cong(\mbox{colim}_{j\in\mathcal{J}} A_{j})\times B$. 
Noting that the simplices of $\overline{L}$ are the same as the simplices of $L$, that is, 
$\overline{\tau}_{\tinyI}$ is the same as $\tau$, $\Theta$ is an instance 
of this general fact and so is a homeomorphism. The projection of $\Theta\circ\Phi$ to 
$(\underline{X},\underline{A})^{L,c}$ or to each $A_{i}$ for $i\notin I$ is essentially the definition 
of $p_{\overline{L},c}^{L,c}$ and $p_{\overline{L},c}^{\{i\},c}$ respectively. This therefore completes 
the inductive step in the case when $L$ is not a simplex. 
\medskip 

\noindent 
\textit{Case 2}: If $L=\tau$ (equivalently, $\overline{L}=\overline{\tau}$) is a simplex then consider the cube
\[
\xymatrix{
(\underline{X}, \underline{A})^{\partial\overline{\tau}}\ar[rr]^{\imath_{\partial\overline{\tau}}^{\overline{\tau}}}
\ar[dr]_{\eta(\partial\overline{\tau})} \ar[dd]_(0.3){f_1}
&& (\underline{X}, \underline{A})^{\overline{\tau}}
\ar[dr]^{\eta(\overline{\tau})}\ar@{-}[d]_(0.6){f_2}\\
& (\underline{X}, \underline{A})^{\partial\overline{\tau},c}\ar[rr]^(0.35){\imath_{\partial\overline{\tau},c}^{\overline{\tau},c}}
\ar[dd]_(0.3){f_3}
&\ar[d] & (\underline{X}, \underline{A})^{\overline{\tau},c} \ar[dd]^(0.3){f_4} \\
(\underline{X}, \underline{A})^{\partial \tau}\times \prod_{i\not\in I} A_i
\ar@{-}[r]^(0.65){\imath_{\partial \tau}^\tau\times id}
\ar[dr]_{\eta(\partial \tau)\times id}
&\ar[r] & (\underline{X}, \underline{A})^{\tau}\times \prod_{i\not\in I} A_i
\ar[dr]^{\eta(\tau)\times id}\\
& (\underline{X}, \underline{A})^{\partial \tau,c}\times \prod_{i\not\in I} A_i
\ar[rr]^{\imath_{\partial \tau,c}^{\tau,c}\times id}
&& (\underline{X}, \underline{A})^{\tau,c}\times \prod_{i\not\in I} A_i
}
\]
The top face is a pushout by Definition~\ref{dfn_weighted poly prod system}. The restriction of 
the bottom face away from $\prod_{i\notin I} A_{i}$ is also a pushout by 
Definition~\ref{dfn_weighted poly prod system}, and the product of a pushout with a fixed space  
in the category of compactly generated spaces is a pushout, so the bottom face is a pushout. Each 
map $f_{i}$ for $1\leq i\leq 4$ is of the form $p\times \prod_{i\not\in I} p$. The maps 
$f_{1}$ and $f_{2}$ are homeomorphisms by the usual polyhedral product case while the map $f_{3}$ 
is a homeomorphism by Case 1. The four sides of the cube commute by the naturality of 
the projection, as stated in Proposition~\ref{prop-projectLK}. The map $f_4$ is therefore 
also a map of pushouts, and so as each of $f_{1}$, $f_{2}$ and $f_{3}$ are homeomorphisms, 
so is $f_{4}$. This completes the inductive step in this case. 
\end{proof}

\begin{proposition}  
   \label{KIcprojection} 
   Let $K$ be a simplicial complex on the vertex set $[m]$ and let $c$ be a power sequence. 
   If $I\subseteq [m]$ then there is a commutative diagram  
   \[\xymatrix  
     { 
        (\underline{X},\underline{A})^{K_{I},c}\times\prod_{i\notin I} A_{i}\rto^-{\cong}\drrto_{\mbox{\tiny\rm project}} 
            & (\underline{X},\underline{A})^{\overline{K}_{I},c}\rto^-{\imath_{\overline{K}_{\scaleI},c}^{K,c}}  
            &  (\underline{X},\underline{A})^{K,c}\dto^{p_{K,c}^{K_{\scaleI},c}} \\ 
         & & (\underline{X},\underline{A})^{K_{I},c}  
      }
     \]  
    where $\pi$ is the projection. 
\end{proposition} 

\begin{proof} 
The homeomorphism is due to Lemma~\ref{chomeo}. The commutativity of the diagram 
follows since $p_{K,c}^{K_{\scaleI},c}$ projects away from the coordinates $i\notin I$. 
\end{proof}

\section{A suspension splitting}  
\label{sec:suspsplitting} 

In this section we prove a suspension splitting for weighted polyhedral products of the form 
$(\underline{X},\underline{\ast})^{K,c}$. 
From now on we will consider sequences 
$(\underline{X},\underline{\ast})=\{(X_i,\underline{\ast})\}^m_{i=1}$, 
where each $(X_{i},\ast)$ is a CW-power couple. For instance, by Examples~\ref{ex_monoid} 
and~\ref{ex_suspension}, this would be the case if each $X_{i}$ is a topological monoid or a suspension. 

For $\sigma\in\Delta^{m-1}$ let $\underline{X}^{\wedge\sigma}$ be the 
smash product $\bigwedge_{i\in\sigma}X_i$ if $\sigma\neq\emptyset$ and be the basepoint $\ast$ if $\sigma=\emptyset$.
By definition of the polyhedral product, $(\underline{X},\underline{\ast})^{\sigma}=\prod_{i\in\sigma} X_{i}$ 
and $(\underline{X},\underline{\ast})^{\partial\sigma}$ is the fat wedge in $\prod_{i\in\sigma} X_{i}$. 
Thus there is a cofibration sequence 
\[\nameddright{(\underline{X},\underline{\ast})^{\partial\sigma}}{}{(\underline{X},\underline{\ast})^{\sigma}} 
     {q_{\sigma}}{\underline{X}^{\wedge\sigma}}\] 
where $q_{\sigma}$ is the standard quotient map from the product $\prod_{i\in\sigma} X_{i}$ to the 
smash product $\bigwedge_{i\in\sigma} X_{i}$.

\begin{lemma}
Let $c$ be a power sequence. For $\sigma\in\Delta^{m-1}$ there is a commutative diagram
\begin{equation}\label{diagram_A=* pushout}
\xymatrix{
(\underline{X},\underline{\ast})^{\partial\sigma}\ar[r]^-{\imath_{\partial\sigma}^{\sigma}}\ar[d]^-{\eta(\partial\sigma)}	&(\underline{X},\underline{\ast})^{\sigma}\ar[d]^-{\eta(\sigma)}\ar[r]^-{q_{\sigma}}	&\underline{X}^{\wedge\sigma}\ar@{=}[d]\\
(\underline{X},\underline{\ast})^{\partial\sigma,c}\ar[r]^-{\imath_{\partial\sigma,c}^{\sigma,c}}	&(\underline{X},\underline{\ast})^{\sigma,c}\ar[r]^-{q_{\sigma,c}}	&\underline{X}^{\wedge\sigma}
}
\end{equation}
where $q_{\sigma,c}$ is a quotient map, and the rows are cofibration sequences.
\end{lemma}

\begin{proof}
The left square is the pushout diagram in Property~(iii) of Definition~\ref{dfn_weighted poly prod system}. Hence the cofibers of $\imath^{\sigma}_{\partial\sigma}$ and $\imath^{\sigma,c}_{\partial\sigma,c}$ are homeomorphic, implying the right square.
\end{proof}

Let $K$ be a simplicial complex on the vertex set $[m]$ and fix a power sequence $c$. If $\tau$ 
is a face of~$K$ then $\tau$ is the full subcomplex of $K$ on the vertex set $V(\tau)$. Therefore 
by Proposition~\ref{prop-projectLK} there is a projection 
$\lnamedright{(\underline{X},\underline{\ast})^{K,c}}{p_{K,c}^{\tau,c}}{(\underline{X},\underline{\ast})^{\tau,c}}$. 
Note here that there is no ambiguity in using the notation $(\underline{X},\underline{\ast})^{\tau,c}$ for 
the range as the subspace of each pair $(X_{i},\ast)$ is the basepoint so  
$(\underline{X},\underline{\ast})^{\tau,c}= 
      (\underline{X},\underline{\ast})^{K_{V(\tau)},c}\times\displaystyle\prod_{i\notin V(\tau)}\ast$.  
      
The existence of $p_{K,c}^{\tau,c}$ lets us form the composite
\[\lnameddright{(\underline{X},\underline{\ast})^{K,c}}{p_{K,c}^{\tau,c}}{(\underline{X},\underline{\ast})^{\tau,c}} 
       {q_{\tau,c}}{\underline{X}^{\wedge\tau}}.\] 
In homology, these maps can be added as $\tau$ ranges over the faces of $K$. 
Let $R$ be a commutative ring and let $\zeta(K,c)$ be the composite
\begin{multline*}
\zeta(K,c)\colon\tilde{H}_*((\underline{X},\underline{\ast})^{K,c};R)
\overset{\triangle}{\longrightarrow}
	\underset{\tau\subseteq K}{\bigoplus} \tilde{H}_*((\underline{X},\underline{\ast})^{K,c};R)\\
	    \overset{\underset{\tau\subseteq K}{\bigoplus}(p^{\tau,c}_{K,c})_*}{\llarrow}
	    \underset{\tau\subseteq K}{\bigoplus}\tilde{H}_*((\underline{X},\underline{\ast})^{\tau,c};R)
	    \overset{\underset{\tau\subseteq K}{\bigoplus}(q_{\tau,c})_*}{\llarrow}
	    \underset{\tau\subseteq K}{\bigoplus}\tilde{H}_*(\underline{X}^{\wedge\tau};R).
\end{multline*} 
We will show that $\zeta(K,c)$ is an isomorphism.

\begin{lemma}\label{lemma_zeta nat incl}
Let $K$ be a simplicial complex on $[m]$ and let $c$ be a power sequence. For any subcomplex $L\subseteq K$ there is a commutative diagram
\begin{equation}\label{diagram_zeta nat incl}
\xymatrix{
\tilde{H}_*((\underline{X},\underline{\ast})^{L,c};R)\ar[r]^-{\zeta(L,c)}\ar[d]_-{(\imath^{K,c}_{L,c})_*}
	&\underset{\tau\subseteq L}{\bigoplus}\tilde{H}_*(\underline{X}^{\wedge\tau};R)\ar[d]^-{\text{incl}}\\
\tilde{H}_*((\underline{X},\underline{\ast})^{K,c};R)\ar[r]^-{\zeta(K,c)}
	&\underset{\tau\subseteq K}{\bigoplus}\tilde{H}_*(\underline{X}^{\wedge\tau};R).
}
\end{equation}
\end{lemma}

\begin{proof}
By Proposition~\ref{prop-projectLK} the projection maps $p_{K,c}^{\tau,c}$ are natural for the 
inclusion of a subcomplex $L\subseteq K$. The naturality of the inclusions $\imath_{\partial\sigma}^{\sigma}$ 
implies, from~(\ref{diagram_A=* pushout}), that the quotient maps $q_{\tau,c}$ to the smash product are also natural.  
Thus for any subcomplex $L\subseteq K$, if $\tau$ is a face of $K$ that is also in $L$ then there is a 
commutative diagram 
\[
\xymatrix{
      (\underline{X},\underline{\ast})^{L,c}\ar[r]^{p_{L,c}^{\tau,c}}\ar[d]^{\imath_{L,c}^{K,c}} 
         & (\underline{X},\underline{\ast})^{\tau,c} 
                  \rto^-{q_{\tau,c}}\ar@{=}[d]
         &  \underline{X}^{\wedge\tau}\ar@{=}[d] \\ 
      (\underline{X},\underline{\ast})^{K,c}\ar[r]^{p_{K,c}^{\tau,c}} 
         & (\underline{X},\underline{\ast})^{\tau,c}\rto^-{q_{\tau,c}}
         &  \underline{X}^{\wedge\tau}  
}
\]
and if $\tau$ is a face of $K$ that is not in $L$ then there is a commutative diagram 
\begin{equation}\label{LKnat2} 
\xymatrix{
      (\underline{X},\underline{\ast})^{L,c}\ar[r]^-{p_{L,c}^{L_{\scaleto{\tau}{2pt}},c}}\ar[d]^{\imath_{L,c}^{K,c}} 
         & (\underline{X},\underline{\ast})^{L_{\tau},c}\ar[r]\ar[d]^{\imath_{L_{\scaleto{\tau}{2pt}},c}^{\tau,c}} & \ast\ar[d] \\ 
      (\underline{X},\underline{\ast})^{K,c}\ar[r]^{p_{K,c}^{\tau,c}} 
         & (\underline{X},\underline{\ast})^{\tau,c}\rto^-{q_{\tau,c}}
         &  \underline{X}^{\wedge\tau}  
}
\end{equation} 
where the right square commutes since 
$\tau$ not being in $L$ implies that there is a vertex $i$ of $\tau$ not in $L$ so coordinate $i$ 
in $(\underline{X},\underline{\ast})^{L_{\tau},c}$ is $\ast$, which then quotients trivially to the 
smash product $\underline{X}^{\wedge\tau}$ . 

After taking homology and summing over such diagrams, 
the definition of $\zeta(K,c)$ implies the diagram~\eqref{diagram_zeta nat incl} commutes.
\end{proof}

Recall that a simplex $\sigma$ in $K$ is maximal if it is a maximal element in the face poset of $K$. For a maximal simplex $\sigma\in K$, let $K-\sigma$ be the subposet of $K$ consisting of all elements in $K$ except $\sigma$.

\begin{lemma}\label{lemma_pushout of K}
Let $K$ be a simplicial complex on $[m]$. If $\sigma$ is a maximal simplex in $K$ then there is a commutative diagram
\begin{equation}
\label{diagram_K-sigma pushout}
\xymatrix{
(\underline{X},\underline{\ast})^{\partial\sigma,c}\ar[r]^-{\imath^{\sigma,c}_{\partial\sigma,c}}\ar[d]_-{\imath^{K-\sigma,c}_{\partial\sigma,c}}	& (\underline{X},\underline{\ast})^{\sigma,c}\ar[d]^-{\imath^{K,c}_{\sigma,c}}\ar[rr]^-{q_{\sigma,c}} & & \underline{X}^{\wedge\sigma}\ar@{=}[d]\\
(\underline{X},\underline{\ast})^{K-\sigma,c}\ar[r]^-{\imath^{K,c}_{K-\sigma,c}}	& (\underline{X},\underline{\ast})^{K,c}\ar[rr]^-{q_{\sigma,c}\circ p^{\sigma,c}_{K,c}} & & \underline{X}^{\wedge\sigma}
}
\end{equation}
where the rows are cofibration sequences and the left square is a pushout.
\end{lemma}

\begin{proof}
Since $K$ is the union of $K-\sigma$ and $\sigma$ along $\partial\sigma$, there is a pushout 
\[
\xymatrix{
\mbox{colim}_{\tau\subseteq\partial\sigma}(\underline{X},\underline{\ast})^{\tau,c}\ar[r]\ar[d]	&\mbox{colim}_{\tau\subseteq\sigma}(\underline{X},\underline{\ast})^{\tau,c}\ar[d]\\
\mbox{colim}_{\tau\subseteq K-\sigma}(\underline{X},\underline{\ast})^{\tau,c}\ar[r]		&\mbox{colim}_{\tau\subseteq K}(\underline{X},\underline{\ast})^{\tau,c}.
}
\]
As $(\underline{X},\underline{\ast})^{L,c}=\mbox{colim}_{\tau\subseteq L}(\underline{X},\underline{\ast})^{\tau,c}$ for any simplicial complex $L$ on $[m]$, this pushout is exactly the left square in~(\ref{diagram_K-sigma pushout}). The right square in~(\ref{diagram_K-sigma pushout}) commutes since $p^{\sigma,c}_{K,c}\circ\imath^{K,c}_{\sigma,c}$ is the identity map on $(\underline{X},\underline{\ast})^{\sigma,c}$ by Corollary~\ref{weightedprojection1}. 

Since $\imath^{\sigma,c}_{\partial\sigma,c}$ is a cofibration, so is $\imath^{K,c}_{K-\sigma,c}$ and their cofibers are homeomorphic. The top row of~\eqref{diagram_K-sigma pushout} is the cofibration sequence given by the bottom row of~\eqref{diagram_A=* pushout}. As the left square in~(\ref{diagram_K-sigma pushout}) is a pushout, this implies that there is a cofibration 
\(\llnameddright{(\underline{X},\underline{\ast})^{K-\sigma,c}}{\imath_{K-\sigma,c}^{K,c}}{(\underline{X},\underline{\ast})^{K,c}}{\gamma} 
   {\underline{X}^{\wedge\sigma}}\) 
for some map $\gamma$ with the property that $\gamma\circ\imath^{K,c}_{\sigma,c}=q_{\sigma,c}$. Observe that the composite $q_{\sigma,c}\circ p^{\sigma,c}_{K,c}\circ\imath^{K,c}_{K-\sigma,c}$ is the trivial map by the commutativity of~(\ref{LKnat2}). Thus 
$q_{\sigma,c}\circ p_{K,c}^{\sigma,c}$ extends across $\gamma$ to a map 
\(e\colon\namedright{\underline{X}^{\wedge\sigma}}{}{\underline{X}^{\wedge\sigma}}\), 
that is, $e\circ\gamma=q_{\sigma,c}\circ p_{K,c}^{\sigma,c}$.  
Since $p_{K,c}^{\sigma,c}\circ\imath^{K,c}_{\sigma,c}$ is the identity map on $(\underline{X},\underline{\ast})^{\sigma,c}$, we obtain $e\circ\gamma\circ\imath^{K,c}_{\sigma,c}=q_{\sigma,c}\circ p_{K,c}^{\sigma,c}\circ\imath^{K,c}_{\sigma,c}=q_{\sigma,c}$. 
But $\gamma\circ\imath^{K,c}_{\sigma,c}=q_{\sigma,c}$, so $e\circ q_{\sigma,c}=q_{\sigma,c}$. As $q_{\sigma,c}$ is an epimorphism, this implies that $e$ is the identity map. Thus the bottom row of~\eqref{diagram_K-sigma pushout} is a cofibration sequence.
\end{proof}

\begin{lemma}\label{lemma_epsilon iso}
The map $\zeta(K,c)$ is an isomorphism of graded $R$-modules.
\end{lemma}

\begin{proof}
The proof is by induction on the number of simplices in $K$. Let $P(n)$ be the statement $\zeta(K,c)$ is an isomorphism for a simplicial complex $K$ having at most $n$ simplices''. When $K=\emptyset$ notice that $\tilde{H}_*((\underline{X},\underline{\ast})^{K,c};R)$ and $\bigoplus_{\tau\subseteq K}\tilde{H}_*(\underline{X}^{\wedge\tau};R)$ are trivial modules. So $\zeta(\emptyset,c)$ is an isomorphism and $P(1)$ is true.

Assume that $P(i)$ is true for $i<n$. Let $K$ be a simplicial complex with $n$ simplices. Let $\sigma$ be a maximal simplex of $K$. Suppose that there is a commutative diagram
\begin{equation}\label{diagram_e(sigma) equiv inductive pf}
\spreaddiagramcolumns{-0.5pc}
\xymatrix{
0\ar[r]
	&\tilde{H}_*((\underline{X},\underline{\ast})^{K-\sigma,c};R)\ar[rr]^-{(\imath^{K,c}_{K-\sigma,c})_*}\ar[d]^-{\zeta(K-\sigma)}
	&&\tilde{H}_*((\underline{X},\underline{\ast})^{K,c};R)\ar[rr]^-{(q_{\sigma,c}\circ p^{\sigma,c}_{K,c})_*}\ar[d]^-{\zeta(K,c)}
	&&\tilde{H}_*(\underline{X}^{\wedge\sigma};R)\ar[r]\ar@{=}[d]
	&0\\
0\ar[r]
	&\underset{\tau\subseteq K-\sigma}{\bigoplus}\tilde{H}_*(\underline{X}^{\wedge\tau};R)\ar[rr]^-{\text{incl}}
	&&\underset{\tau\subseteq K}{\bigoplus}\tilde{H}_*(\underline{X}^{\wedge\tau};R)\ar[rr]^-{\text{project}}
	&&\tilde{H}_*(\underline{X}^{\wedge\sigma};R)\ar[r]
	&0
}
\end{equation}
where the rows are exact sequences. Then the inductive hypothesis implies that $\zeta(K-\sigma,c)$ is an isomorphism and the Five Lemma therefore implies that $\zeta(K,c)$ is an isomorphism. Thus $P(n)$ is true and the induction is complete.

It remains to show the existence of Diagram~\eqref{diagram_e(sigma) equiv inductive pf}. The left square commutes by~\eqref{diagram_zeta nat incl}. The right square commutes by the definition of $\zeta(K,c)$. The bottom row is clearly an exact sequence. Thus it remains to show that the top row is an exact sequence.

By Lemma~\ref{lemma_pushout of K}, the top row of~\eqref{diagram_K-sigma pushout} is a cofibration sequence. Taking homology therefore implies that there is a long exact sequence
\[
\cdots\to
	H_{n+1}(\underline{X}^{\wedge\sigma};R)\longrightarrow
	H_n((\underline{X},\underline{\ast})^{K-\sigma,c};R)\overset{(\imath^{K,c}_{K-\sigma,c})_*}{\llarrow}
	H_n((\underline{X},\underline{\ast})^{K,c};R)
\overset{(q_{\sigma,c}\circ p^{\sigma,c}_{K,c})_*}{\llarrow}
	H_n(\underline{X}^{\wedge\sigma};R)
	\to\cdots
\]
By~\cite{BBCG1}, the quotient map $q_{\sigma}$ in~\eqref{diagram_A=* pushout} has the property 
that $\Sigma q_{\sigma}$ has a right homotopy inverse. The commutativity of the right square 
in~\eqref{diagram_A=* pushout} therefore implies that $\Sigma q_{\sigma,c}$ has a right homotopy 
inverse. By Corollary~\ref{weightedprojection1}, the map $p_{K,c}^{\sigma,c}$ has 
a right inverse. Hence $\Sigma(q_{\sigma,c}\circ p^{\sigma,c}_{K,c})$ has a right homotopy inverse. 
Therefore $(q_{\sigma,c}\circ p^{\sigma,c}_{K,c})_*$ is surjective and the long exact sequence in 
homology breaks into short exact sequences, resulting in the top row 
of~\eqref{diagram_e(sigma) equiv inductive pf} being an exact sequence.
\end{proof}

The isomorphism in Lemma~\ref{lemma_epsilon iso} can be geometrically realized by a suspension splitting. Let 
\(s\colon\namedright{\Sigma(\underline{X},\underline{\ast})^{K,c}}{}{\bigvee_{\tau\subseteq K}\Sigma(\underline{X},\underline{\ast})^{K,c}}\) 
be an iterated comultiplication.  

\begin{theorem}\label{suspXS}
For any simplicial complex $K$ on $[m]$, the composite
\[
e(K,c):\Sigma(\underline{X},\underline{\ast})^{K,c}\overset{s}{\longrightarrow}\bigvee_{\tau\subseteq K}\Sigma(\underline{X},\underline{\ast})^{K,c}\overset{\bigvee_{\tau\subseteq K}\Sigma p^{\tau,c}_{K,c}}{\llllarrow}\bigvee_{\tau\subseteq K}\Sigma(\underline{X},\underline{\ast})^{\tau,c}\overset{\bigvee_{\tau\subseteq K}\Sigma q_{\tau,c}}{\llllarrow}\bigvee_{\tau\subseteq K}\Sigma\underline{X}^{\wedge\tau}
\]
is a homotopy equivalence.
\end{theorem}

\begin{proof}
Let $R$ be a commutative ring. By its definition, $\zeta(K,c)$ equals the composite
\[
\tilde{H}_{*-1}((\underline{X},\underline{\ast})^{K,c};R)
	\overset{\cong}{\longrightarrow}
\tilde{H}_*(\Sigma(\underline{X},\underline{\ast})^{K,c};R)
	\overset{e(K,c)_*}{\llarrow}
\bigoplus_{\tau\subseteq K}\tilde{H}_*(\Sigma\underline{X}^{\wedge\sigma};R)
	\overset{\cong}{\longrightarrow}
\bigoplus_{\tau\subseteq K}\tilde{H}_{*-1}(\underline{X}^{\wedge\sigma};R)
\]
where the first and the third maps are the suspension isomorphism. By Lemma~\ref{lemma_epsilon iso}, $\zeta(K,c)$ is an isomorphism of graded $R$-modules. Thus $e(K,c)_{\ast}$ induces an isomorphism in integral homology. Since $\Sigma(\underline{X},\underline{\ast})^{K,c}$ and $\bigvee_{\tau\subseteq K}\Sigma\underline{X}^{\wedge\tau}$ are simply connected CW-complexes, Whitehead's Theorem implies that $e(K,c)$ is a homotopy equivalence.
\end{proof} 

Theorem~\ref{suspXS} recovers the suspension splitting in~\cite{BBCG1} for polyhedral 
products of the form $(\underline{X},\underline{\ast})^{K}$ when each $X_{i}$ is a CW-complex. By Example~\ref{polyprod_reduction}, when $c(\sigma)=(1,\ldots,1)$ then 
$(\underline{X},\underline{\ast})^{K,c}=(\underline{X},\underline{\ast})^{K}$ and the map 
\(\namedright{(\underline{X},\underline{\ast})^{K}}{\eta(K)}{(\underline{X},\underline{\ast})^{K,c}}\) 
is the identity map. Therefore, if $e(K)$ is the composite
\[
e(K):\Sigma(\underline{X},\underline{\ast})^{K}\overset{s}{\longrightarrow}\bigvee_{\tau\subseteq K}\Sigma(\underline{X},\underline{\ast})^{K}\overset{\bigvee_{\tau\subseteq K}\Sigma p^{\tau}_K}{\llarrow}\bigvee_{\tau\subseteq K}\Sigma(\underline{X},\underline{\ast})^{\tau}\overset{\bigvee_{\tau\subseteq K}\Sigma q_{\tau}}{\llarrow}\bigvee_{\tau\subseteq K}\Sigma\underline{X}^{\wedge\tau}
\]
then it is a homotopy equivalence.

Dually working with cohomology, let $\epsilon(K,c)$ be the composite
\[
\epsilon(K,c)\colon\underset{\tau\subseteq K}{\bigoplus}\tilde{H}^*(\underline{X}^{\wedge\tau};R)
\overset{\underset{\tau\subseteq K}{\bigoplus}(q_{\tau,c}\circ p^{\tau,c}_{\sigma,c})^*}{\llarrow}
     \underset{\tau\subseteq K}{\bigoplus} \tilde{H}^*((\underline{X},\underline{\ast})^{K,c};R)
     \overset{\text{sum}}{\longrightarrow}\tilde{H}^*((\underline{X},\underline{\ast})^{K,c};R)
\]
and let $\epsilon(K)\colon\underset{\tau\subseteq K}{\bigoplus}\tilde{H}^*(\underline{X}^{\wedge\tau};R)\to\tilde{H}^*((\underline{X},\underline{\ast})^{K};R)$ be defined similarly.

\begin{corollary} \label{cor_epsilon iso}
   \label{suspXScor} 
   The maps $\epsilon(K,c)$ and $\epsilon(K)$ are isomorphisms of graded $R$-modules. In addition, for any subcomplex $L\subseteq K$ there is a commutative diagram
\begin{equation}\label{diagram_epsilon nat incl}
\xymatrix{
\underset{\tau\subseteq K}{\bigoplus}\tilde{H}^*(\underline{X}^{\wedge\tau};R)\ar[r]^-{\epsilon(K,c)}\ar[d]_-{\text{project}}	&\tilde{H}^*((\underline{X},\underline{\ast})^{K,c};R)\ar[d]^-{(\imath^K_L)^*}\\
\underset{\tau\subseteq L}{\bigoplus}\tilde{H}^*(\underline{X}^{\wedge\tau};R)\ar[r]^-{\epsilon(L,c)}	&\tilde{H}^*((\underline{X},\underline{\ast})^{L,c};R). 
}
\end{equation}
\end{corollary}

\begin{proof}
Since $\epsilon(K,c)$ equals the composite
\[
\bigoplus_{\tau\subseteq K}\tilde{H}^*(\underline{X}^{\wedge\sigma};R)\overset{\cong}{\longrightarrow}\bigoplus_{\tau\subseteq K}\tilde{H}^{*+1}(\Sigma\underline{X}^{\wedge\sigma};R)\overset{e(K,c)^*}{\llarrow}\tilde{H}^{*+1}(\Sigma(\underline{X},\underline{\ast})^{K,c};R)\overset{\cong}{\longrightarrow}\tilde{H}^*((\underline{X},\underline{\ast})^{K,c};R),
\]
it is an isomorphism of graded $R$-modules. Similarly so is $\epsilon(K)$.

Modifying the proof of Lemma~\ref{lemma_zeta nat incl} shows that the diagram~\eqref{diagram_epsilon nat incl} commutes.
\end{proof}

\section{Cohomology}\label{coho_compute}

In this section we calculate the cohomology of weighted polyhedral products of the form 
$(\underline{X},\underline{\ast})^{K,c}$ when all the spaces $X_i$ for $1\leq i\leq m$  
are suspensions and the power sequence $c$ arises from the power maps on $X_{i}$ 
induced by the suspension structure. Let $R$ be a subring of $\mathbb{Q}$; cohomology 
will be taken with coefficients in $R$. Throughout we will assume that $H^{\ast}(X_{i};R)$ 
for $1\leq i\leq m$ is a free graded $R$-module of finite type so that the K\"{u}nneth Isomorphism 
holds. 

We begin by recalling a description of the cohomology ring of 
$(\underline{X},\underline{\ast})^{\Delta^{m-1}}$ in~\cite{BBCG2}. For 
$\sigma=(i_{1},\ldots,i_{k})\in\Delta^{m-1}$, recall that 
$(\underline{X},\underline{\ast})=X_{i_{1}}\times\cdots\times X_{i_{k}}$  
and $\underline{X}^{\wedge\sigma}=X_{i_{1}}\wedge\cdots\wedge X_{i_{k}}$. Let 
\[d_{\sigma}\colon\namedright{(\underline{X},\underline{\ast})^{\sigma}}{} 
      {(\underline{X},\underline{\ast})^{\sigma}\wedge (\underline{X},\underline{\ast})^{\sigma}}\] 
\[\widehat{d}_{\sigma}\colon\namedright{\underline{X}^{\wedge\sigma}}{} 
     {\underline{X}^{\wedge\sigma}\wedge\underline{X}^{\wedge\sigma}}\] 
be the reduced diagonals. (These are denoted by $\Delta_{I}$ and $\widehat{\Delta}_{I}$ 
respectively in~\cite{BBCG2}, where $I=\{i_{1},\ldots,i_{k}\}$.)
Let $r_{\Delta^{m-1}}^{\sigma}$ be the composite 
\[r_{\Delta^{m-1}}^{\sigma}\colon\llnameddright{(\underline{X},\underline{\ast})^{\Delta^{m-1}}} 
     {p_{\Delta^{m-1}}^{\sigma}}{(\underline{X},\underline{\ast})^{\sigma}}{q_{\sigma}}{\underline{X}^{\wedge\sigma}},\]
where it may be helpful to recall that $p_{\Delta^{m-1}}^{\sigma}$ is the projection 
onto the full subcomplex and $q_{\sigma}$ is the quotient map.
If $\sigma$ is the disjoint union of $\tau$ and $\omega$,$\tau\sqcup\omega=\sigma$, 
then there are partial diagonals 
\[\widehat{d}_{\sigma}^{\tau,\omega}\colon\namedright{\underline{X}^{\wedge\sigma}} 
     {}{\underline{X}^{\wedge\tau}\wedge\underline{X}^{\wedge\omega}}\] 
that satisfy a commutative diagram 
\begin{equation} 
  \label{deltadgrm1} 
  \diagram 
      (\underline{X},\underline{\ast})^{\Delta^{m-1}}\rrto^-{d_{\Delta^{m-1}}}\dto^{r_{\Delta^{m-1}}^{\sigma}} 
          & & (\underline{X},\underline{\ast})^{\Delta^{m-1}}\wedge (\underline{X},\underline{\ast})^{\Delta^{m-1}} 
                 \dto^{r_{\Delta^{m-1}}^{\tau}\wedge r_{\Delta^{m-1}}^{\omega}} \\ 
      \underline{X}^{\wedge\sigma}\rrto^-{\widehat{d}_{\sigma}^{\tau,\omega}} 
          & & \underline{X}^{\wedge\tau}\wedge\underline{X}^{\wedge\omega}. 
  \enddiagram 
\end{equation}  

Given cohomology classes $u\in H^{p}(\underline{X}^{\wedge\tau};R)$ and 
$v\in H^{q}(\underline{X}^{\wedge\omega};R)$, we obtain the class 
$u\otimes v\in H^{p}(\underline{X}^{\wedge\tau};R)\otimes H^{q}(\underline{X}^{\wedge\omega};R)$. 
The image of $u\otimes v$ in $H^{p+q}(\underline{X}^{\wedge\tau}\wedge\underline{X}^{\wedge\omega}  ;R )$ 
under the Kunneth isomorphism will also be denoted by $u\otimes v$. Define the \emph{star product} 
$u\ast v\in H^{p+q}(\underline{X}^{\wedge\sigma};R)$ by  
\[u\ast v=(\widehat{d}_{\sigma}^{\tau,\omega})^{\ast}(u\otimes v).\] 
Observe that the commutativity of~(\ref{deltadgrm1}) implies that 
\[(r^{\sigma}_{\Delta^{m-1}})^{\ast}(u\ast v)= 
      (r^{\tau}_{\Delta^{m-1}})^{\ast}(u)\cup ( r^{\omega}_{\Delta^{m-1}})^{\ast}(v)\] 
where $\cup$ is the cup product in $H^{\ast}((\underline{X},\underline{\ast})^{\Delta^{m-1}};R)$.

Let  
\[\mathcal{H}^{\ast}((\underline{X},\underline{\ast})^{\Delta^{m-1}};R)= 
       \bigoplus_{\sigma\in\Delta^{m-1}}\tilde{H}^{\ast}(\underline{X}^{\wedge\sigma};R)\oplus R\] 
and define a ring structure on $\mathcal{H}^{\ast}((\underline{X},\underline{\ast})^{\Delta^{m-1}};R)$ 
by the star product. Define 
\[h\colon\namedright{\mathcal{H}^{\ast}((\underline{X},\underline{\ast})^{\Delta^{m-1}};R)} 
      {}{H^{\ast}((\underline{X},\underline{\ast})^{\Delta^{m-1}};R)}\] 
by requiring that the restriction of $h$ to $H^{\ast}(\underline{X}^{\wedge\sigma};R)$ is 
$(r_{\Delta^{m-1}}^{\sigma})^{\ast}$. In the notation of Corollary~\ref{suspXScor}, 
$h=\epsilon(\Delta^{m-1})$ and Corollary~\ref{suspXScor} states that $h$ 
is an additive isomorphism. A multiplicative result was proved in~\cite[Theorem 1.4]{BBCG2}. 

\begin{theorem} 
   \label{BBCGstar} 
   The map $h$ is a ring isomorphism.~$\qqed$ 
\end{theorem} 

In fact, a much more general result was proved in~\cite[Theorem 1.4]{BBCG2}, 
where the star product was defined and shown to describe the cup product 
structure  in $H^{\ast}((\underline{X},\underline{\ast});R)$ for any polyhedral 
product. However, we only need the stated special case. 

\begin{remark} 
\label{starremark} 
In the case when each $X_{i}$ is a suspension it is shown in~\cite[Theorem 1.6]{BBCG2} 
that if $\tau\cup\omega=\sigma$ but $\tau\cap\omega\neq\emptyset$ then $u\ast v=0$. Thus, by 
Theorem~\ref{BBCGstar}, if each $X_{i}$ is a suspension then nontrivial cup products exist in 
$\mathcal{H}^{\ast}((\underline{X},\underline{\ast})^{\Delta^{m-1}};R)\cong 
H^{\ast}((\underline{X},\underline{\ast})^{\Delta^{m-1}};R)$ 
only when $\tau\cup\omega=\sigma$ and $\tau\cap\omega=\emptyset$. 
\end{remark} 

Recall that it is assumed that the underlying graded $R$-module of $H^*(X_i;R)$ is free and of finite type 
for $1\leq i\leq m$. Define a basis for $H^{\ast}((\underline{X},\underline{\ast})^{\Delta^{m-1}};R)$ as follows. 
Let $H^{\ast}(X_{i};R)$ have basis $\{x'_{i,j}\mid j\in J_i\}$ for a countable index set $J_i$. 
Define $x_{i,j}\in H^{\ast}((\underline{X},\underline{\ast})^{\Delta^{m-1}};R)$ by 
\[x_{i,j}=h(x'_{i,j}).\] 
If $\tau=(i_1,\ldots,i_k)\in\Delta^{m-1}$, let $J_{\tau}=J_{i_1}\times\cdots\times J_{i_k}$. 
For $\mathfrak{u}=(j_1,\ldots,j_k)\in J_{\tau}$, define $x_{\tau,\mathfrak{u}}\in H^{\ast}((\underline{X},\underline{\ast})^{\Delta^{m-1}};R)$ by 
\[x_{\tau,\mathfrak{u}}=x_{i_{1},j_{1}}\cup\cdots\cup x_{i_{k},j_{k}}.\] 
Equivalently, since $h$ is a ring homomorphism, in terms of the star product we have 
\[x_{\tau,\mathfrak{u}}=h(x'_{i_{1},j_{1}}\ast\cdots\ast x'_{i_{k},j_{k}}).\] 
As a graded $R$-module the image of 
\(\llnamedright{\tilde{H}^{\ast}(\underline{X}^{\wedge\tau};R)}{(r_{\Delta^{m-1}}^{\tau})^{\ast}} 
      {\tilde{H}^*((\underline{X},\underline{\ast})^{\Delta^{m-1}};R)}\) 
is spanned by \mbox{$\{x_{\tau,\mathfrak{u}}\mid \mathfrak{u}\in J_{\tau}\}$}, so Theorem~\ref{BBCGstar} implies 
that $\tilde{H}^*((\underline{X},\underline{\ast})^{\Delta^{m-1}};R)$ is spanned by 
$\{x_{\tau,\mathfrak{u}}|\tau\in\Delta^{m-1},\mathfrak{u}\in J_{\tau}\}$.

\begin{remark} 
\label{BBCGstarsigmaremark} 
If $\sigma\in\Delta^{m-1}$ then $\sigma=\Delta^{m'-1}$ for some $m'\leq m$. The map $h$ in 
this case identifies with the map $\epsilon(\sigma)$ in Corollary~\ref{suspXScor} and takes the form  
 \(\epsilon(\sigma)\colon\namedright{\mathcal{H}^{\ast}((\underline{X},\underline{\ast})^{\sigma};R)}  
          {}{H^{\ast}((\underline{X},\underline{\ast})^{\sigma};R)}\).  
It is a ring isomorphism and as a graded $R$-module 
$\tilde{H}^*((\underline{X},\underline{\ast})^{\sigma};R)$ is spanned by 
$\{x_{\tau,\mathfrak{u}}|\tau\subseteq\sigma,\mathfrak{u}\in J_{\tau}\}$.~$\qqed$ 
\end{remark}

%

We now define a basis for $H^{\ast}((\underline{X},\underline{\ast})^{\Delta^{m-1},c};R)$. The definition is somewhat 
parallel to that for the basis of $H^{\ast}((\underline{X},\underline{\ast})^{\Delta^{m-1}};R)$, which was determined by 
$h=\epsilon(\Delta^{m-1})$ and the isomorphism in Corollary~\ref{suspXScor}. Again, $H^{\ast}(X_{i};R)$ has 
basis $\{x'_{i,j}\mid j\in J_i\}$ for a countable index set $J_i$. Let $c$ be a power sequence and recall 
the graded $R$-module isomorphism $\epsilon(\Delta^{m-1},c)$ defined above Lemma~\ref{cor_epsilon iso}. 
Define $y_{i,j}\in H^{\ast}((\underline{X},\underline{\ast})^{\Delta^{m-1},c};R)$ by 
\[y_{i,j}=\epsilon(\Delta^{m-1},c)(x'_{i,j}).\] 
If $\tau=(i_1,\ldots,i_k)\in\Delta^{m-1}$, let $J_{\tau}=J_{i_1}\times\cdots\times J_{i_k}$. 
For $\mathfrak{u}=(j_1,\ldots,j_k)\in J_{\tau}$, define $y_{\tau,\mathfrak{u}}\in H^{\ast}((\underline{X},\underline{\ast})^{\Delta^{m-1},c};R)$ by 
\[y_{\tau,\mathfrak{u}}=\epsilon(\Delta^{m-1},c)(x'_{i_{1},i_{1}}\otimes\cdots\otimes x'_{i_{k},j_{k}}).\] 
Since $\epsilon(\Delta^{m-1},c)$ is an isomorphism of graded $R$-modules,
$H^{\ast}((\underline{X},\underline{\ast})^{\Delta^{m-1},c};R)$ is 
spanned by the elements $\{y_{\tau,\mathfrak{u}}\mid\tau\in\Delta^{m-1},\mathfrak{u}\in J_{\tau}\}$. 

More generally, for any $\sigma\in\Delta^{m-1}$, the map $\epsilon(\sigma,c)$ 
gives $H^{\ast}((\underline{X},\underline{\ast})^{\sigma,c};R)$ a graded $R$-module basis  
spanned by the elements $\{y_{\tau,\mathfrak{u}}\mid\tau\subseteq\sigma,\mathfrak{u}\in J_{\tau}\}$.

Notice at this point that the definition of $y_{i,j}$ exactly parallels the definition of $x_{i,j}$, but 
the elements $y_{\tau,\mathfrak{u}}$ are defined in an additive manner while the elements $x_{\tau,\mathfrak{u}}$ 
are defined by the cup product. The elements $x_{\tau,\mathfrak{u}}$ could also have been defined 
in an additive manner, but the fact that $h=\epsilon(\Delta^{m-1})$ (or $\epsilon(\sigma)$ in the case 
when $\tau\subseteq\sigma$) is a ring homomorphism implies that 
the two descriptions are equivalent. The goal is now to determine the cup product structure 
on the basis elements $y_{\tau,\mathfrak{u}}$. To do this we first describe the effect in cohomology of the map 
\(\namedright{(\underline{X},\underline{\ast})^{\sigma}}{\eta(\sigma)}{(\underline{X},\underline{\ast})^{\sigma,c}}\). 

\begin{definition} 
\label{c^tau/sigma}
For any $\tau\subseteq\sigma$ let $\underline{c}^{\wedge\sigma/\tau}:\underline{X}^{\wedge\tau}\to\underline{X}^{\wedge\tau}$ be the wedge of power maps
\[
\underline{c}^{\wedge\sigma/\tau}=\bigwedge_{i\in\tau}\left(\frac{c^{\sigma}_i}{c^{\tau}_i}\right):\bigwedge_{i\in\tau}X_i\to\bigwedge_{i\in\tau}X_i.
\] 
\end{definition}

\begin{proposition}\label{prop_(sigma) and eta(sigma)}
For any $\tau\subseteq\sigma$, there is a commutative diagram
\[
\xymatrix{
\underset{\tau\subseteq\sigma}{\bigoplus}\tilde{H}^*(\underline{X}^{\wedge\tau};R)\ar[r]^-{\epsilon(\sigma,c)}\ar[d]_-{\underset{\tau\subseteq\sigma}{\bigoplus}(\underline{c}^{\wedge\sigma/\tau})^*}	&\tilde{H}^*((\underline{X},\underline{\ast})^{\sigma,c};R)\ar[d]^-{\eta(\sigma)^*}\\
\underset{\tau\subseteq\sigma}{\bigoplus}\tilde{H}^*(\underline{X}^{\wedge\tau};R)\ar[r]^-{\epsilon(\sigma)}	&\tilde{H}^*((\underline{X},\underline{\ast})^{\sigma};R)
}
\]
\end{proposition}

\begin{proof}
By their definitions, $\epsilon(\sigma)$ and $\epsilon(\sigma,c)$ are the sums of $(q_{\tau}\circ p^{\tau}_{\sigma})^*$ and $(q_{\tau,c}\circ p^{\tau,c}_{\sigma,c})^*$ for all simplices $\tau\subseteq\sigma$. Therefore it suffices to show that for each $\tau\subseteq\sigma$ there is a commutative diagram
\begin{equation}\label{diagram_local eta and epsilon}
\xymatrix{
\tilde{H}^*(\underline{X}^{\wedge\tau};R)\ar[rr]^-{(q_{\tau,c}\circ p^{\tau,c}_{\sigma,c})^*}\ar[d]_-{(\underline{c}^{\wedge\sigma/\tau})^*} & &\tilde{H}^*((\underline{X},\underline{\ast})^{\sigma,c};R)\ar[d]^-{\eta(\sigma)^*}\\
\tilde{H}^*(\underline{X}^{\wedge\tau};R)\ar[rr]^-{(q_{\tau}\circ p^{\tau}_{\sigma})^*} & &\tilde{H}^*((\underline{X},\underline{\ast})^{\sigma};R).
}
\end{equation}

For $\tau\subseteq\sigma$ we have $(\underline{X},\underline{\ast})^{\tau}=\prod_{i\in\tau}X_i$ and $\underline{X}^{\wedge\tau}=\bigwedge_{i\in\tau}X_i$. So $\underline{c}^{\sigma/\tau}$ applied to $(\underline{X},\underline{\ast})^{\tau}$ is the map
\[
\prod_{i\in\tau}\frac{c^{\sigma}_i}{c^{\tau}_i}\colon\prod_{i\in\tau}X_i\llarrow\prod_{i\in\tau}X_i,
\]
and $\underline{c}^{\wedge\sigma/\tau}$ applied to $\underline{X}^{\wedge\tau}$ is the map
\[
\bigwedge_{i\in\tau}\frac{c^{\sigma}_i}{c^{\tau}_i}\colon\bigwedge_{i\in\tau}X_i\llarrow\bigwedge_{i\in\tau}X_i.
\]
Thus there is a commutative diagram
\begin{equation}\label{diagram_quotient map and power map}
\xymatrix{
(\underline{X},\underline{\ast})^{\tau}\ar[r]^-{q_{\tau}}\ar[d]_-{\underline{c}^{\sigma/\tau}}	&\underline{X}^{\wedge\tau}\ar[d]^-{\underline{c}^{\wedge\sigma/\tau}}\\
(\underline{X},\underline{\ast})^{\tau}\ar[r]^-{q_{\tau}}	&\underline{X}^{\wedge\tau}.
}
\end{equation}

To prove the commutativity of~\eqref{diagram_local eta and epsilon}, consider the diagram
\[
\xymatrix{
(\underline{X},\underline{\ast})^{\sigma}\ar[r]^-{p^{\tau}_{\sigma}}\ar[dd]_-{\eta(\sigma)}	&(\underline{X},\underline{\ast})^{\tau}\ar[r]^-{q_{\tau}}\ar[d]_-{\underline{c}^{\sigma/\tau}}	&\underline{X}^{\wedge\tau}\ar[d]^-{\underline{c}^{\wedge\sigma/\tau}}\\
	&(\underline{X},\underline{\ast})^{\tau}\ar[r]^-{q_{\tau}}\ar[d]_-{\eta(\tau)}	&\underline{X}^{\wedge\tau}\ar@{=}[d]\\
(\underline{X},\underline{\ast})^{\sigma,c}\ar[r]^-{p^{\tau,c}_{\sigma,c}}	&(\underline{X},\underline{\ast})^{\tau,c}\ar[r]^-{q_{\tau,c}}	&\underline{X}^{\wedge\tau}
}
\]
The left rectangle commutes by Remark~\ref{3diagramremark}, the top right square commutes by~\eqref{diagram_quotient map and power map}, and the bottom right square commutes by the right square in~\eqref{diagram_A=* pushout}. Hence the outer rectangle commutes. Taking reduced cohomology then gives~\eqref{diagram_local eta and epsilon}, as required.  
\end{proof} 

Proposition~\ref{prop_(sigma) and eta(sigma)} has two useful corollaries. 
 
\begin{corollary} 
   \label{xycup} 
   The map 
   \(\namedright{(\underline{X},\underline{\ast})^{\sigma}}{\eta(\sigma)}{(\underline{X},\underline{\ast})^{\sigma,c}}\) 
   relates the bases for $H^{\ast}((\underline{X},\underline{\ast})^{\sigma,c};R)$ and  
   $H^{\ast}((\underline{X},\underline{\ast})^{\sigma};R)$ as follows. If 
    $\tau=(i_{1},\ldots,i_{k})$ and $\mathfrak{u}\in J_{\tau}$ then
   \[
      \eta(\sigma)^*(y_{\tau,\mathfrak{u}})= 
           \left(\prod^k_{\ell=1}\frac{c_{i_{\ell}}^{\sigma}}{c_{i_{\ell}}^{\tau}}\right)x_{\tau,\mathfrak{u}}. 
   \] 
\end{corollary} 

\begin{proof} 
By definition, $y_{\tau,\mathfrak{u}}=\epsilon(\sigma,c)(x'_{i_{1},j_{1}}\otimes\cdots\otimes x'_{i_{k},j_{k}})$. Therefore, 
Proposition~\ref{prop_(sigma) and eta(sigma)} implies that 
\[\eta(\sigma)^{\ast}(y_{\tau,\mathfrak{u}})= 
    \left(\prod^k_{\ell=1}\frac{c_{i_{\ell}}^{\sigma}}{c_{i_{\ell}}^{\tau}}\right) 
        \epsilon(\sigma)(x'_{i_{1},j_{1}}\otimes\cdots\otimes x'_{i_{k},j_{k}}).\] 
On the other hand, as $\epsilon(\sigma)$ is a ring homomorphism by 
Remark~\ref{BBCGstarsigmaremark}, we have 
\[\epsilon(\sigma)(x'_{i_{1},j_{1}}\otimes\cdots\otimes x'_{i_{k},j_{k}})=
      \epsilon(\sigma)(x'_{i_{1},j_{1}})\cup\cdots\cup\epsilon(\sigma)(x'_{i_{k},j_{k}}).\] 
By definition, $x_{i,j}=\epsilon(\sigma)(x'_{i,j})$ and $x_{\tau,\mathfrak{u}}=x_{i_{1},j_{1}}\cup\cdots\cup x_{i_{k},j_{k}}$ for $\tau\subseteq\sigma$. Thus 
$\eta(\sigma)^*(y_{\tau,\mathfrak{u}})= 
           \left(\prod^k_{\ell=1}\frac{c_{i_{\ell}}^{\sigma}}{c_{i_{\ell}}^{\tau}}\right)x_{\tau,\mathfrak{u}}$, 
as asserted. 
\end{proof}    

\begin{corollary} 
   \label{etasigmainjection} 
   The map 
   \(\namedright{(\underline{X},\underline{\ast})^{\sigma}}{\eta(\sigma)}{(\underline{X},\underline{\ast})^{\sigma,c}}\) 
   induces an injection in cohomology. 
\end{corollary} 

\begin{proof} 
By its definition, 
$\underline{c}^{\wedge\sigma/\tau}\colon\underline{X}^{\wedge\tau}\to \underline{X}^{\wedge\tau}$ 
induces in cohomology the multiplication
\[
(\underline{c}^{\wedge\sigma/\tau})^*(x)=\left(\prod_{i\in\tau}\frac{c_i^{\sigma}}{c_i^{\tau}}\right)x
\]
for any $x\in\tilde{H}^*(\underline{X}^{\wedge\tau};R)$. Since each $H^{\ast}(X_{i};R)$ for $1\leq i\leq m$ 
is assumed to be free as a graded $R$-module, so is $H^{\ast}(\underline{X}^{\wedge\tau};R)$. 
So since $R\subset \mathbb{Q}$ is an integral domain $(\underline{c}^{\wedge\sigma/\tau})^{\ast}$ is an injection. As this is true for any $\tau\subseteq\sigma$, 
the map 
\[\lllnamedright{\bigoplus_{\tau\subseteq\sigma} H^{\ast}(\underline{X}^{\tau};R)} 
     {\bigoplus_{\tau\subseteq\sigma}(\underline{c}^{\wedge\sigma/\tau})^{\ast}} 
     {\bigoplus_{\tau\subseteq\sigma} H^{\ast}(\underline{X}^{\tau};R)}\] 
is an injection. 

By Lemma~\ref{lemma_epsilon iso} and Corollary~\ref{suspXScor} respectively, the maps $\epsilon(\sigma,c)$ 
and $\epsilon(\sigma)$ are isomorphisms of graded $R$-modules. Thus the commutative 
diagram in the statement of Proposition~\ref{prop_(sigma) and eta(sigma)} implies that 
$\eta(\sigma)^{\ast}=\epsilon(\sigma)\circ(\bigoplus_{\tau\subseteq\sigma}(\underline{c}^{\wedge\sigma/\tau})^{\ast} 
     \circ(\epsilon(\sigma,c))^{-1}$, 
and therefore $\eta(\sigma)^{\ast}$ is also an injection. 
\end{proof}   

Putting all this together, we describe the cohomology ring of $(\underline{X},\underline{\ast})^{\sigma,c}$ in 
the special case when each $X_i$ is a suspension and $H^*(X_i;R)$ is a free $R$-module. Recall the additive basis  $\{y_{\tau,\mathfrak{u}}\mid\tau\subseteq\sigma, \mathfrak{u}\in J_{\tau}\}$ for $H^*((\underline{X},\underline{\ast})^{\sigma,c};R)$.

\begin{proposition}\label{cohomology ring sigma}
Let $R$ be a subring of $\Q$ and let $c$ be a  power sequence. For $1\leq i\leq m$ suppose that each $X_i$ is a suspension, the power maps on $X_{i}$ are induced by the suspension structure, and $H^*(X_i;R)$ is free as a graded $R$-module. If $\sigma\in\Delta^{m-1}$ then the cup product for $H^*((\underline{X},\underline{\ast})^{\sigma,c};R)$ is determined 
by graded commutativity, linearity and the equations: 
\[y_{i_1,j_1}\cup\cdots\cup y_{i_k,j_k}=\begin{cases}
\left(\displaystyle\prod_{\ell=1}^{k}\frac{c_{i_\ell}^{\tau}}{c_{i_\ell}^{\{i_\ell\}}}\right)y_{\tau,\mathfrak{u}}	
      &\text{if }\tau=(i_1,\ldots,i_k)\text{ for $i_{1}<\cdots<i_{k}$ and }\mathfrak{u}=\{j_1,\ldots,j_k\}\\
0	&\text{if some }i_s=i_t. 
\end{cases}\] 
\end{proposition}

\begin{proof} 
First observe that, if $\tau=(i_{1},\ldots,i_{k})\subseteq\sigma$ with $i_{1}<\cdots <i_{k}$ and 
$\mathfrak{u}=(j_{1},\ldots,j_{k})\in J_{\tau}$, then by Corollary~\ref{xycup}, 
\[
\eta(\sigma)^*(y_{\tau,\mathfrak{u}})=\left(\prod^k_{\ell=1}\frac{c_{i_{\ell}}^{\sigma}}{c_{i_{\ell}}^{\tau}}\right)x_{\tau,\mathfrak{u}}\in H^*((\underline{X},\underline{\ast})^{\sigma};R).
\]
In particular, 
\begin{equation} 
  \label{etayx} 
  \mbox{if $\tau=\{i\}$ and $\mathfrak{u}=\{j\}$ then $\eta(\sigma)^{\ast}(y_{i,j})=\displaystyle\frac{c^{\sigma}_{i}}{c^{\{i\}}_{i}}x_{i,j}$}. 
\end{equation}  
Now consider the sequence of equalities 
\[
\eta(\sigma)^*(y_{i_1,j_1}\cup\cdots\cup y_{i_k,j_k})
=\prod_{\ell=1}^{k}\eta(\sigma)^{\ast}(y_{i_{\ell},j_{\ell}}) 
     =\prod^k_{\ell=1}\left(\frac{c_{i_{\ell}}^{\sigma}}{c_{i_{\ell}}^{\{i_{\ell}\}}}x_{i_{\ell},j_{\ell}}\right) 
     =\left(\prod^k_{\ell=1}\frac{c_{i_{\ell}}^{\sigma}}{c_{i_{\ell}}^{\{i_{\ell}\}}}\right)x_{\tau,\mathfrak{u}}. 
\] 
The first equality holds since $\eta(\sigma)^*$ is a ring homomorphism, the second equality holds 
by~(\ref{etayx}), and the third holds by the cup product structure in 
$H^{\ast}((\underline{X},\underline{\ast})^{\sigma};R)$, noting that the indices 
$i_{1},\ldots,i_{k}$ are all distinct. Since $\eta(\sigma)^{\ast}$ is an injection by 
Corollary~\ref{etasigmainjection}, we obtain 
\[
y_{i_1,j_1}\cdots y_{i_k,j_k}=\left(\prod^k_{\ell=1}\frac{c_{i_{\ell}}^{\tau}}{c_{i_{\ell}}^{\{i_{\ell}\}}} \right)y_{\tau,\mathfrak{u}}.
\] 

Second, we show that $y_{i_1,j_1}\cdots y_{i_k,j_k}=0$ if some $i_s=i_t$. It suffices to show 
that $y_{i,j}\cdot y_{i,j'}=0$. Since $\eta(\sigma)^{\ast}$ is an injection by Corollary~\ref{etasigmainjection}, 
it suffices to show that $\eta(\sigma)^{\ast}(y_{i,j}\cdot y_{i,j'})=0$. By~(\ref{etayx}) and the fact that 
$\eta(\sigma)^{\ast}$ is a ring homomorphism, we obtain 
\[\eta(\sigma)^{\ast}(y_{i,j}\cdot y_{i,j'})=
    \left(\frac{c^{\sigma}_i}{c^{\{i\}}_i}\right)\cdot\left(\frac{c^{\sigma}_i}{c^{\{i\}}_i}\right) x_{i,j}\cdot x_{i,j'}.\] 
By assumption, each $X_{i}$ for $1\leq i\leq m$ is a suspension, so Remark~\ref{starremark} implies 
that $x_{i,j}\cdot x_{i,j'}=0$.  
\end{proof} 

It will be useful to abstract the properties that will describe the cup product structure in 
$H^{\ast}((\underline{X},\underline{\ast})^{K,c};R)$. 

\begin{definition}\label{def_lambda_weighted_graded}
Let $R$ be a subring of $\Q$, let $c$ be a power sequence and let $\underline{Y}=\{Y_i\}^m_{i=1}$ be a sequence of positively graded free $R$-modules. Suppose for $1\leq i\leq m$ each $Y_i$ has basis $\{\bar{y}_{i,j}|j\in J_i\}$ for a countable index set $J_i$. If $\tau=(i_1,\ldots,i_k)\in\Delta^{m-1}$ let $J_{\tau}=J_{i_1}\times\cdots\times J_{i_k}$. The \emph{weighted algebra} $\Lambda(\underline{Y},c)$ associated to $\underline{Y}$ and~$c$ is defined as the following connected graded commutative $R$-algebra. An $R$-module basis of $\Lambda(\underline{Y},c)^{>0}$ is a collection of elements  
\[
\{\bar{y}_{\tau,\mathfrak{u}}\mid\tau\in\Delta^{m-1}, \mathfrak{u}\in J_{\tau}\}
\]
where, if $\tau=(i_{1},\ldots,i_{k})$ and $\mathfrak{u}=(j_1,\ldots,j_k)$, the degree of $\bar{y}_{\tau,\mathfrak{u}}$ is $\sum^k_{\ell=1}|\bar{y}_{i_{\ell},j_{\ell}}|$. For any elements $i_{1},\ldots,i_{k}$ of $[m]$, not necessarily distinct, 
define the product
\begin{equation}\label{multiplication of Lambda Y 1}
\bar{y}_{i_1,j_1}\cdots \bar{y}_{i_k,j_k}=\begin{cases}
\left(\displaystyle\prod_{\ell=1}^{k}\frac{c_{i_\ell}^{\tau}}{c_{i_\ell}^{\{i_\ell\}}}\right)\bar{y}_{\tau,\mathfrak{u}}	
        &\text{if }\tau=(i_1<\cdots<i_k)\text{ if $i_{1}<\cdots<i_{k}$ and }\mathfrak{u}=\{j_1,\ldots,j_k\}\\
0	&\text{if some }i_s=i_t. 
\end{cases}
\end{equation}
The multiplication on the basis $\{\bar{y}_{\tau,\mathfrak{u}}\mid\tau\in\Delta^{m-1}, \mathfrak{u}\in J_{\tau}\}$ is determined by~(\ref{multiplication of Lambda Y 1}),
graded commutativity and linearity. 
\end{definition} 

\begin{definition} 
Given $\sigma\in\Delta^{m-1}$ let $\Lambda(\underline{Y},c)^{\sigma}$ be the subring of $\Lambda(\underline{Y},c)$ generated by 
\[ 
\{\bar{y}_{\tau,\mathfrak{u}}\mid\tau\subseteq\sigma, \mathfrak{u}\in J_{\tau}\}. 
\] 
\end{definition} 


We wish to show that the cohomology of a weighted polyhedral product is a weighted algebra, in the case when each underlying space $X_{i}$ is a suspension and the power maps are determined by the suspension structure. This begins with the case when $K$ is a simplex. 

\begin{proposition}\label{prop_cohomology ring sigma}
Let $R$ be a subring of $\Q$ and let $c$ be a  power sequence. For $1\leq i\leq m$ suppose that each $X_i$ is a suspension, the power maps on $X_{i}$ are induced by the suspension structure, and $H^*(X_i;R)$ is free as a graded $R$-module; let $\underline{Y}=\{\tilde{H}^*(X_i;R)\}^m_{i=1}$. If $\sigma\in\Delta^{m-1}$ then there is an isomorphism of graded commutative $R$-algebras
\[
H^*((\underline{X},\underline{\ast})^{\sigma,c};R)\cong\Lambda(\underline{Y},c)^{\sigma}.
\] 
\end{proposition}

\begin{proof} 
The graded $R$-module basis 
$\{y_{\tau,\mathfrak{u}}\mid\tau\subseteq\sigma,\mathfrak{u}\in J_{\tau}\}$ of 
$H^*((\underline{X},\underline{\ast})^{\sigma,c};R)$ is in one-to-one correspondence with the graded $R$-module basis 
$\{\bar{y}_{\tau,\mathfrak{u}}\mid\tau\subseteq\sigma,\mathfrak{u}\in J_{\tau}\}$ 
for $\Lambda(\underline{Y},c)^{\sigma}$. So to prove that there is an isomorphism 
$H^*((\underline{X},\underline{\ast})^{\sigma,c};R)\cong\Lambda(\underline{Y},c)^{\sigma}$ 
of graded commutative $R$-algebras 
we need to show that cup products for the basis elements in $H^*((\underline{X},\underline{\ast})^{\sigma,c};R)$ satisfy~(\ref{multiplication of Lambda Y 1}). But this holds by Proposition~\ref{cohomology ring sigma}.  
\end{proof}

\begin{theorem}\label{weighted_exterior_realization} 
Assume the hypothesis in Proposition~\ref{prop_cohomology ring sigma}. If $K$ is a simplicial complex on the vertex set $[m]$ then there is an isomorphism of graded commutative $R$-algebras
\[
H^*((\underline{X},\underline{\ast})^{K,c};R)\cong \Lambda(\underline{Y},c)/I_K,
\]
where $I_K$ is the ideal generated by the elements $y_{\tau,\mathfrak{u}}$ satisfying $\tau\notin K$.
\end{theorem}

\begin{proof} 
By~(\ref{diagram_epsilon nat incl}) the simplicial inclusion 
\(\namedright{K}{}{\Delta^{m-1}}\) 
induces a commutative diagram 
\begin{equation}\label{diagram_epsilon nat incl KDelta}
\xymatrix{
\underset{\tau\subseteq\Delta^{m-1}}{\bigoplus}\tilde{H}^*(\underline{X}^{\wedge\tau};R)\ar[rr]^-{\epsilon(\Delta^{m-1},c)}\ar[d]_-{\pi}	& &\tilde{H}^*((\underline{X},\underline{\ast})^{\Delta^{m-1},c};R)\ar[d]^-{\imath^{\ast}}\\
\underset{\tau\subseteq K}{\bigoplus}\tilde{H}^*(\underline{X}^{\wedge\tau};R)\ar[rr]^-{\epsilon(K,c)}	& &\tilde{H}^*((\underline{X},\underline{\ast})^{K,c};R) 
}
\end{equation}
where $\pi$ is the projection and $\imath$ is shorthand for $\imath^{\Delta^{m-1}}_{K}$. 
By Lemma~\ref{lemma_epsilon iso}, both $\epsilon(\Delta^{m-1},c)$ and $\epsilon(K,c)$ are isomorphisms of 
graded $R$-modules. Thus the kernel of $\iota^{\ast}$ is isomorphic to the kernel of $\pi$, as 
graded $R$-modules. Further, since $\pi$ has a right inverse, its kernel is a direct summand of 
$\underset{\tau\subseteq\Delta^{m-1}}{\bigoplus}\tilde{H}^*(\underline{X}^{\wedge\tau};R)$, 
implying that the kernel of $\iota^{\ast}$ is a direct summand of 
$\tilde{H}^*((\underline{X},\underline{\ast})^{\Delta^{m-1},c};R)$. Thus there is a split short exact sequence 
of graded $R$-modules 
\[
0\longrightarrow ker(\imath^{\ast})\overset{q}{\longrightarrow}\tilde{H}^*((\underline{X},\underline{\ast})^{\Delta^{m-1},c};R)\overset{\imath^*}{\longrightarrow}\tilde{H}^*((\underline{X},\underline{\ast})^{K,c};R)\longrightarrow0.
\]
Since $\imath^*$ is a multiplicative epimorphism, we obtain an isomorphism of graded commutative $R$-algebras 
\[
H^*((\underline{X},\underline{\ast})^{K,c};R)\cong H^*((\underline{X},\underline{\ast})^{\Delta^{m-1},c};R)/ker(\imath^*). 
\] 
By Proposition~\ref{prop_cohomology ring sigma}, there is an isomorphism of graded commutative $R$-algebras 
$H^*((\underline{X},\underline{\ast})^{\Delta^{m-1},c};R)\cong\Lambda(\underline{Y},c)$, implying that there 
is an isomorphism of graded commutative $R$-algebras 
\[
H^*((\underline{X},\underline{\ast})^{K,c};R)\cong\Lambda(\underline{Y},c)/ker(\imath^*).
\] 
Since $ker(\imath^{\ast})=ker(\pi)$, and $ker(\pi)$ is generated by those elements $y_{\tau,\mathfrak{u}}$ 
such that $\tau\notin K$, that is $I_K$, we obtain an isomorphism of graded commutative $R$-algebras 
$H^*((\underline{X},\underline{\ast})^{K,c};R)\cong \Lambda(\underline{Y},c)/I_K$.
\end{proof}

\section{Spaces with weighted sphere product cohomology} 

In this section we return to Steenrod's problem and prove Theorem~\ref{geomrealization}. This begins 
by defining a family of weighted sphere product algebras whose product structure is governed by 
a certain ``coefficient sequence", then showing that a power sequence produces a 
coefficient sequence and the corresponding algebra can be realized as 
the cohomology of a weighted polyhedral product. We then go on to compare the family 
of coefficient sequences to the family of power sequences. 

Throughout this section, let $R$ be a subring of $\mathbb{Q}$. 

\begin{definition} \label{coeff_sequ}
\label{coeffseq} 
A {\em coefficient sequence} (of length $m$) is a map
$\mathfrak{c}\colon \Delta^{m-1}\rightarrow\mathbb{N}$, $\sigma\mapsto\mathfrak{c}_{\sigma}$ satisfying:  
\begin{itemize} 
   \item [1)] $\mathfrak{c}_{\emptyset}=\mathfrak{c}_{\{i\}}=1$;  
   \item [2)]  if $\sigma=\sigma'\cup \sigma''$ and 
                   $\sigma'\cap \sigma''=\emptyset$ then 
                   $\mathfrak{c}_{\sigma'}\mathfrak{c}_{\sigma''}\mid\mathfrak{c}_{\sigma}$;
   \item [3)] $(\mathfrak{c}_{\sigma},p)=1$ if $p$ is invertible in $R$. 
\end{itemize} 
\end{definition} 

\begin{definition} 
\label{Acdef} 
Let $\mathfrak{c}$ be a coefficient sequence and associate to each $i\in [m]$ a degree $d_i$ in $\mathbb{N}$. 
A \emph{weighted sphere product algebra} is a graded commutative algebra $A(\mathfrak{c})$ defined as follows. 
As a graded module, let 
\[A(\mathfrak{c})\cong \oplus_{\sigma\subseteq [m]}R\langle a_{\sigma}\rangle\] 
where $|a_{\sigma}|=\Sigma_{i\in \sigma} d_i$. Let the product be determined by the formula 
\[\prod_{i\in \sigma}a_{\{i\}}=\mathfrak{c}_{\sigma} a_{\sigma}\] 
where $\sigma\in\Delta^{m-1}$. 
\end{definition} 

In terms of the weighted sphere product algebra $A(\mathfrak{c})$, property 1) defining a coefficient sequence 
can be thought of as a normalizing condition and property 2) is a compatibility condition between 
subsets that is necessary for the multiplication on $A(\mathfrak{c})$ to be well-defined. Property 3) is 
odd looking on first glance but is a streamlining condition that is justified by the following lemma. 

\begin{lemma}\label{all_coeff}
Let $\mathfrak{c}\colon\namedright{\Delta^{m-1}}{}{\mathbb{N}}$ be a sequence satisfying 1) and 2) in 
Definition~\ref{coeffseq}. If we let $A(\mathfrak{c})$ be defined as in Definition~\ref{Acdef} then there is a 
coefficient sequence $\mathfrak{c}'$ such that $A(\mathfrak{c})\cong A(\mathfrak{c}')$. 
\end{lemma}

\begin{proof}
Given $\sigma\in\Delta^{m-1}$, write $\mathfrak{c}_{\sigma}=u_{\sigma}v_{\sigma}$ where $u_{\sigma}$ is 
invertible in $R$ and $(u_{\sigma},v_{\sigma})=1$. Let $\mathfrak{c}'_{\sigma}=v_{\sigma}$. Then the 
sequence $\mathfrak{c}'$ is a coefficient sequence. As graded modules, 
$A(\mathfrak{c})$ and $A(\mathfrak{c}')$ are identical. Define a ring homomorphism 
\(\namedright{A(\mathfrak{c})}{}{A(\mathfrak{c}')}\) 
by sending $a_{\sigma}$ to $\frac{1}{u_{\sigma}}a_{\sigma}$ and extending multiplicatively. 
This clearly has an inverse 
\(\namedright{A(\mathfrak{c}')}{}{A(\mathfrak{c})}\) 
given by sending $a_{\sigma}$ to $u_{\sigma}a_{\sigma}$. 
\end{proof} 

Steenrod's problem in the case of a weighted sphere product algebra $A(\mathfrak{c})$ asks when 
this algebra can be realized as the cohomology of a space. Phrased slightly differently, 
for which coefficient sequences $\mathfrak{c}$ can $A(\mathfrak{c})$ be realized as the cohomology of a space? 
We can also ask a slightly more general question. Suppose $A$ is a torsion free finitely generated graded commutative algebra such that 
$A\otimes \mathbb{Q}$ is a rational sphere product algebra, or alternatively suppose $A$ is an order in a 
rational sphere product algebra. Is $A$ realizable?
Some cases of this are discussed in \cite{SSTW}.

To address this we first relate power sequences to coefficient sequences. Let $PS$ denote 
the set of power sequences (of length $m$) with the property that $c\in PS$ satisfies $c_{i}^{\{i\}}=1$ 
for $1\leq i\leq m$, and let $CS$ denote the set of coefficient sequences 
(of length $m$). Define a map
\begin{equation}\label{eqn_Phi map}
\Phi\colon \mbox{PS} \rightarrow \mbox{CS}
\end{equation}
by sending a power sequence $c$ to the coefficient sequence $\Phi(c)$ defined by 
$\Phi(c)_{\sigma}=\prod_{i\in \sigma} c_i^{\sigma}$. Note that if $\sigma=\{i\}$ then 
$\Phi(c)_{\{i\}}=c_{i}^{\{i\}}$, so the requirement that $c_{i}^{\{i\}}=1$ implies that 
$\Phi(c)_{\{1\}}$ satisfies condition 1) of Definition~\ref{coeffseq}. 

\begin{theorem} 
\label{Steenrodpowerseq} 
Let $R$ be a subring of $\mathbb{Q}$. Suppose that $c$ is a power sequence of length $m$ and 
$\{d_{1},\ldots,d_{m}\}$ is a set of degrees. If $X_i=S^{d_i}$ for $1\leq i\leq m$ then there is an 
isomorphism of graded commutative $R$-algebras
\[
H^*((\underline{X},\underline{\ast})^{\Delta^{m-1},c};R)\cong A(\Phi(c)). 
\] 
\end{theorem}

\begin{proof} 
Applying Theorem~\ref{weighted_exterior_realization} to the power sequence $c$, and noting 
that with $K=\Delta^{m-1}$ the ideal~$I_{K}$ is empty since every possible face is in $K$, 
there is an isomorphism of graded commutative $R$-algebras 
\[
H^*((\underline{X},\underline{\ast})^{\Delta^{m-1},c};R)\cong \Lambda(\underline{Y},c)
\] 
where $\underline{Y}=\{\tilde{H}^*(S^{d_{i}};R)\}_{i=1}^{m}$. This description of $\underline{Y}$ 
implies that the definitions of $\Lambda(\underline{Y},c)$ and $A(\Phi(c))$ coincide. 
\end{proof} 

Theorem~\ref{Steenrodpowerseq} implies that any coefficient sequence that is the image 
of a power sequence has the property that the corresponding weighted sphere product algebra  
can be geometrically realized as the cohomology of a space. This 
leads to the question of whether $\Phi$ is an isomorphism. We will show that this is not 
true, $\Phi$ is neither injective nor surjective. Therefore there are other cases of coefficient 
sequences for which Steenrod's problem remains open, and addressing these cases will 
require different techniques. When $m=3$ additional methods have been used to show all weighted sphere product algebras can be realized \cite{SSTW}.
\medskip 

\noindent 
\textbf{Comparing power and coefficient sequences}. 
To close the paper we compare power sequences and coefficient sequences, showing 
that $\Phi$ is neither injective nor surjective, while showing that power sequences and 
coefficient sequences are not too dissimilar since they both generate monoids with the 
same group completion. 

\begin{lemma} 
   The map 
   \(\Phi\colon\namedright{PS}{}{CS}\) 
   is not an injection. 
\end{lemma} 

\begin{proof} 
Take $m=3$. Define a power sequence $c$ of size $3$ by the following data: by 
definition of a power sequence, $c_{i}^{\{i\}}=1$ for $i\in\{1,2,3\}$, and let 
$$
\begin{matrix}
c_1^{( 1,2)}=p & c_2^{(1,2)}=1 \\

c_1^{( 1,3)}=p & c_3^{( 1,3)}=1 \\

c_2^{(2,3)}=p & c_3^{(2,3)}=1 \\

c_1^{(1,2,3)}=p & c_2^{(1,2,3)}=p & c_3^{( 1,2,3)}=1
\end{matrix}
$$
Then $\Phi(c)$ satisfies 
$$
\Phi(c)_{\sigma}=
\begin{cases}
1 & |\sigma|=1 \\
p & |\sigma|=2 \\
p^2 & |\sigma|=3.
\end{cases}
$$
In particular the coefficients of $\Phi(c)$ are invariant under the action of the symmetric group. 
However, $c$ is not invariant under the action of the symmetric group: define the power sequence 
$\bar{c}$ by $\bar{c}^{\sigma}_{i}=c^{\sigma}_{i}$ for all $\sigma\in\Delta^{m-1}$ except $\sigma=(1,2)$, 
where $\bar{c}_{1}^{(1,2)}=1$ and $\bar{c}_{2}^{(1,2)}=p$. Then $\bar{c}\neq c$ but $\Phi(\bar{c})=\Phi(c)$. 
\end{proof}

\begin{lemma}\label{no_surj}
   The map 
   \(\Phi\colon\namedright{PS}{}{CS}\) 
   is not a surjection.
\end{lemma}
 
\begin{proof} 
Let $\mathfrak{c}$ be the coefficient sequence of size $3$ defined by 
$$
\mathfrak{c}_{\sigma}=
\begin{cases}
2 & |\sigma|=2 \mbox{ or } 3 \\
1 & |\sigma|=1.
\end{cases}
$$
Suppose that $c\in PS$ satisfies $\Phi(c)=\mathfrak{c}$. By definition, 
$\Phi(c)_{(1,2)}=c_{1}^{(1,2)} c_{2}^{(1,2)}$, so $\Phi(c)_{(1,2)}=2$ implies that 
either $c_1^{(1,2)}=2$ or $c_2^{(1,2)}=2$. Similarly, either 
$c_2^{(2,3)}=2$ or $c_3^{(2,3)}=2$. The definition of a power sequence then implies that 
at least two of $c_1^{(1,2,3)}$, $c_2^{(1,2,3)}$ and $c_3^{(1,2,3)}$ are divisible 
by~$2$ and hence
$\Phi(c)_{(1,2,3)}=c_{1}^{(1,2,3)}c_{2}^{(1,2,3)}c_{3}^{(1,2,3)}$ implies 
that $\Phi(c)_{(1,2,3)}$ is divisible by $4$. But $\Phi(c)_{(1,2,3)}=\mathfrak{c}_{(1,2,3)}=2$, a 
contradiction. 
\end{proof}

\begin{remark}\label{use2}
In the size $2$ case (using vertex set $\{1, 2\}$), $A(\mathfrak{c})$ is determined by $\mathfrak{c}_{12}$.
As a $\Z$-module $A(\mathfrak{c})=\Z\langle 1, a_1, a_2, a_{12}\rangle$ and the product is determined by 
$a_1a_2=\mathfrak{c}_{12}a_{12}$. In that case any power sequence 
with $c^{1}=c^{2}=(1,1)$ and $c^{(1,2)}=(\alpha, \beta)$ with $\alpha\beta=c_{12}$ will give $\Phi(c)=\mathfrak{c}$ 
and so $\Phi$ is surjective when restricted to this case, and all the size $2$ $A(\mathfrak{c})$ can be realized by Theorem
\ref{Steenrodpowerseq}. This can also be seen directly using Whitehead products. 
\end{remark}

Next, we consider how $PS$ and $CS$ are similar. 
For a prime $p$ let $PS(p)$ denote the subset of $PS$ such that for all $\sigma\in\Delta^{m-1}$ 
and $i\in \sigma$, $c_i^{\sigma}=p^k$ for some $k\in \mathbb{N}$. Similarly define $CS(p)$. Let 
\[\Phi(p)\colon \mbox{PS(p)} \rightarrow \mbox{CS(p)}\] 
be the restriction of $\Phi$ to $PS(p)$. Define a pointwise multiplication on $PS$ as  
follows. If $c,d\in PS$ let $cd\in PS$ be the power sequence defined by 
$$ (cd)^{\sigma}_i=c^{\sigma}_id^{\sigma}_i.$$
Similarly define a pointwise multiplication on $CS$: if $\mathfrak{c},\mathfrak{d}\in CS$ let 
$\mathfrak{c}\mathfrak{d}\in CS$ be the coefficient sequence defined by 
$$(\mathfrak{c}\mathfrak{d})_{\sigma}=\mathfrak{c}_{\sigma}\mathfrak{d}_{\sigma}.$$ 
Observe that these multiplications restrict to multiplications on $PS(p)$ and $CS(p)$. 

\begin{proposition}
With the pointwise multiplication, $PS$, $CS$, $PS(p)$ and $CS(p)$ are all torsion free 
commutative monoids and the maps $\Phi\colon PS\rightarrow CS$ and 
$\Phi(p)\colon PS(p)\rightarrow CS(p)$
are all maps of monoids. 
\end{proposition}

\begin{proof}
This follows immediately from $\mathbb{N}$ being a torsion free commutative monoid and the 
definition of the pointwise multiplication on $PS$ and $CS$.  
\end{proof}

We specify some distinguished power and coefficient sequences that will generate 
the group completions of $PS(p)$ and $CS(p)$ but not the monoids themselves. 
Fix $\tau\in\Delta^{m-1}$ and $j\in [m]$. Let $c(\tau, j)\in PS(p)$ be defined by
\begin{equation} 
\label{ctauj} 
c(\tau,j)_i^{\sigma}=
\begin{cases}
p & j=i \mbox{ and } \tau\subseteq\sigma \\
1 & \mbox{else}
\end{cases}
\end{equation}
and let $\mathfrak{c}(\tau)\in CS(p)$ be defined by 
\begin{equation}\label{ctau_CS}
\mathfrak{c}(\tau)_{\sigma}=
\begin{cases}
p & \tau\subseteq\sigma \\
1 & \mbox{else}.
\end{cases}
\end{equation}
Similarly, let $d(\tau, j)\in PS(p)$ be defined by
\begin{equation}
d(\tau,j)_i^{\sigma}=
\begin{cases}
p & j=i \mbox{ and } \tau=\sigma \\
1 & \mbox{else}
\end{cases}
\end{equation}
and let $\mathfrak{d}(\tau)\in CS(p)$ be defined by 
\begin{equation}
\mathfrak{d}(\tau)_{\sigma}=
\begin{cases}
p & \tau=\sigma \\
1 & \mbox{else}.
\end{cases}
\end{equation} 

\begin{lemma} 
\label{Phicd} 
For $\tau\in\Delta^{m-1}$ and $j\in [m]$, we have $\Phi(c(\tau,j))=\mathfrak{c}(\tau)$ 
and $\Phi(d(\tau,j))=\mathfrak{d}(\tau)$. 
\end{lemma} 

\begin{proof} 
Let $\sigma\in\Delta^{m-1}$. By definition,  
$\Phi(c(\tau,j))_{\sigma}=\prod_{i\in\sigma} c(\tau,j)_{i}^{\sigma}$. There are two cases.
First, if $\tau\not\subseteq\sigma$ then by~(\ref{ctauj}) each $c(\tau,j)_{i}^{\sigma}=1$, so 
$\Phi(c(\tau,j))_{\sigma}=1$. Second, if $\tau\subseteq\sigma$ then by~(\ref{ctauj}) we have  
$c(\tau,j)_{i}^{\sigma}=p$ for the one instance when $i=j$ and $c(\tau,j)_{i}^{\sigma}=1$ 
for all $j\neq i$. Thus $\Phi(c(\tau,j))_{\sigma}=p$. In either case we obtain 
$\Phi(c(\tau,j))_{\sigma}=\mathfrak{c}(\tau)_{\sigma}$. As this is true for all $\sigma\in\Delta^{m-1}$ 
we obtain $\Phi(c(\tau,j))=\mathfrak{c}(\tau)$. Similarly, $\Phi(d(\tau,j))=\mathfrak{d}(\tau)$. 
\end{proof} 

\begin{proposition} 
\label{cgeneratesd} 
For $\tau\in\Delta^{m-1}$ and $j\in [m]$, we have
$d(\tau, j)=\prod_{\tau\subseteq \sigma} c(\sigma,j)^{(-1)^{|\sigma\setminus \tau|}}$ and 
$\mathfrak{d}(\tau)=\prod_{\tau\subseteq \sigma} \mathfrak{c}(\sigma)^{(-1)^{|\sigma\setminus \tau|}}$.
\end{proposition}

\begin{proof} 
Since each $d(\tau,j)$ is defined in terms of its components $d(\tau,j)^{\delta}_{i}$, it 
suffices to show that 
\begin{equation} 
  \label{ddelta} 
  d(\tau, j)^{\delta}_{i}=\prod_{\tau\subseteq\sigma} (c(\sigma,j)^{\delta}_{i})^{(-1)^{|\sigma\setminus \tau|}} 
\end{equation}  
for all $\delta\in\Delta^{m-1}$ and $i\in [m]$. 

If $j\neq i$ then, by definition, $d(\tau,j)^{\delta}_{i}=1$ and each $c(\sigma,j)^{\delta}_{i}=1$ 
so~(\ref{ddelta}) holds. If $j=i$, observe that, by its definition, $c(\sigma,j)^{\delta}_{j}=1$ if 
$\sigma\not\subset\delta$. Therefore 
\begin{equation} 
  \label{ddelta2} 
  d(\tau,j)^{\delta}_{j}=\prod_{\tau\subseteq \sigma} (c(\sigma,j)^{\delta}_{j})^{(-1)^{|\sigma\setminus \tau|}}= 
       \prod_{\tau\subseteq \sigma\subseteq\delta} (c(\sigma,j)^{\delta}_{j})^{(-1)^{|\sigma\setminus \tau|}}. 
\end{equation}  
\textit{Case 1}: If $\delta=\tau$ then, by its definition, $d(\tau,j)^{\tau}_{j}=p$ while the right side of~(\ref{ddelta2}) 
has the product indexed over a single factor since $\tau=\sigma=\delta$ and that factor 
is $(c(\tau,j)^{\tau}_{j})^{(-1)^{|\tau\setminus \tau|}}=c(\tau,j)^{\tau}_{j}$, which by definition is~$p$. 
Hence~(\ref{ddelta2}) holds in this case. 

\noindent 
\textit{Case 2}: If $\delta\neq\tau$ then, as $\tau$ is a subset of $\delta$, 
it is a proper subset. Fix $k\in\delta\backslash\tau$. Observe that there is a disjoint union
\[\{\sigma\mid\tau\subseteq\sigma\subseteq\delta\}= 
     \{\sigma\mid\tau\subseteq\sigma\subseteq\delta\ \mbox{and}\ k\in\sigma\}\cup 
     \{\sigma\mid\tau\subseteq\sigma\subseteq\delta\ \mbox{and}\ k\notin\sigma\}.\]
Since $k\in\delta\backslash\tau$, any $\sigma$ with $\tau\subseteq\sigma\subseteq\delta$ 
and $k\in\sigma$ may be written as $\sigma'\cup\{k\}$ for $\tau\subseteq\sigma'\subseteq\delta$ 
and $k\notin\sigma'$. Thus the disjoint union above may be rewritten as 
\[\{\sigma\mid\tau\subseteq\sigma\subseteq\delta\}= 
     \{\sigma\cup\{k\}\mid\tau\subseteq\sigma\subseteq\delta\ \mbox{and}\ k\notin\sigma\}\cup 
     \{\sigma\mid\tau\subseteq\sigma\subseteq\delta\ \mbox{and}\ k\notin\sigma\}.\]
In particular, the sets $A=\{\sigma\cup\{k\}\mid\tau\subseteq\sigma\subseteq\delta\ \mbox{and}\ k\notin\sigma\}$ 
and $B=\{\sigma\mid\tau\subseteq\sigma\subseteq\delta\ \mbox{and}\ k\notin\sigma\}$ have the 
same cardinality. If $\sigma\in A$ then $\sigma\cup\{k\}\in B$ and the signs of $(-1)^{|\sigma/\tau|}$ 
and $(-1)^{|\sigma\cup\{k\}/\tau|}$ are opposite. Further, since $\sigma\subseteq\delta$ and 
$\sigma\cup\{k\}\subseteq\delta$, we have $c(\sigma,j)^{\delta}_{j}=p=c(\sigma\cup\{k\},j)^{\delta}_{j}$. 
Therefore $(c(\sigma,j)^{\delta}_{j})^{(-1)^{|\sigma/\tau|}}$ is the inverse of  
$(c(\sigma\cup\{k\},j)^{\delta}_{j})^{(-1)^{|\sigma\cup\{k\}/\tau}|}$. Thus in~(\ref{ddelta2}), the 
right side of the equation equals $1$. The left side is also equal to $1$ by the definition of 
$d(\tau,j)^{\delta}_{j}$. Hence~(\ref{ddelta2}) holds in this case as well. 

Consequently, ~(\ref{ddelta}) holds in all cases. 

Finally, the second equation for $\mathfrak{d}(\tau)$ follows from applying $\Phi$ to the first equation and 
using Lemma~\ref{Phicd}. 
\end{proof} 

Combining Lemma~\ref{Phicd} and Proposition~\ref{cgeneratesd} implies the following. 

\begin{corollary}\label{completion}
$PS(p)$ is a sub-monoid of the free monoid generated by the $d(\tau,j)$ and 
its group completion is isomorphic to the free 
abelian group $\mathbb{Z}\langle d(\tau,j) \rangle$. 
Similarly, $CS(p)$ is a sub-monoid of the free monoid generated by the 
$\mathfrak{d}(\tau)$ and its group completion is isomorphic to the free 
abelian group $\mathbb{Z}\langle \mathfrak{d}(\tau) \rangle$. 
Consequently, the image of $PS(p)$ generates the group completion of $CS(p)$.~$\qqed$ 
\end{corollary}

\bibliographystyle{amsalpha}

\end{document}